\documentclass[a4paper,11pt]{article}
\pdfoutput=1 % if your are submitting a pdflatex (i.e. if you have images in pdf, png or jpg format)

\usepackage[paper=a4paper,left=3cm,right=3cm,top=3cm,bottom=3cm,footskip=1cm]{geometry}

\usepackage{graphicx}
\usepackage{graphics}
\usepackage{bbm}
\usepackage{amsmath,amssymb}
\usepackage{setspace}
\usepackage[T1]{fontenc}
\usepackage{bbm}

\usepackage[export]{adjustbox}
\usepackage{subfig}
\usepackage{xspace}

\newcommand{\bR}{\mathbb{R}}

\newcommand{\C}{{\mathbb C}}

\newcommand{\cP}{\mathcal{P}}

\newcommand{\bea}{\begin{eqnarray}}
\newcommand{\eea}{\end{eqnarray}}
\newcommand{\be}{\begin{equation}}
\newcommand{\ee}{\end{equation}}

\newcommand{\randommarkovmove}{\textsc{RandomMarkovMove}\xspace}
\newcommand{\smartcollapse}{\textsc{SmartCollapse}\xspace}
\newcommand{\knotify}{\textsc{Knotify}\xspace}
\newcommand{\randomunknot}{\textsc{RandomUnknot}\xspace}
\newcommand{\randomknot}{\textsc{RandomKnot}\xspace}
\newcommand{\braidrelationone}{\textsc{BraidRelation}1\xspace}

\usepackage{algpseudocode}
\usepackage{algorithm}

\usepackage{hyperref}

\begin{document}
\begin{titlepage}
\begin{flushright}
CALT-2020-046\\
CERN-TH-2020-179
\end{flushright}

\medskip

\begin{center}
{\LARGE {\bf Learning to Unknot}}\\[12pt]
\bigskip

Sergei Gukov$^1$, James Halverson$^{2,3}$, Fabian Ruehle$^{4,5}$, Piotr Su\l{}kowski$^{1,6}$
\bigskip
\setstretch{1.0}

\vspace{0.cm}
{
	{\it ${}^{\text{1}}$ Walter Burke Institute for Theoretical Physics,\\California Institute of Technology, Pasadena, CA 91125, USA}\\[.5em]
	{\it ${}^{\text{2}}$ Department of Physics, Northeastern University, Boston, MA 02115}\\[.5em]
	{\it ${}^{\text{3}}$ The NSF AI Institute for Artificial Intelligence and Fundamental Interactions}\\[.5em]
	{\it ${}^{\text{4}}$ CERN Theory Department, 1 Esplanade des Particules,\\ CH-1211 Geneva, Switzerland}\\[.5em]
	{\it ${}^{\text{5}}$ Rudolf Peierls Centre for Theoretical Physics, Department of Physics, \\ University of Oxford, Parks Road, Oxford OX1 3PU, United Kingdom}\\[.5em]
	{\it ${}^{\text{6}}$ Faculty of Physics, University of Warsaw, ul. Pasteura 5, 02-093 Warsaw, Poland}\\[.5em]
}
\medskip
Emails: \href{mailto:gukov@theory.caltech.edu}{gukov@theory.caltech.edu},~\href{mailto:j.halverson@northeastern.edu}{j.halverson@northeastern.edu},\\ \href{mailto:fabian.ruehle@cern.ch}{fabian.ruehle@cern.ch}, \href{mailto:psulkows@fuw.edu.pl}{psulkows@fuw.edu.pl}
\end{center}
\setstretch{1.3}
\bigskip\bigskip

\begin{abstract}\noindent
We introduce natural language processing into the study of knot theory, as made natural by the braid word representation of knots. We study the UNKNOT problem of determining whether or not a given knot is the unknot. After describing an algorithm to randomly generate $N$-crossing braids and their knot closures and discussing the induced prior on the distribution of knots, we apply binary classification to the UNKNOT decision problem. We find that the Reformer and shared-QK Transformer network architectures outperform fully-connected networks, though all perform well. Perhaps surprisingly, we find that accuracy increases with the length of the braid word, and that the networks learn a direct correlation between the confidence of their predictions and the degree of the Jones polynomial.  Finally, we utilize reinforcement learning (RL) to find sequences of Markov moves and braid relations that simplify knots and can identify unknots by explicitly giving the sequence of unknotting actions. Trust region policy optimization (TRPO) performs consistently well for a wide range of crossing numbers and thoroughly outperformed other RL algorithms and random walkers. Studying these actions, we find that braid relations are more useful in simplifying to the unknot than one of the Markov moves.
\end{abstract}
\end{titlepage}
\clearpage
\setcounter{footnote}{0}
\setcounter{tocdepth}{2}
\tableofcontents
\clearpage

\section{Introduction\label{sec:intro}}
In work and in play, some of the most difficult or even unsolvable problems can be formulated by using a fairly small set of rules. Indeed, even when the rules of the game are simple, the state space of all possible configurations can be extremely large, way too large for a human brain or a deterministic algorithm to identify a given configuration and tell where in a big scheme of things it belongs. This is precisely the domain where machine learning and artificial intelligence hold a consistent record of winning the game, growing stronger each year and outperforming the best chess grand masters~\cite{Silver2017:aab} and go players~\cite{Silver:2016aaa, Silver2017:aaa}.

There are many such ``games'' in fundamental science too, with simple rules and a vast landscape of possible outcomes. 

The one considered in this paper involves three Reidemeister moves (or, equivalently, Markov moves) as ``rules of the game'' and the rich state space is spanned by many different knots or, more precisely, by different presentations of knots. Although these basic rules can be counted on one hand and encode all possible equivalences, the richness of the state space immediately gets in the way of identifying whether two different presentations are equivalent or not. It is rather ironic that this is an obstacle to several fundamental problems in low-dimensional topology, including the smooth 4-dimensional Poincar\'e conjecture. Other areas where finding the simplest representation of a knot will be beneficial are for example the knots-quivers correspondence~\cite{Kucharski:2017poe,Kucharski:2017ogk} in physics, or protein folding in biology \cite{10.1093/nar/gku1059}.

In the field of string theory, it has been realized~\cite{He:2017aed,Krefl:2017yox,Ruehle:2017mzq,Carifio:2017bov} within the last three years that machine learning can also be applied to the large state space of string vacua and compactification spaces; see~\cite{Ruehle:2020jrk} for an introduction and overview. In particular, in~\cite{Halverson:2019tkf}, a Reinforcement learning was applied to find solutions to a set of coupled quartic Diophantine equations that describe consistent string vacua, of which there are many more~\cite{Douglas:2003um,Ashok:2003gk,Taylor:2015xtz,Halverson:2017ffz} than configurations in Go.

From the AI/ML point of view, the problem of identifying equivalence classes, i.e.\ different presentations of the same knot, is very similar to the problem of completing the sentence ``I grew up in France\ldots I speak fluent \ldots''. Roughly, the reason is that the latter task requires identifying the meaning of the sentence and placing it next to other sentences with a similar meaning in a large space of possibilities. This is a classical problem in Natural Language Understanding (NLU) or Natural Language Processing (NLP). Therefore, the question we wish to ask here is: How quickly and how well can a neural network learn to speak the language of knots?

This question was asked before, however, not from the NLP perspective, which is one novelty of this paper. For example in~\cite{Hughes:2016aaa}, Hughes uses a simple feedforward neural network to predict knot invariants such as quasi-positivity, the slice genus, and the Ozsv\'ath-Szab\'o $\tau$-invariant. In~\cite{Jejjala:2019kio} the authors also use a simple feedforward network to compute the hyperbolic knot volume from the Jones polynomial.

The knot theory problem we are studying is the UNKNOT problem, i.e.\ recognition of whether a given knot is the unknot. In addition to using NLP tools for the binary classification task, we also employ reinforcement learning to explicitly find a sequence of moves that allow to transform a (potentially complicated) representation of the unknot to its simplest representative, a circle with no crossings. Since the algorithm finds the necessary Reidemeister moves, rather than just predicting a probability for the knot being the unknot, the results can serve to prove that a given knot is the unknot.

Another novelty is that, for the NLP itself, the example of the ``knot language problem'' studied here presents new twists and opportunities. For example, the role of equivalence classes so central to this example could be also useful in other problems, not only in fundamental science.

This paper is organized as follows. In Section \ref{sec:review} we review the basics of NLP and knot theory and introduce how the braid representation of knots yields an NLP description of knots. In Section \ref{sec:generating_data} we introduce an algorithm by which trivial and non-trivial knots may be generated, represented by braids with a fixed number of crossings. In Section \ref{sec:decision} we utilize a variety of neural networks to apply binary classification to the UNKNOT problem, and use the trained networks to study correlations with the Jones polynomial and notions of hardness. In Section \ref{sec:unkotting} we utilized reinforcement learning to find sequences of Reidmeister moves, represented by braid relations and Markov moves on the braid, that simplify a non-trivial representation of the unknot to the trivial one. In Section \ref{sec:conclusion} we summarize the main results of this work and discuss. In Appendix \ref{sec:app_algos} we provide pseudo-code for some algorithms used in this paper and in Appendix \ref{sec:game} we provide an unknotting game.

\section{Knots and Natural Language\label{sec:review}}

In this section we review NLP and introduce
its application to knot theory.

\subsection{Embedding Layers for Semantics}
\label{sec:Embedding-Layers}

A language $L$ is composed of words from a vocabulary $V(L)$. 
In NLP it is useful to have an embedding of a word 
into a vector space that ideally encodes its meaning:
\begin{equation}   
E: V(L) \to \bR^{d},
\end{equation}
where $d$ is the embedding dimension. 

Since the vocabulary 
is a discrete set of words, one embedding, known as the one-hot encoding,
maps the $i^\text{th}$ word $w_i \in V(L)$ as $w_i \mapsto e_i$, where $e_i$
is a unit vector and $d=|V(L)|$. 
From the NLP perspective, this embedding has 
a number of issues. First, the dimension of the target
vector space is $|V(L)|$, which for any non-trivial
language will be quite large. Second,
all but one of the entries is zero; the vector is sparse. 
Finally, the embedding
only contains the information of the index in the
set $V(L)$, which is arbitrary and can be permuted;
no useful information is encoded in the embedding.

One would like a better technique for associating a vector to a word.
The problem of sparseness may be solved by choosing $d<|V(L)|$,
typically $d\ll |V(L)|$. In some cases $E$ is fixed by
using pre-trained word vectors for the embedding, while
in others $E$ has randomly initialized parameters and
a useful embedding is learned by training on some task.
In the process, semantics may be learned that encoded
meaning into the vector representatives of words. 
(e.g. \cite{Vylomova:2015aaa}) A famous example is 
\begin{equation}
E(\text{king}) - E(\text{man}) + E(\text{woman}) \simeq E(\text{queen}),
\end{equation}
an approximate equivalence at the level of the
vector relationships that encodes an actual semantic
relationship in the language. Other semantic
relations have also been learned, e.g.
related to capitals
\begin{equation}
E(\text{Paris}) - E(\text{France}) + E(\text{Poland}) \simeq E(\text{Warsaw}),
\end{equation}
and pluralization
\begin{equation}
    E(\text{cars}) - E(\text{car}) + E(\text{apple}) \simeq E(\text{apples}).
\end{equation}
Clearly, word embeddings that capture semantic features of 
a word or language could be useful in a variety of machine
learning tasks with respect to that language.
    
In what follows we will be discussing queries and keys, and it will be assumed
that each word in a sequence of length $l$ has been mapped to $d$-vector
via an embedding layer, so that each embedded sequence has shape $[l,d]$.

\subsection{Attention and Transformers}
Recent years have seen great progress in NLP
with the evolution
of the attention mechanism and its introduction into
various architectures. It works as the
name suggests: by training the neural network
to pay attention to the most important parts of 
sentences.

To explain the mechanism we will utilize the 
notion of queries, keys, and values 
\cite{Vaswani2017AttentionIA}. This notion is
used because the mechanism mimics the retrieval
of a value $v_q$ for a query $q$ based on
a key $k_i$ in a database, each of which
has its own value $v_i$. In normal database
retrieval, one finds the key $k_i$ that is 
identical to the query and returns the value. 
In attention, we wish instead to have a similarity
measure $s(q,k_i)$ between the query and key,
which is used as the weight to determine the attention paid to the different elements in a weighted
sum of values,
\begin{equation}
\text{Attention}(q,k,v)=v_q=\sum_i s(q,k_i) v_i.
\end{equation}
In this formulation, the case of normal database retrieval is the
case where $s(q,k_i) = 1$ if $q=k_i$ and $0$ 
otherwise.
The different types of attention that exist
in the literature \cite{Graves2014NeuralTM,Bahdanau2014NeuralMT,Luong2015EffectiveAT,Vaswani2017AttentionIA}
correspond to different
choices for similarity function $s$, which
is chosen to be differentiable (unlike usual
database retrieval) to allow for backpropagation in
a neural network. The similarity is usually softmax
applied to some score function, so that the
weights sum to one. 

The attention mechanism is a crucial component
of the so-called Transformer architecture \cite{Vaswani2017AttentionIA},
where the version of attention used is known
as scaled dot-product attention,
\begin{equation}
    \text{Attention}(Q,K,V) = \text{softmax}\left(\frac{QK^T}{\sqrt{d_k}}\right) V
    \label{eqn:scaled_dot_product_attention}
\end{equation}
where $Q$ is a set of queries and the keys and values
are packed into matrices $K$ and $V$, and $d_k$ is the dimension
of the keys.  The softmax function of a vector $x$ is defined as
\begin{align}
\begin{split}
    \text{softmax}: \mathbbm{R}^n&\to  \mathbbm{R}^n\\
    x_i&\mapsto\frac{e^{x_i}}{\sum_{j=1}^n e^{x_j}}\,,
\end{split}
\end{align}
which is applied to the dot product of the queries with the keys. The scaling in the softmax in~\eqref{eqn:scaled_dot_product_attention} by a factor of $1/\sqrt{d_k}$ improves stability of the gradients.

\emph{Multi-head attention} \cite{Vaswani2017AttentionIA}
is a simple variant of attention that can lead to improved training.
In multi-head attention, $h\in \mathbb{N}$ different 
linear projections of the $d$-dimensional queries, keys, and values are
learned, to $d_q$, $d_k$, and $d_v$ dimensions, respectively. Attention 
is then computed for each of the projected queries, keys, and values,
which are then concatenated and projected again. The result is known as 
multi-head attention, with $h$ heads.

The Transformer \cite{Vaswani2017AttentionIA} is an encoder-decoder
language translation architecture that uses stacked multi-head attention
layers. Since we will be utilizing a memory-efficient modification 
of the Transformer, we refer the reader to the original literature
for further details.

\subsection{Reformer}

The Reformer is a new architecture, an efficient transformer, that makes 
a number of memory improvements with respect to the original Transformer
and related follow-ups.
In this section we review the essential elements of the Reformer, as presented
in \cite{Kitaev2019ReformerTE}.

\medskip

Perhaps the key improvement in the Reformer is the use of \emph{locality sensitive hashing} (LSH)
attention. The essential idea behind LSH attention is that, due to the exponential
dependence in the softmax in \eqref{eqn:scaled_dot_product_attention}, some keys
contribute much stronger to attention (for fixed query) than others. This means that the matrix $\text{softmax}(QK^T)$ is sparse and dominated by a few entries, and we want to only compute these dominant ones. This will improve the complexity from $\mathcal{O}(l^2)$ to  $\mathcal{O}(l\log l)$, which becomes especially important for long sequences.
In more detail, the softmax of a key $k_j$ contributes a factor $\text{exp}(q_i\cdot k_j)$ to the attention of a query $q_i$. One now wishes to find the keys $k_j$ with maximal 
$q_i\cdot k_j = |q_i| \, |k_j| \, \cos(\theta_{ij})$, i.e.\ finding keys that are nearest neighbors to $q_i$ in a high-dimensional vector space.

Formulated abstractly, a hashing function (or scheme) $h: \mathcal{V} \to \{1,\dots,b\}$ assigns a vector $x\in \mathcal{V}$ to one of $b$ hash values. In cryptography, $h$ is chosen such that the hash values $h(x)$ of nearby values $x$ are as uncorrelated as possible in order to avoid revealing whether a guessed secret $x$ is close to the actual secret. Here, we want the inverse situation: nearby values $x$ should be mapped to nearby hashes $h(x)$. Such a hashing scheme is called
locality-sensitive. An example for an LSH scheme uses 
\begin{equation}
h(x) = \text{argmax}([xR;-xR]),
\end{equation}
where $[u;v]$ denotes the concatenation of two vectors $u$ and $v$, $R$ is a random
matrix of shape $\text{dim}(x) \times b/2$, and argmax returns the index of the largest
vector component~\cite{Andoni2015PracticalAO}. The idea is that under the random projection, nearby vectors
will map to nearby vectors and thus receive the same hash with high probability.

Returning to computing the attention~\eqref{eqn:scaled_dot_product_attention}, we can now only evaluate those scalar products in $QK^T$ that contribute the most. The attention $a_i$ of a query $q_i$ is given by
\begin{equation}
a_i = \sum_{j \in \cP_i} \text{exp}(q_i\cdot k_j - z(i,\cP_i)) \, v_j\,.
\end{equation}
Here, $\cP_i:= \{j:i\geq j\}$ is the set that the query at position $i$ attends to, the exponential structure comes from the softmax, $z$ is a normalizing term for the softmax, and we have omitted the factor $1/\sqrt{d_k}$ for clarity. Note that the structure of $\cP_i$ ensures that the $i^\text{th}$ position in the query may only attend to itself and the prior positions~\cite{Vaswani2017AttentionIA}.

We now change this attention scheme by only paying attention to elements within the same hash bucket, i.e.\ we set 
\begin{equation}
\cP_i^\text{LSH} = \{j: h(q_i) = h(k_j)\}\,.
\end{equation}
As discussed above, the computational and memory gains arise because $|\cP_i^\text{LSH}|\ll |\cP_i|$. Sometimes (but rarely), similar vectors will fall in different hash buckets. The chance that this happens can be reduced by performing multi-round LSH attention, i.e.\ the Reformer uses $n_\text{hashes}$ distinct hashing functions, defined by distinct, random matrices $R$. 

Additional details of LSH attention in the Reformer include causal masking 
that ensures positions may only attend to prior positions, and also a chunking scheme
that allows for efficient batch processing. In practice, the input with batch-size
$N$ is a tensor of shape $[N,l,d]$ which the Transformer then turns 
into $Q, K,$ and $V$ via three different linear layers. However, for LSH attention
in the Reformer to make sense we need $Q=K$. Similarly, a shared-QK Transformer
is a Transformer that has $Q=K$, and it turns out \cite{Kitaev2019ReformerTE}
that this has little effect on performance. Further improvements are achieved by 
using reversible layers.

In summary, the Reformer is a modern NLP architecture where improvements relative
to the Transformer allow sophisticated sequence data to be trained effectively
on a single GPU, bypassing the need for extensive computational resources and therefore
allowing easy exploration of new domains with NLP techniques. The most important hyperparameters
introduced by the Reformer are the number of hashes $b$ in LSH attention, and also
the choice of LSH attention or full attention, for the sake of comparison.

\subsection{Knots as Language}

Knots have various data presentation as words in appropriate sets of letters, which makes it natural to think of them as language.\footnote{An NLP that deals with letters and words would be to predict the next letter to be typed based on the letters that have already been input.} In this section we develop the idea in the context of 
natural language processing. We start by briefly summarizing some basics of knot theory, and then introduce the braid representation of a knot, which we use in most of our analysis and which can be interpreted as language.

\begin{figure}[t]
\centering
\includegraphics[width=0.18\textwidth]{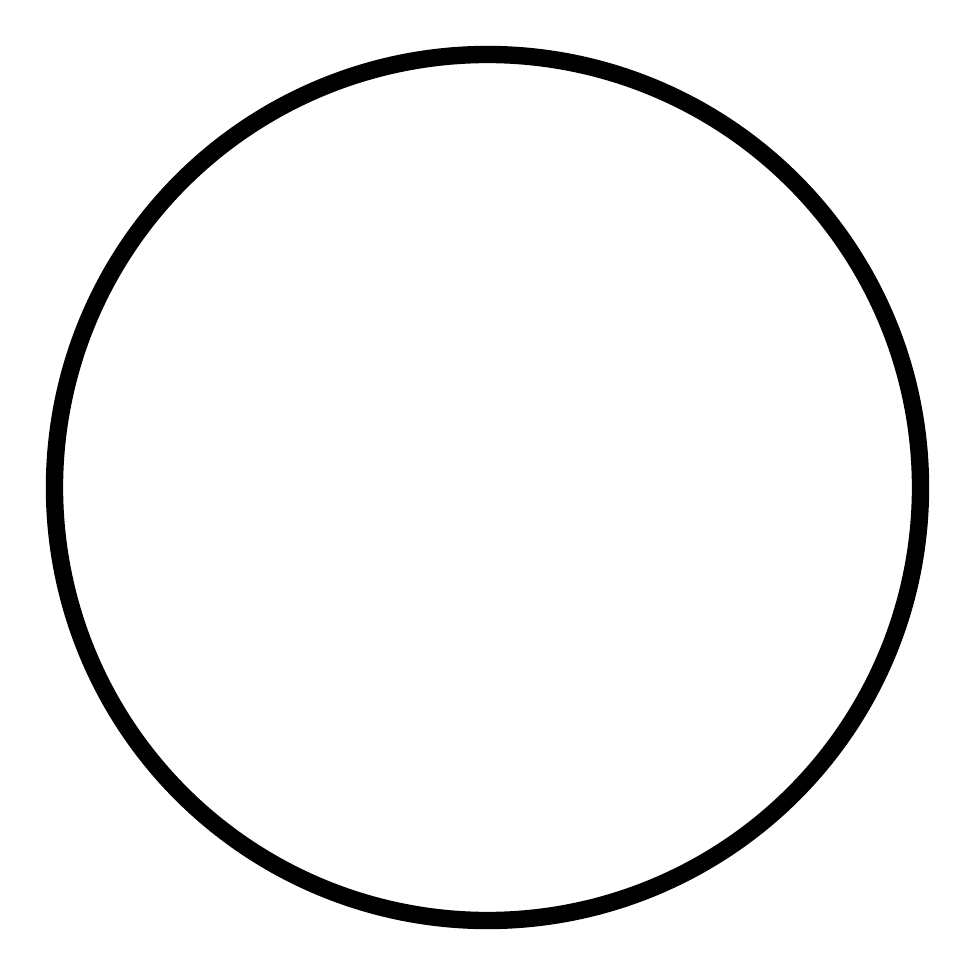}~
\includegraphics[width=0.18\textwidth]{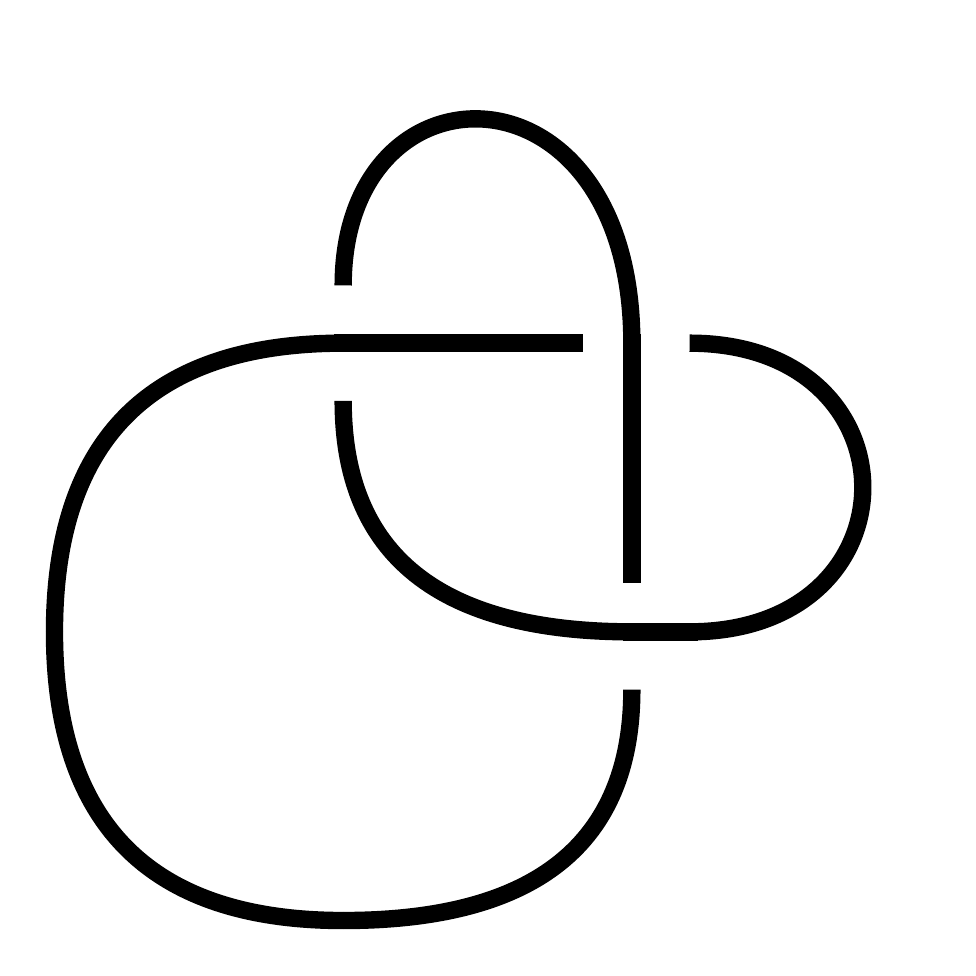}~
\includegraphics[width=0.18\textwidth]{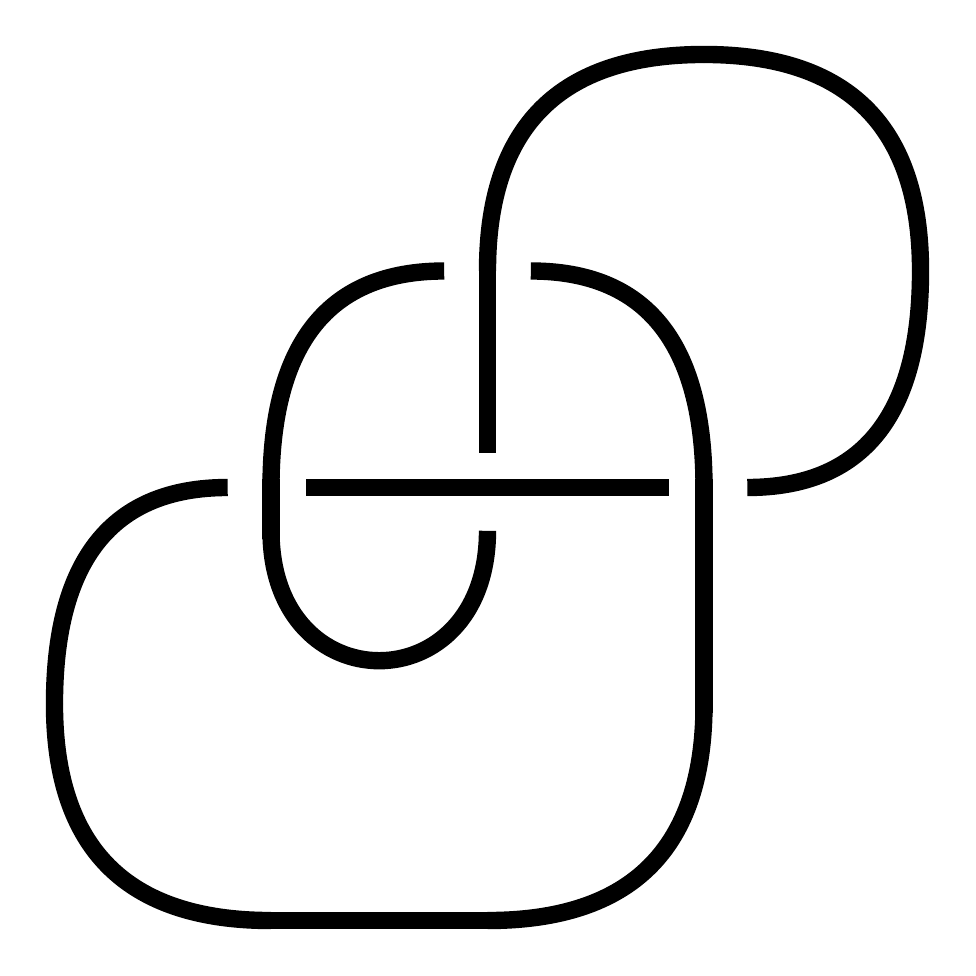}~
\includegraphics[width=0.18\textwidth]{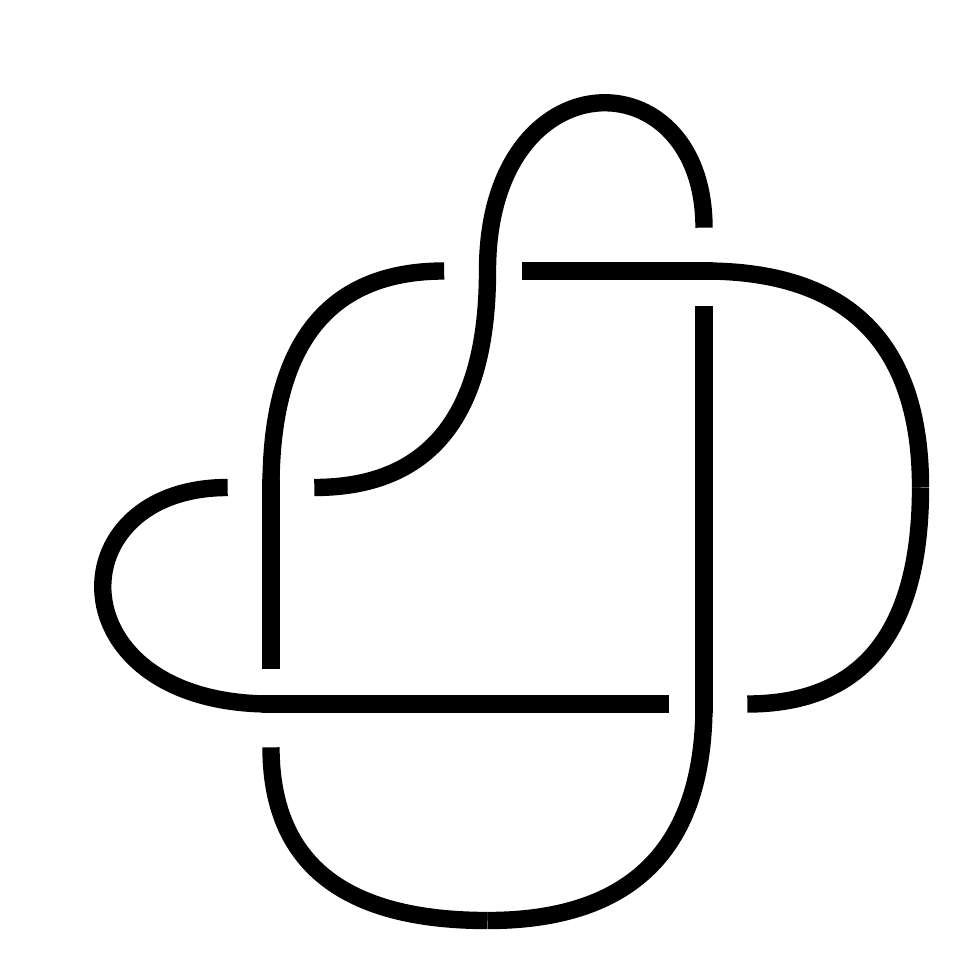}~
\includegraphics[width=0.18\textwidth]{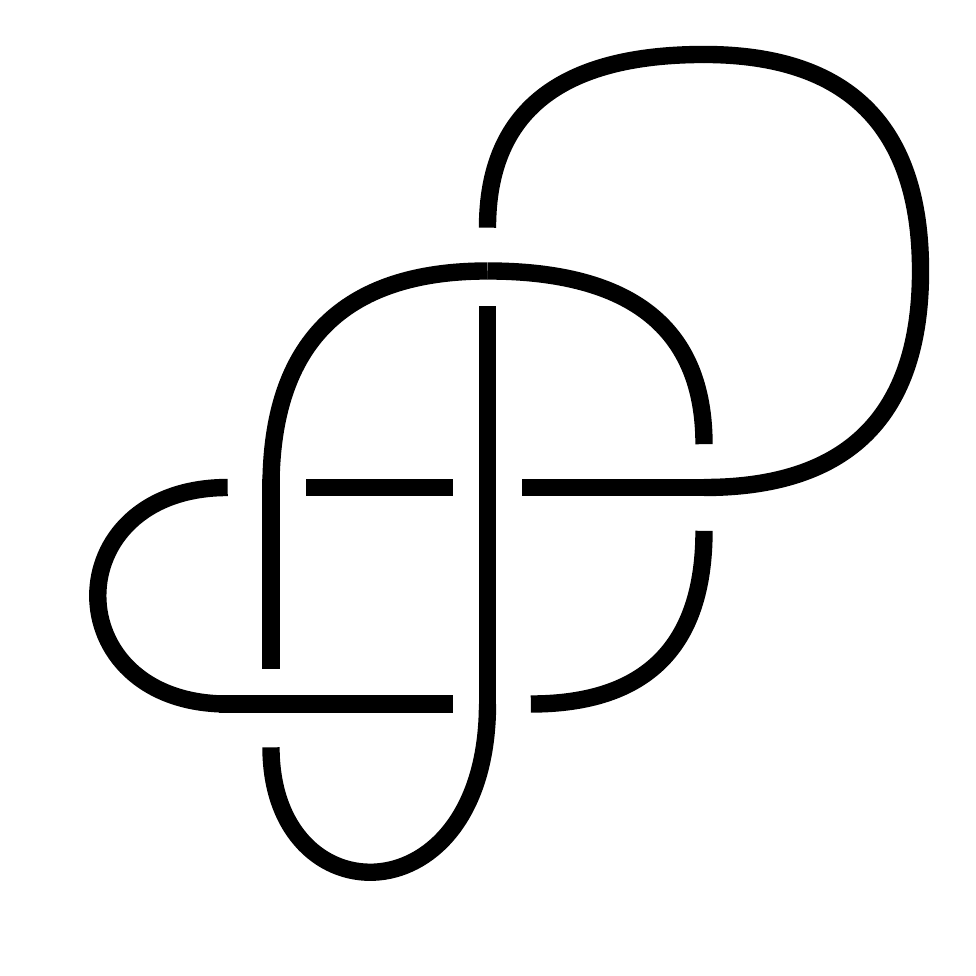}
\caption{Examples of knots. From left to right: unknot ($0_1$), trefoil ($3_1$), figure-eight ($4_1$), $5_1$, and $5_2$.}  \label{fig-knots}
\end{figure}

\begin{figure}[t]
\centering
\subfloat[][Type I: Twist]{\includegraphics[width=0.29\textwidth,valign=c]{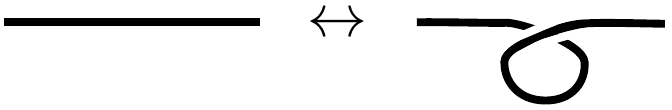}\label{fig:Reidemeister_1}}\qquad
\subfloat[][Type II: Poke]{\includegraphics[width=0.29\textwidth,valign=c]{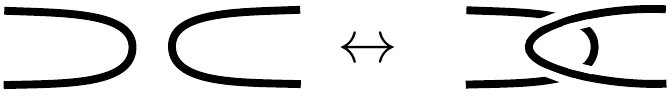}\label{fig:Reidemeister_2}}\qquad
\subfloat[][Type III: Slide]{\includegraphics[width=0.29\textwidth,valign=c]{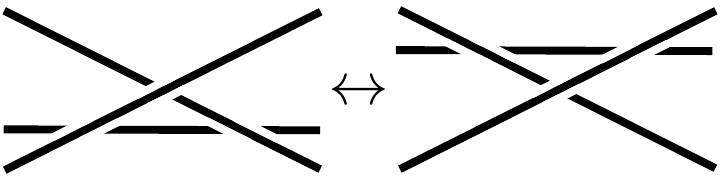}\label{fig:Reidemeister_3}}
\caption{Reidemeister moves.}  \label{fig-reide}
\end{figure}

A knot is an embedding of $S^1$ in 3-dimensional space, without self-intersections and up to ambient isotopy. The main goal of knot theory is to classify all knots, and to develop tools that enable to determine whether two different embeddings of $S^1$ are topologically equivalent, i.e.\ whether they represent the same knot -- in other words, whether one can be transformed onto the other without cutting. An important specialization of this problem that we address in this paper is to determine whether a given knot is topologically equivalent to the unknot, i.e.\ an unknotted loop, also referred to as the trivial knot. A collection of several possibly entangled knots is called a link. 

\begin{figure}[t]
\centering
\includegraphics[width=0.2\textwidth]{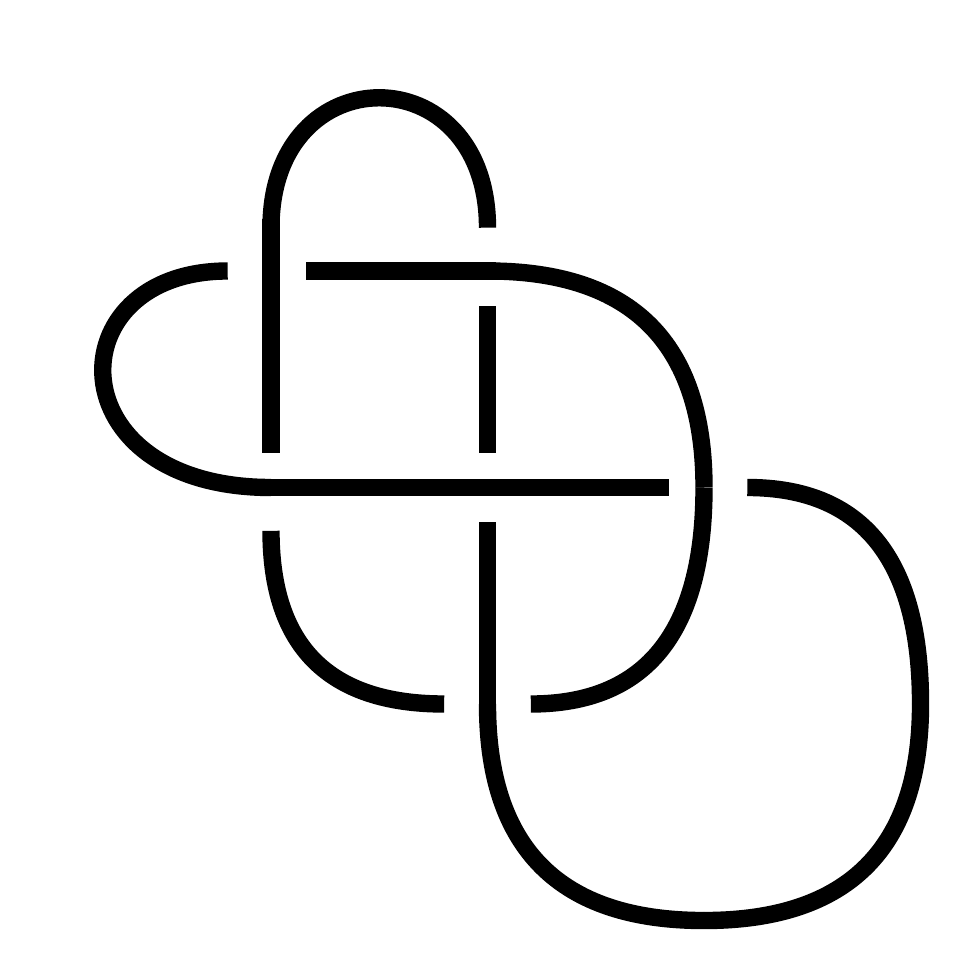}~~\raisebox{1.2cm}{+}~~\includegraphics[width=0.2\textwidth]{3_1}~~\raisebox{1.2cm}{=}~~\includegraphics[width=0.2\textwidth]{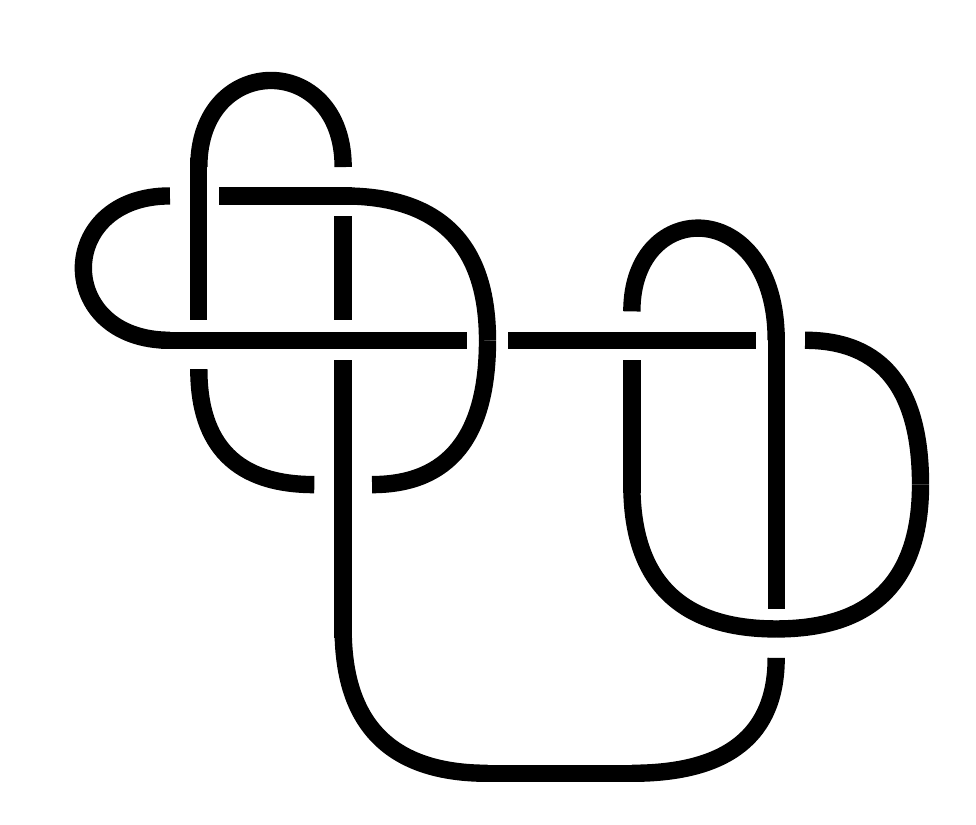}
\caption{Obtaining new knots as the sum of prime knots. This knot is the sum of the the knot $5_2$ (left) and the trefoil $3_1$ (right).}\label{fig:knot-sum}
\end{figure}

One useful approach to analyze knots is to consider their projections on a plane, see Figure~\ref{fig-knots}. Two knots are topologically equivalent if and only if their projections can be related to each other by a sequence of Reidemeister moves. These are three special moves that involve one, two, or three strands, see Figure \ref{fig-reide}: 
\begin{itemize}
\item A twist (Figure~\ref{fig:Reidemeister_1}) takes a strand and twists it, changing the crossing number by 1,
\item A poke (Figure~\ref{fig:Reidemeister_2}) pulls one strand over another, changing the crossing number by~2,
\item A slide (Figure~\ref{fig:Reidemeister_3}) slides a strand over (or under) a crossing of two strands, not changing the crossing number.
\end{itemize}
Furthermore, the most basic characteristic of a knot is the minimal number of crossings that one gets upon its projection onto an (appropriately chosen) plane. The simplest knots are the unknot, trefoil and figure-eight knot, denoted respectively $0_1$, $3_1$ and $4_1$, whose (minimal) numbers of crossings are given by the main number in this notation (i.e.\ 0, 3 and 4), while the subscript labels inequivalent knots with the same number of crossings. The unknot, trefoil and figure-eight are the only knots with less than 5 crossings. For a fixed, larger number of crossings there are many topologically inequivalent knots, e.g. there are 2 knots with 5 crossings (denoted $5_1$ and $5_2$). In addition to these unique prime knots, new ``composite'' knots can be formed as the sum of two or more prime knots. This can be thought of as taking two or more prime knots, cutting them open at one position, and tieing the open ends of each knot together, c.f.\ Figure~\ref{fig:knot-sum}.

The number of inequivalent knots (and indeed already the number of inequivalent prime knots) with a given number of crossings grows rapidly, so more elaborate characteristics must be employed to encode their structure and to distinguish them. For example, there are 165 prime knots with 10 crossings, 1,388,705 prime knots with 16 crossings, etc. 

A given knot clearly has many representations; for example projections on various planes typically look different, and in particular may yield different numbers of crossings. Therefore, one issue we have to deal with is how to represent the structure of a given projection. The second issue one needs to deal with is how to determine whether different representations represent topologically the same type of knot. Let us briefly discuss these two points.

In order to determine a type of a knot, so-called knot invariants are constructed. Knot invariants are various mathematical objects (numbers, polynomials, groups, homologies, etc.) which depend only on the topological type of a knot, and have the same form irrespective of the representative used to compute it. To prove that a given quantity is a knot invariant, it is sufficient to show that it is invariant under each of the Reidemeister moves. Note that if an invariant computed for two knots yields two different values, it means that these knots are inequivalent. On the other hand, if two knots yield the same invariant, they may be either equivalent or inequivalent. More powerful invariants distinguish more knots from each other, and a knot theorist's dream is to find a simple and practically computable invariant that would distinguish all knots.

For various purposes, in particular in order to compute various invariants, one needs to encode topological structure of a knot succinctly. The most common strategy to this end is to capture the pattern of crossings in a projection of a knot on a plane; it is clear that such a pattern determines a type of knot under consideration. Note that there are two types of crossings: once we traverse a knot, we may pass under or over each crossing that we come across. Keeping track of this information while we travel along the knot enables us to reconstruct its structure, and one way to capture this information is to use the Dowker-Thistlethwaite notation. 

\subsubsection*{Dowker-Thistlethwaite} 
To encode the structure of a knot in this notation, we traverse the knot and label each of the $n$ crossings from $1$ to $2n$, since each crossing is visited twice. We subject this labelling to the additional rule that the even label gets a minus sign when the strand followed crosses over at the crossing. At the end of this process, each crossing is labeled by one even and one odd number (and the even numbers are either positive or negative). Order these two-tuples in order of increasing odd numbers. The Dowker-Thistlethwaite notation is defined to be the sequence of the signed, even numbers in these ordered tuples. 

While the Dowker-Thistlethwaite notation can be easily determined for a given diagram, it also has certain disadvantages; for example, it is difficult to implement Reidemeister moves in terms of this notation, and in order to analyze links some additional information must be provided. For these reasons, in most of this work we represent the structure of knot projections in another way, namely representing knots as braids and using braid notation.

\subsubsection*{Braids} 
Let us therefore summarize what braids are and how to use them to encode the structure of knots.  Recall that the (Artin) braid group Br$_n$ is a non-Abelian, infinite, finitely generated group acting on $n$ strands with generators $\sigma_1,\ldots, \sigma_{n-1}$ and their inverses $\sigma_1^{-1},\ldots, \sigma_{n-1}^{-1}$, which satisfy the relations
\begin{subequations}
\label{braid-relations}
\begin{align}
\text{Braid relation 1:~}&\qquad\sigma_i \sigma_{i+1} \sigma_i = \sigma_{i+1} \sigma_i \sigma_{i+1}\,,\label{eq:BR1}\\
\text{Braid relation 2:~}&\qquad\qquad\! \sigma_i \sigma_j = \sigma_j \sigma_i \qquad \text{for}\ |i-j|\geq 2\,,\label{eq:BR2}
\end{align}
\end{subequations}
and similarly for the inverses.

\begin{figure}[t]
\centering
\includegraphics[width=0.3\textwidth,valign=b]{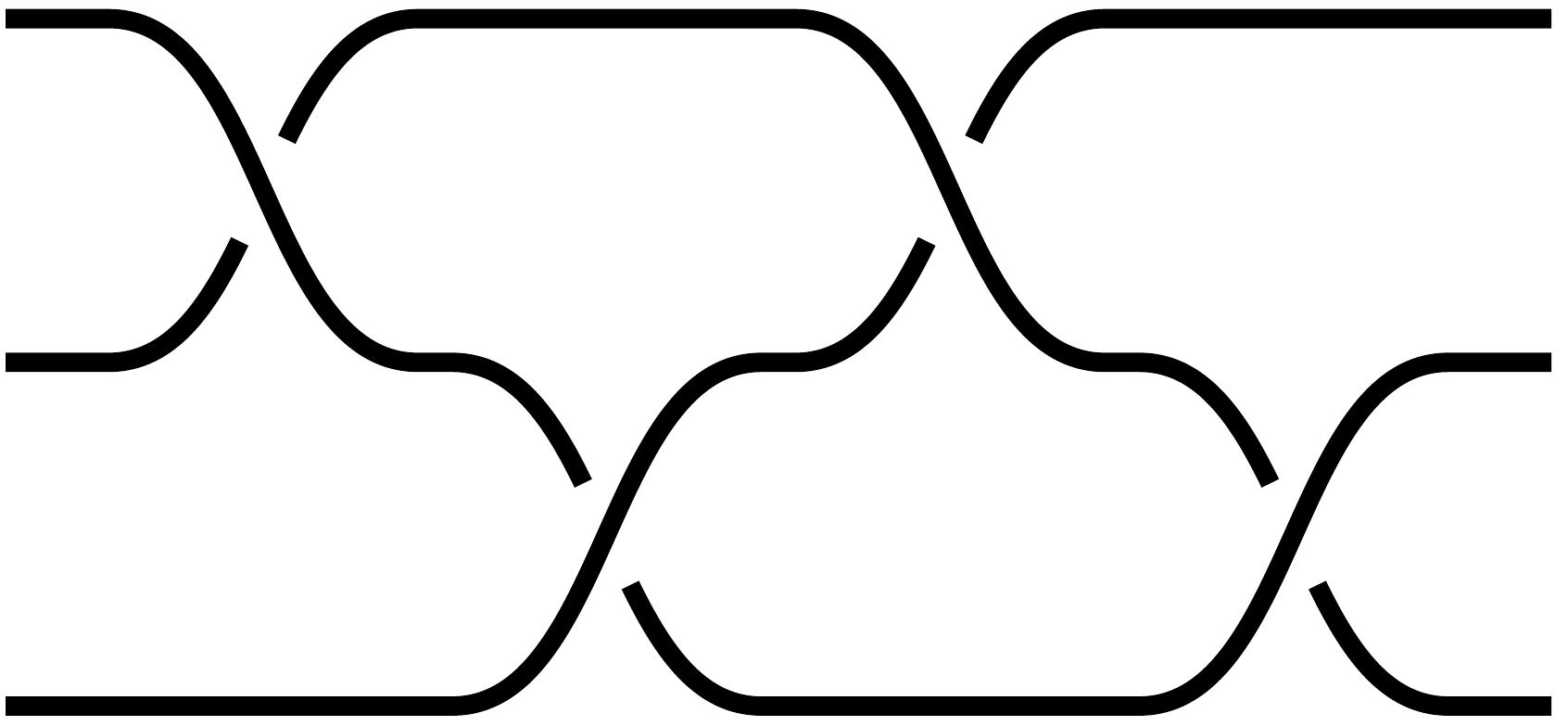} \qquad\qquad
\includegraphics[width=0.4\textwidth,valign=b]{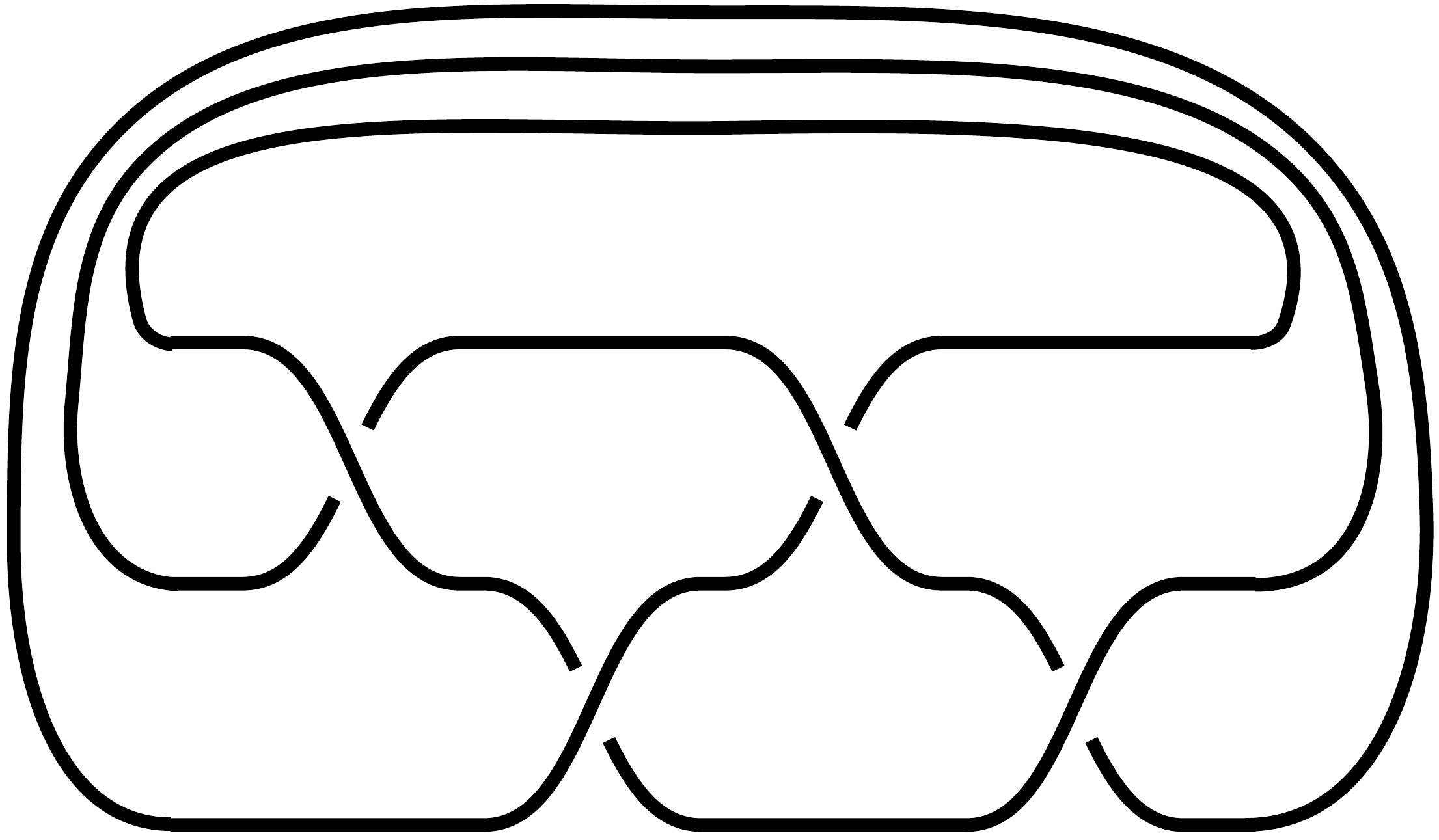} 
\caption{A braid $\sigma_1 \sigma_2^{-1} \sigma_1 \sigma_2^{-1}$ (left) and its closure (right).}  \label{fig-braid}
\end{figure}

\begin{figure}[t]
\centering
\subfloat[][Relation 1: $\sigma_1\sigma_2\sigma_1=\sigma_2\sigma_1\sigma_2$]{\includegraphics[width=0.45\textwidth,valign=c]{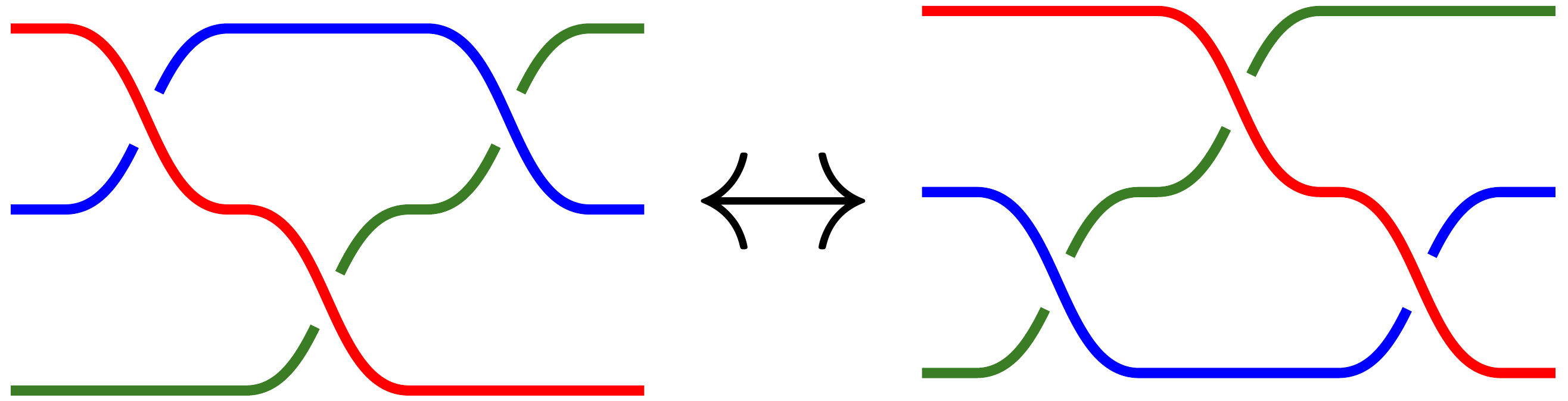}\label{fig:br1}} \qquad
\subfloat[][Relation 2: $\sigma_1\sigma_3\sigma_2=\sigma_3\sigma_1\sigma_2$]{\includegraphics[width=0.45\textwidth,valign=c]{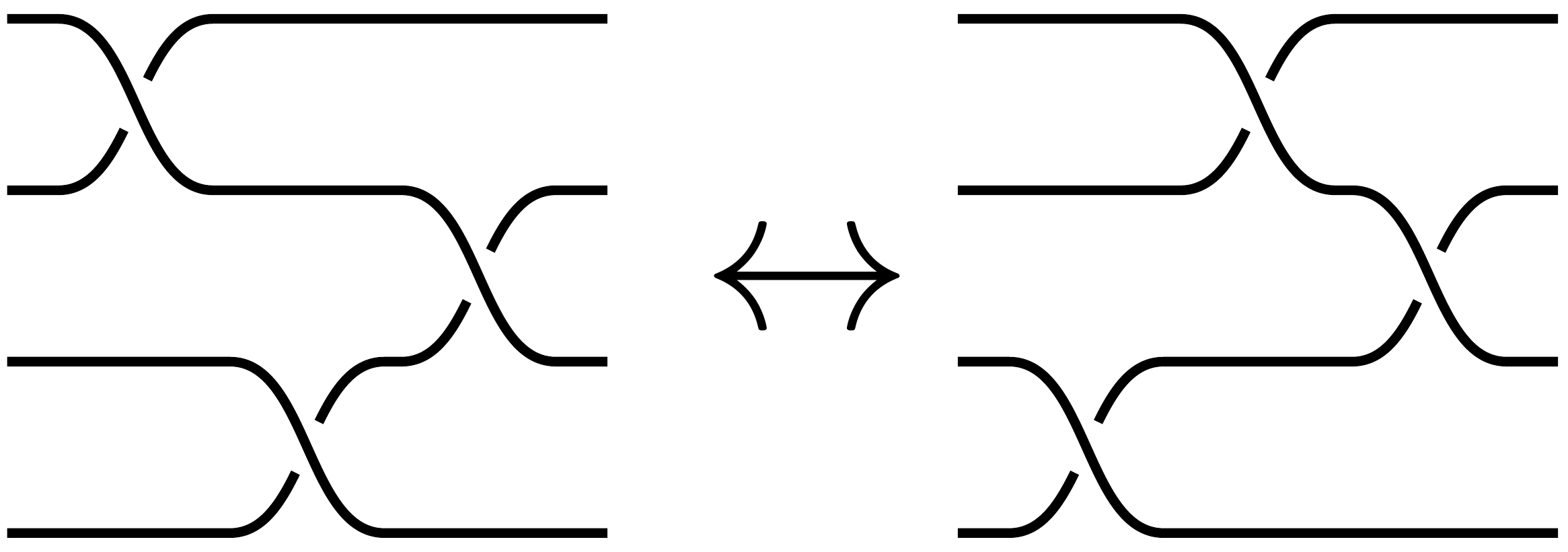} \label{fig:br2}}
\caption{Braid relations}  \label{fig:brain-relations}
\end{figure}

For a set of $n$ parallel strands, the generator $\sigma_i$ can be thought of as moving the $i^\text{th}$ strand over the $(i+1)^\text{st}$, and its inverse $\sigma^{-1}_i$ as moving the $(i+1)^\text{st}$ strand over the $i^\text{th}$ strand. A group element $\sigma_{i_1}^{\pm 1} \sigma_{i_2}^{\pm 1} \sigma_{i_3}^{\pm 1}\cdots$ can be represented as a pattern of interlacing strands and is referred to as a braid, see Figure~\ref{fig-braid} (left). A braid can be turned into a knot diagram by connecting beginnings and endpoints of all strands by a set of $n$ parallel arcs, as in Figure~\ref{fig-braid} (right). This operation is also referred to as closure. Furthermore, a theorem by Alexander states that each knot can be represented as a braid, and there is an effective algorithm that turns a knot into a braid. Note that a braid that we obtain upon such an operation may be regarded as a different knot projection, which of course represents the same knot type. 

\begin{figure}[t]
\centering
\subfloat[][Move 1 (conj.): $\sigma_1 \sigma_2^{-1} \sigma_1 \sigma_2^{-1}\leftrightarrow  \sigma_2^{-1} \sigma_1 \sigma_2^{-1} \sigma_1$]{~\includegraphics[width=0.45\textwidth,valign=b]{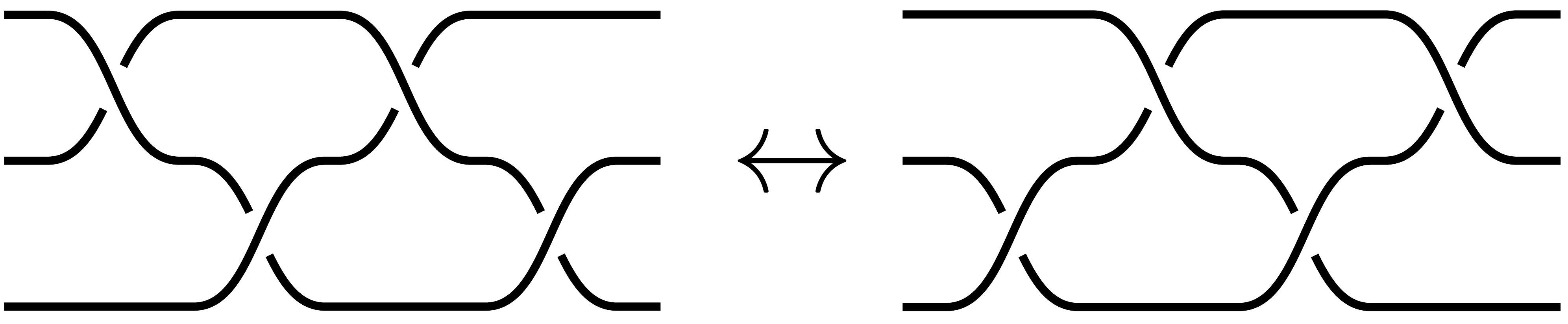}~\label{fig:mm1}} \qquad
\subfloat[][Move 2 (stabilization): $\sigma_1\sigma_2\sigma_1\leftrightarrow\sigma_1\sigma_2\sigma_1\sigma_3$]{\includegraphics[width=0.45\textwidth,valign=b]{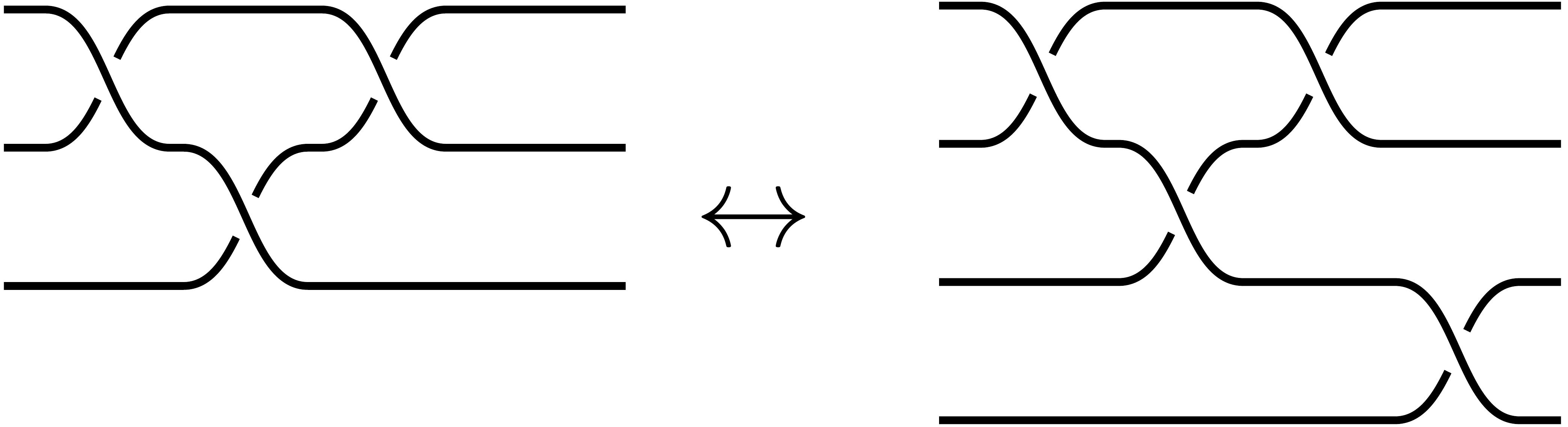} \label{fig:mm2}}
\caption{Markov moves.}  \label{fig-markov}
\end{figure}

As mentioned above, two different projections of the same knot can be related by a series of Reidemeister moves. There is an analogous statement on the level of braids, which is formalized in Markov's theorem. This theorem states that two braids that represent the same knot can be transformed into each other by a series of Markov moves. There are two types of Markov moves, shown in Figure~\ref{fig-markov}:
\begin{itemize}
\item conjugation and
\item stabilization/destabilization.
\end{itemize}
Conjugation sends a braid word $ww'$ to $w' w$. This can be achieved by repeating the following two-step procedure for every generator $\sigma_{i_j}$ in $w=\sigma_{i_1}\sigma_{i_2}\cdots \sigma_{i_l}$: First, we multiply $ww'$ by $\sigma_{i_1}^{-1}$ on the left and $\sigma_{i_1}$ on the right, giving
\begin{align}
\label{eq:conjugation_2step}
\sigma_{i_1}^{-1}ww'\sigma_{i_1}=\sigma_{i_1}^{-1}\sigma_{i_1}\sigma_{i_2}\cdots \sigma_{i_l} w'\sigma_{i_1} = \sigma_{i_2}\cdots \sigma_{i_l} w'\sigma_{i_1}\,.
\end{align} 
In the last step, we removed the consecutive inverses $\sigma_{i_1}^{-1}\sigma_{i_1}$. Repeating this for $\sigma_{i_j}$, $j=2,3,\ldots,l$ will result in sending $ww'$ to $w'w$. Note that the conjugation move corresponds, on the level of the knot or the braid closure, to inserting consecutive pairs of inverse generators, i.e.\ it acts trivially on the braid closure. Conjugation can be thought of as turning a knot into a braid by ``cutting it open at a different position''. Of course, one can in general  conjugate a braid word $w=\sigma_{i_1}\sigma_{i_2}\cdots \sigma_{i_l}$ with any generator $\sigma_{i_k}$, not just with $\sigma_{i_1}$, then $\sigma_{i_2}$, etc., sending it to $\sigma_{i_k}^{-1} w \sigma_{i_k}$ (or $\sigma_{i_k} w \sigma_{i_k}^{-1}$). On the level of the braid closure or the knot, this inserts the identity $\sigma_{i_k}^{-1}\sigma_{i_k}$.

Stabilization and destabilization are given by 
\begin{align}
\label{eq:destabilization}
\text{Stabilization:~} w&\to w\sigma_{n}\,,\qquad \text{Destabilization:~} w\sigma_{n}\to w
\end{align}
for $w\in\text{Br}_{n}$. This changes the braid (in fact, it even changes the underlying braid group from Br$_{n-1}$ to Br$_{n}$, or the other way around), but not the knot. Note that if we did not change the braid group from Br$_{n}$ to Br$_{n-1}$ in a destabilization move, we would be left with one strand (the $n^\text{th}$) which would not be acted on by any of the braid generators. As a consequence, we would have a two-component link: the first component would be a braid that describes a knot equivalent to the one we started with, and the second component would be the unknot, corresponding to the closure of the $n^\text{th}$ strand. Since this is not desired, we take destabilization to remove the generator $\sigma_{n}$ and change the braid group.

Note that there is a close connection between the Reidemeister moves and the Markov moves together with the braid relations:
\begin{itemize}
\item Reidemeister move 1 (twist) corresponds to Markov move 2, i.e.\ (de-)stabilization.
\item Reidemeister move 2 (poke) corresponds to adding a trivial element $\sigma_i\sigma_i^{-1}$ at some position in a braid word.
\item Reidemeister move 3 (slide) corresponds to the action~\eqref{eq:BR1} of first braid relation on the closure of the braid.
\end{itemize}

At this stage, we can finally relate braid representations of knots to language. We simply interpret generators $\sigma_i^{\pm 1}$ of Br$_n$ as letters, and braids of the form  $\sigma_{i_1}^{\pm 1} \sigma_{i_2}^{\pm 1} \sigma_{i_3}^{\pm 1}\cdots$ (which represent knots after the closure) as words. In practice and in what follows, we represent a braid generator $\sigma_i^{\pm}$ simply by $\pm i$, so that a word is represented by a string of integer numbers $i\in[-(n-1),n-1]$ that represents the braid.

There are several crucial points from the language perspective that should be stressed. First, we wish to identify, and treat as equivalent, different words (braids) that represent the same knot. This means that the AI needs to learn to identify such equivalent words. To conduct such a learning process we also need to generate equivalent words. To this end, we can take advantage of Markov moves. In particular, in the unknot problem, we can generate various representations of the unknot by applying a series of Markov moves to the empty braid. Furthermore, the topological character of knottedness makes the problem global rather than local: even a single change of a word in the sentence may change the type of knot under consideration; there is no notion of a ``small'' error, which makes the learning process hard; but this is also true for applictions to NLP.

\subsection{The UNKNOT Problem}

In this section we introduce the main problem that we study: the UNKNOT problem.

\subsubsection{Why unknotting?}
\label{sec:Unknot-motivation}

While the problem of distinguishing knots is interesting in its own right, much of our motivation comes from the smooth 4-dimensional Poincar\'e conjecture (or, SPC4, as it is often called).
Indeed, many problems in topology of 4-manifolds, including SPC4, can be described (and, sometimes, completely reduced) to the language of knots in $S^3$.

At the most basic level, the reason is that every closed smooth 4-manifold $M_4$ can be represented by a Kirby diagram, which basically consists of knots drawn on a 3-sphere $S^3$. More precisely, to build an $M_4$ one starts with a 0-handle, {\it i.e.} a 4-ball $B^4$, then attaches 1-handles, then 2-handles, 3-handles, and finally a 4-handle, which is also a 4-ball. In fact, since this last step involves no ambiguity, we don't need to attach the 4-handle. Either way, the relation to knots in $S^3$ comes after all $k$-handles with $k \le 3$ are attached.\footnote{Recall, that a four-dimensional $k$-handle is $B^k \times B^{4-k}$, which attaches onto the boundary of lower-index handles along $\partial B^k \times B^{4-k}$.}

There are many candidate counterexamples to SPC4, {\it i.e.} ``exotic'' spheres $M_4$ homeomorphic to $S^4$ which are not known to be diffeomorphic to $S^4$. One way to show that such an $M_4$ is the standard 4-sphere is to use equivalence relations (Kirby moves) to reduce its Kirby diagram to that of $S^4$, which has no $k$-handles with $k=1,2,3$. This problem is basically the unknotting problem, or a close variant of it.

Since, as mentioned earlier, adding a 4-handle is a fairly unambiguous operation, one often works with close relatives of SPC4 that involve $B^4$ with $S^3$ boundary in place of $S^4$. For example, the corresponding version of SPC4 is known as the smooth {\it relative} 4-dimensional Poincar\'e conjecture. If true, it implies the original SPC4. As in the case of SPC4 itself, there are many candidate exotic\footnote{An exotic 4-ball has no smooth radius function with 3-sphere levels.} 4-balls, {\it i.e.} $M_4$ homeomorphic to $B^4$ which are not known to be diffeomorphic to it.
For example, every knot $K \subset S^3$ which is fibered and ribbon gives such a candidate $M_4$ since, according to Casson and Gordon \cite{MR722728}, it bounds a fibered disk $D \subset M_4$ in some $M_4$ which is homeomorphic to a 4-ball $B^4$ but is not known to be diffeomorphic to it.
Therefore, if a fibered ribbon knot $K$ does not bound any fibered disk in $B^4$, then the smooth relative 4-dimensional Poincar\'e conjecture is false.

Conceptually, this is the same reason why knots in $S^3$ can tell us about smooth structures in one dimension higher that we already mentioned earlier. A knot $K = \partial \Sigma$ appears as a boundary of a surface $\Sigma \subset M_4$, and the question is whether $\Sigma$ can be a disk in $B^4$ or only in homotopy-$B^4$.
Whether $K \subset S^3 = \partial B^4$ bounds a disk in $B^4$ is controlled by the 4-ball genus (a.k.a. slice genus), $g_4 (K)$, which is defined to be the minimal value of $g(\Sigma)$, such that $\Sigma \subset B^4$ is bounded by $K$. A knot $K$ with $g_4 (K)=0$ is called {\it slice}.

Then, the strategy \cite{MR2657647} to disprove (relative) SPC4 could be to take a knot $K \subset S^3$ that is slice ({\it i.e.} bounds a disk) in a homotopy 4-ball $M_4$, with $M_4 \ne B^4$, and show that $K$ is not slice in $B^4$. For this, one needs obstructions to sliceness, {\it i.e.} lower bounds on $g_4 (K)$. One such bound comes from deformations and spectral sequences in Khovanov homology, namely the Rasmussen's $s$-invariant \cite{MR2729272}.
It bounds the 4-ball genus
\be
\frac{|s(K)|}{2} \; \le \; g_4 (K)
\ee
More generally, one may hope to find exotic 4-balls by looking for knots that exhibit different genus bounds in $B^4$ and in $M_4 \simeq B^4$. In \cite{MR2657647}, this strategy was applied to co-cores of 2-handles, which are disks in $M_4 \simeq B^4$ bounding knots and links in $S^3$. All those $M_4$ were soon shown to be standard \cite{MR2680408}.

One can also consider knots with the trivial Alexander polynomial, $\Delta_K (x) = 1$. In the early 1980's Freedman showed that all such knots are topologically slice \cite{MR1201584}. Therefore, demonstrating that any such knot has $g_4 (K) > 0$ would immediately imply the existence of an exotic 4-ball.
A similar conclusion follows if any fibered ribbon knot, as discussed above, has $g_4 (K) > 0$.

\subsubsection{Complexity}

After describing some motivation for unknotting, let us see how hard it can be.

More than 20 years ago, Hass-Lagarias-Pippenger \cite{MR1693203} proved that the unknotting problem, {\it i.e.} the decision problem whether a given knot $K$ is actually an unknot, is in complexity class NP (``Nondeterministic Polynomial-time'' Turing machine).
This is the complexity class that, famously, contains P (class of problems\footnote{It includes problems like multiplication and sorting.} for which ``Polynomial-time'' algorithms are possible) but is not known to (and, in fact, widely not believed to) be equal to it.
Problems in class NP are like Sudoku puzzles; they may not have a simple algorithm to solve, but a proposed solution can be verified in polynomial time.
In other words, while problems in class P are the ones for which an answer can be found in polynomial time, problems in class NP are the ones for which checking the answer can be done in polynomial time, provided that the answer is yes.
The result of \cite{MR1693203} means that the unknotting problem joins the class of problems like protein folding, SAT (satisfying truth assignment), or the traveling salesman problem, which are also in class NP.

A close cousin of the class NP --- which, though not too likely, may be equal to it --- is the class coNP. It consists of decision problems whose negative answers can be checked in polynomial time, i.e.\ if the answer is no. If $\text{NP} \ne \text{coNP}$, then $\text{NP} \ne \text{P}$ (but the other direction is not known).
The unknot recognition problem turns out to be not only in class NP but also in the complexity class coNP. This was first shown by Kuperberg \cite{MR3177300} assuming the generalized Riemann hypothesis (GRH).
This assumption was later relaxed in \cite{Lackenby}, where it was also pointed out that, in the unlikely event that either the unknotting problem or its negation (called knottedness) is NP-complete, then $\text{NP} = \text{coNP}$.
To summarize,
\be
\text{unknot recognition} \; \in \; \text{NP} \, \cap \, \text{coNP}
\ee
This result is particularly interesting because many decision problems that originally started in this intersection --- {\it e.g.} deciding whether an integer number is prime or composite --- were later found to be in class P \cite{Agrawal02primesis}.
Therefore, there is a chance that the unknotting problem we are trying to tackle here actually admits a polynomial time algorithm. Approaching this problem via AI/ML can hopefully help us find such an algorithm, if it exists.

In fact, it has been a long standing problem whether the unknot recognition is truly more difficult than a similar problem for braids, the braid word problem. The latter is known to be in class P according to the Garside-Thurston theorem, which says that one can identify the trivial braid in polynomial time, $O(|\text{word length}|^2 \;n \log n)$ for the Artin braid group $\text{Br}_n$.
This can be improved to $O(|\text{word length}|^2 \;n)$ with the BKL algorithm~\cite{MR1654165}.\footnote{Note that the closure of a trivial braid group element is the unknot, but there can be non-trivial braid elements, whose closure is still the unknot. This is why the BKL algorithm does not solve the unknot problem in polynomial time.}

At the same time, perhaps one should not be overly optimistic.
For example, it was shown recently that imposing an upper bound on the number of Reidemeister moves immediately makes the unknot recognition problem NP-hard \cite{MR3968635}.
This paper also helps to understand how the unknotting problem compares to deciding whether two vertices of a given finite graph are connected or not, which is in class P.
Indeed, if we think about knot diagrams as vertices of an abstract graph, with edges representing Reidemeister moves, then the unknotting problem is equivalent to deciding whether a vertex belongs to the same component of the graph as the ``origin'' (the vertex associated with a trivial diagram of the unknot). If this abstract graph was finite and explicitly presented, then the unknotting problem would be in class P, but \cite{MR3968635} can be viewed as an indication that these two problems are qualitatively different.

Finally, since earlier we talked about computation of delicate knot invariants, it should be noted that many closely related problems were recently shown to be parsimoniously $\#$P-complete \cite{KuperbergSamperton}.
This is one of the more esoteric complexity classes, based on $\#$P which is larger than NP but is contained in PSPACE (``Polynomial-space''). And, ``parsimoniously complete'' refers to a more specific version of the completeness relation, such that for every solution of problem A there is a unique solution of problem B.
Note, the class PSPACE also contains NP and coNP that we discussed earlier, as well as the probabilistic version of the polynomial time solver (BPP).
Interestingly, both \cite{MR3177300} and \cite{KuperbergSamperton} use the representation variety $\pi_1 (S^3 \setminus K) \to G$ in a crucial way. When $G=SL(2,\C)$, this is the familiar A-polynomial that plays an important role in Chern-Simons theory~\cite{Gukov:2003na}.

\section{Generating Knots and Unknots \label{sec:generating_data}}
In this Section we describe the algorithm that we use to generate representatives of non-trivial knots and unknots, or 
alternatively the prior from which they are drawn. 
Details of all of the algorithms and subroutines are presented in Appendix \ref{sec:app_algos}.

In describing the prior we attempt to find a balance between being explicit about our
subroutines and explaining how they are sewn together to form our databases consisting of non-trivial knots and unknots.
Crucial subroutines include:
\begin{itemize}
\item \randommarkovmove, Algorithm \ref{alg:RandomMarkovMove}, performs a random Markov move drawn from a uniform distribution, changing the braid but not the topology of its closure.
\item \braidrelationone, Algorithm \ref{alg:BraidRelation1}, applies the first braid relation in~\eqref{eq:BR1}.
\item \smartcollapse, Algorithm \ref{alg:SmartCollapse}, iteratively removes consecutive inverses, free strands, twists (i.e.\ performs a destabilization move), and non-consecutive inverses (associated with inverses on opposite ends of the braid) until the braid no longer changes. 
\item \knotify, Algorithm \ref{alg:Knotify}, performs a sum over link components (as illustrated in Figure~\ref{fig:knot-sum}) by iteratively interweaving link components associated with a braid closure until only one component is left, i.e.\ the braid closure is a knot.
\end{itemize}
These play a role in the algorithms used to draw random non-trivial knots and unknots:
\begin{itemize}
\item \randomunknot, Algorithm \ref{alg:RandomUnknot}. Starts with the empty braid and iteratively applies \randommarkovmove and \braidrelationone a total of $M$ times before applying \smartcollapse, until the braid has length $n_\text{letters}$. 
\item \randomknot, Algorithm \ref{alg:RandomKnot}. As long as the generated braid $B$ does not have $n_\text{letters}$, draws a new $B$ of length $n_\text{letters}$ from a uniform distribution on the generators of the braid group of an input number of strands and then applies \randommarkovmove and \braidrelationone a total of $M$ times, followed by \smartcollapse.
\end{itemize}
Braids produced by \randomunknot have topologically trivial closure; we perform checks 
of necessary conditions by computing the Arf invariant, Alexander polynomial, and whether or not the knot is alternating and comparing to the unknot values $(0,1,\text{False})$. Braids produced by \randomknot could potentially have topologically trivial closure, though the probability of this occurring should be exponentially suppressed in $n_\text{letters}$. We check that the braids are topologically non-trivial by computing the same invariants and ensuring that at least one of them differs from the unknot values.

We use \randomknot and \randomunknot together with the topological checks to produce databases of over $10^4$ non-trivial knots and unknots of various length, with
\begin{equation}
n_\text{letters} \in \{12,24,36,48,72,96\}.
\end{equation}
Any value is possible, though generation time goes up significantly with $n_\text{letters}$.

\subsection*{Distribution of knots}

It would be interesting to study the distribution of knots introduced by the flat prior from which the Markov moves are drawn. In particular for small crossing numbers, all topologically inequivalent knots have been classified. This means one should in principle be able to check for our database of knots of this length, which of the inequivalent knots are produced (and how often). However, identifying the specific knot associated to our randomly generated braids is a hard problem that requires computing and comparing knot invariants (or ML). In general this problem is difficult and beyond the scope of our paper.

\begin{figure}
\centering
\includegraphics[width=.75\textwidth]{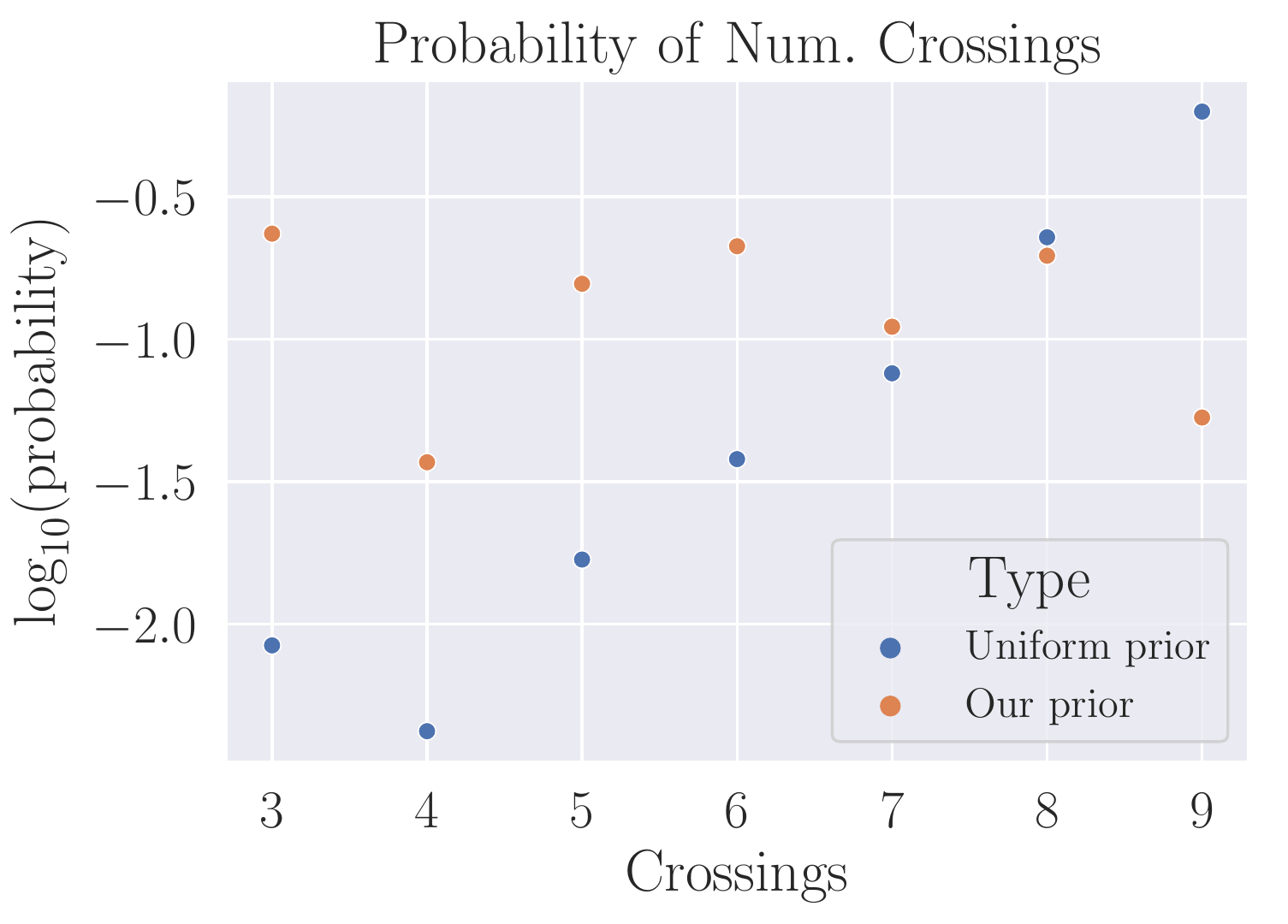}
\caption{Distribution on the number of crossings induced by our prior, and also by a uniform distribution on all knots of $9$ or fewer crossings.}
\label{fig:crossing_prob}
\end{figure}

However, knots with $9$ or fewer crossings have a particularly nice property: they may be uniquely identified by their Jones polynomial, and therefore $n_\text{letters}=9$ braids drawn from our prior have knots closures that may be identified. By computing the Jones polynomials of the prime knots of $9$ or fewer crossings (e.g., using knots in the Rolfsen table and their mirrors) one may compute the Jones polynomials of all knots of $9$ or fewer crossings by taking products. Drawing $6455$ $n_\text{letters}=9$ braids from our prior, we may identify the knot closure of each by its Jones polynomial. We see from Figure \ref{fig:crossing_prob} that the distribution on the number of crossings induced by our prior is much flatter than the one induced by a uniform distribution on knots with $9$ or fewer crossings, which grows exponentially. 

\begin{figure}[t]
    \centering
    \makebox[\textwidth][c]{
            \includegraphics[width=.49\textwidth]{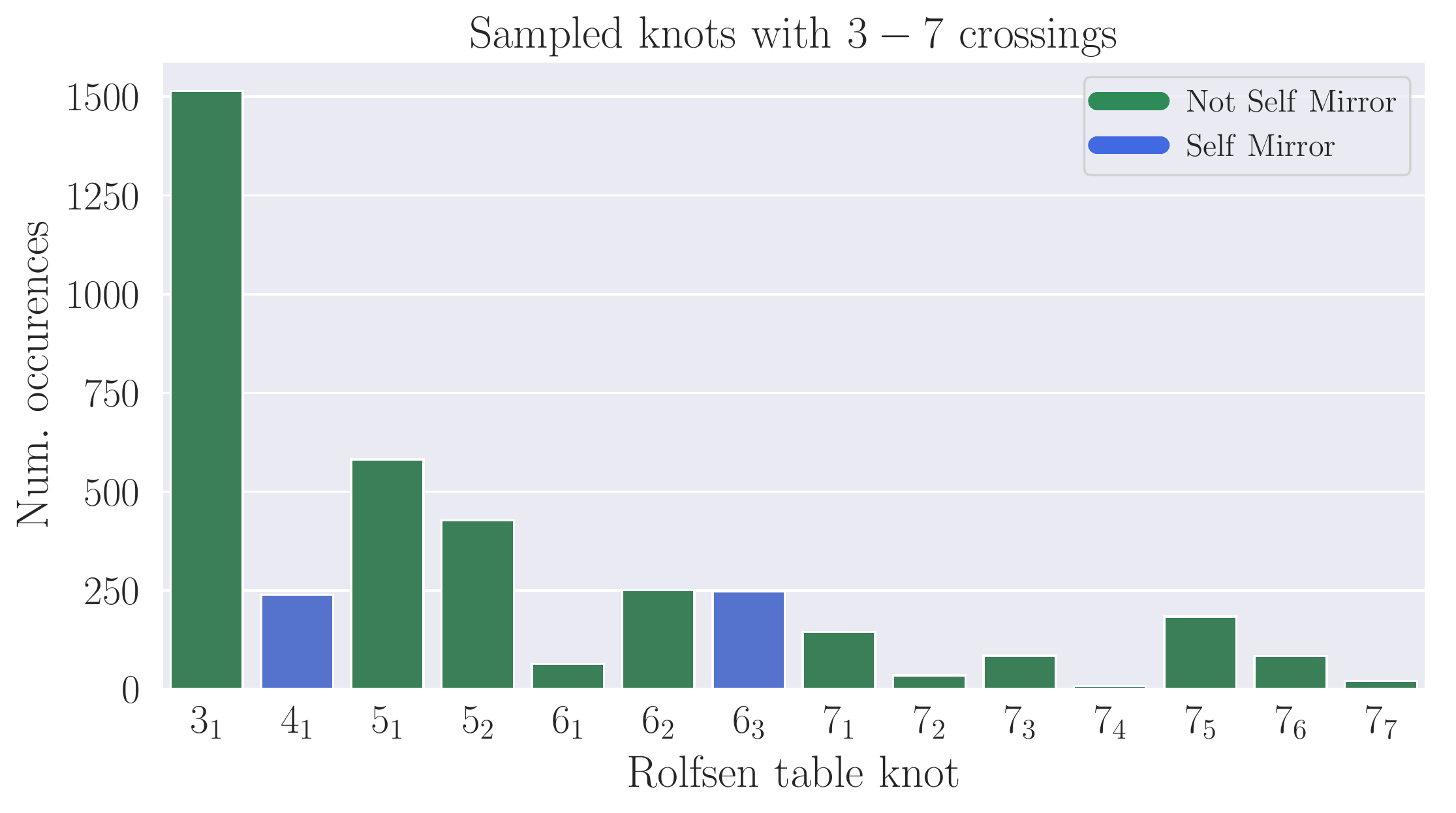}
            \includegraphics[width=.49\textwidth]{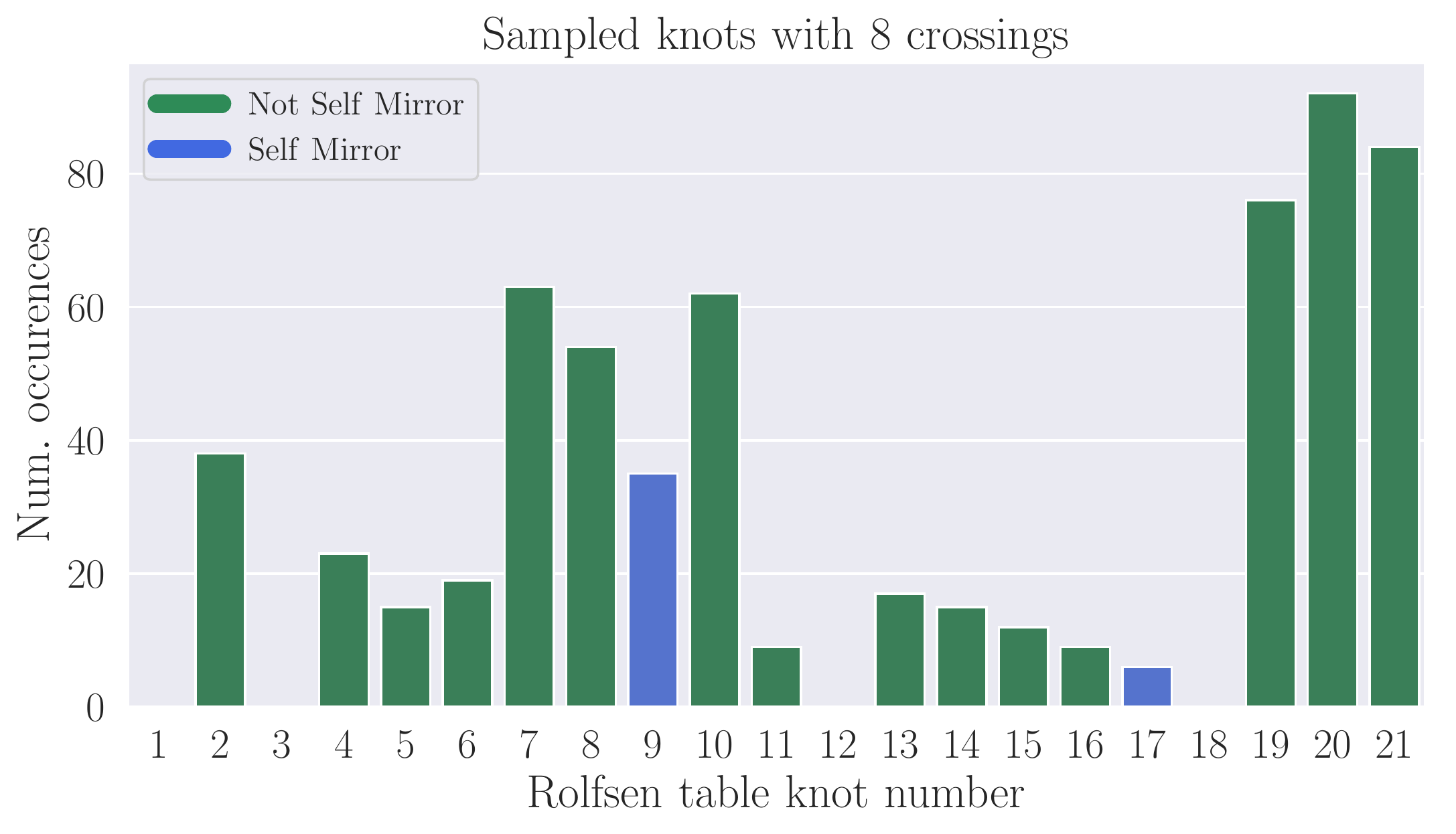}
    }
    \makebox[\textwidth][c]{
            \includegraphics[width=0.98\textwidth]{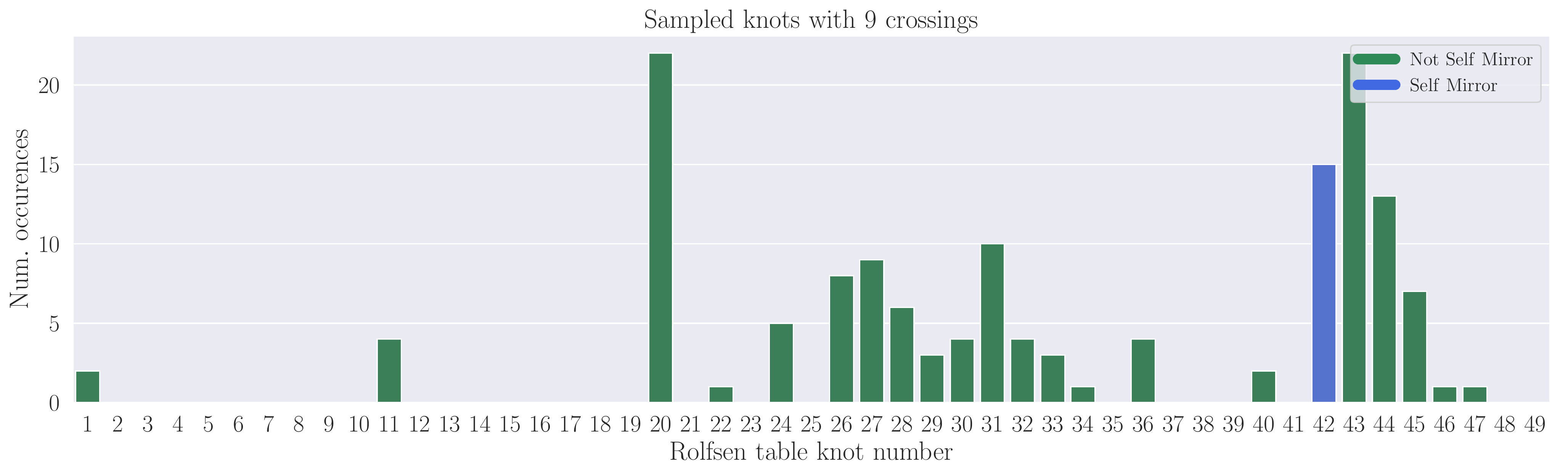}
    }
    \caption{Drawing $6455$ $N=9$ braids from our prior yields knots with $9$ or fewer crossings, $4664$ of which are prime. Plotted are the number of occurences of knots in the Rolfsen table for knots with $3$ through $9$ crossings, with mirrors counted for knots that are not self-mirror.
    }
    \label{fig:rolfsen_figs}
\end{figure}

\section{Unknot Decision Problem\label{sec:decision}}

\begin{table}[t]
    \centering
    \begin{tabular}{|c||c|c|}
	\hline    
    Parameters & Reformer Value &  Feedforward NN Values\\ \hline
    Full Attention & \{True, False\} & N/A \\
    $n_\text{hashes}$ & $\{1, 4\}$& N/A\\
    Bucket Size & \{$2$, $N/2$\} & N/A \\
    Causal & \{True, False\} & False (dense layers) \\
    Width & fixed by other hyperparameters & $850$\\
    Depth & $10$ & $10$\\
    Embedding Dimension & $250$ &250 \\
    Epochs & 50 & 50 \\ 
    Optimizer & RMSProp& RMSProp \\
    Learning Rate & $.0001$ & $0.01$ \\ \hline
    \end{tabular}     %\hspace {1cm}
    \caption{Parameter values for reformer and feedforward runs. Depth is the number of attention modules in the reformer, which themselves consist of multiple layers. Full attention means that a shared-QK transformer is used, rather than a reformer, i.e., full shared-QK attention is used rather than LSH attention. $n_\text{hashes}$ and bucket size are only meaningful when not using full attention.}
    \label{tab:decision_parameters}
    \end{table}

Given the motivations in Section~\ref{sec:Unknot-motivation}, in this section we study the UNKNOT decision problem. That is, given a representative of a knot, we wish to use supervised learning to determine whether or not it is the unknot. In Section \ref{sec:unkotting}, we will study braid representatives of knots and utilize reinforcement learning to find a sequence of Markov moves and braid relations that explicitly reduces it to the unknot (if possible), or to a braid word that is as short as possible otherwise.

We train Reformers, shared-QK transformers, and feedforward networks, on $10,000$ braids drawn from the prior of Section \ref{sec:generating_data} with $N\in\{12,24,36,48\}$, using the parameter values presented in Table \ref{tab:decision_parameters}. We label non-trivial knots as class $0$ and unknots as class $1$ and allow network outputs to vary between $0$ and $1$. 
For the Reformers and shared-QK transformers we use \texttt{reformer-pytorch}~\cite{Wang2020:reformer}. Specifically, for the Reformer runs we use the \texttt{ReformerLM} class, which applies an embedding layer and a Reformer module, which we then follow with a fully connected layer to map it to a single output, and finally a Sigmoid activation to ensure the output is between $0$ and $1$. Binary cross-entropy is used for the loss function, and we pick a decision threshold of 0.5. 

Let us comment on the embedding layers used. As described in Section~\ref{sec:Embedding-Layers}, NLP (which is one of the main areas where the reformer architecture is used) is dealing with words or letters (i.e.\ categorical data). These need to be converted to numerical data. This can be done via a one-hot encoding, where each word/letter is a vector of zeros (whose length corresponds to the number of words/letters in the dictionary, with a single one at the $i^\text{th}$ position, indicating the position of the word/letter in the dictionary). If a linear layer is applied to this one-hot encoded vector, it picks out the $i^\text{th}$ column of the corresponding matrix. An embedding layer combines the one-hot encoding with the linear layer, i.e.\ it corresponds to a matrix of dimensions (\texttt{embedding\_size})$\times$(\texttt{dictionary\_size}). Looking up the $i^\text{th}$ word/letter hence returns a vector of length \texttt{embedding\_size}. Note that in the case of braids, the numbers we assign to the generators are not completely arbitrary: While the numbering of the generators is arbitrary (although it is conventional to have $\sigma_i$ operate on strands $i$ and $i+1$), there is information in the fact that $\sigma_i$ and $\sigma_i^{-1}$ are inverses. We encode these generators as $i$ and $-i$, which are also inverse under addition. This fact could of course also be learned by the embedding layer during training (in the same way that it can learn that hot and cold are opposites in NLP), but we nevertheless noted that using an embedding layer is not necessary and we could use the (normalized) input directly as well. 

A second point is that the attention layers have no notion of how the (embedded) input was ordered in the original NLP. It can hence be beneficial to add this information to the chosen embedding. How beneficial this is depends on the language and its grammar; it is for example less important in Latin as compared to English. In the case of braids, most braid generators do commute but not all, see~\eqref{eq:BR2}. Encoding the position can be done via another embedding layer, which adds some vector to the embedded input. In this way, the positional encoding is done via another embedding layer of dimensions (\texttt{max\_input\_vectors}) $\times$ (\texttt{embedding\_dimension}), where \texttt{max\_input\_vectors} is the (maximum number of) input vectors that the reformer attends to. The Reformer can now reconstruct the position of the (embedded) input by looking up the vector that has been added to the input.\footnote{In practice, the reformer uses a more memory-efficient positional encoding known as Axial Positional Encoding. This works similarly but uses some tricks (factorizations) that allow to not store the full positional embedding matrix.}

We studied the dependence of performance on various parameters. We find that the performance difference between full attention (shared-QK Transformers) and Reformers is negligible. Moreover, the non-autoregressive model (i.e.\ \texttt{causal=False}) performs marginally better than the autoregressive one. Here autoregressive means that the future tokens are masked in the attention and the Reformer has to predict the next token only based on the previous ones (often in NLP, one wants to predict the next letter/word following the user input up to now). Since knots live in the closure of braids, there is no well-defined future words (the knot is turned into a braid by ``cutting it open'' at an arbitrary position), and hence one would expect the non-autoregressive transformer to perform better.  

\begin{figure}[t]
\centering
\subfloat[][Performance dependence on the number of locality sensitive hashes.]{\includegraphics[width=.47\textwidth]{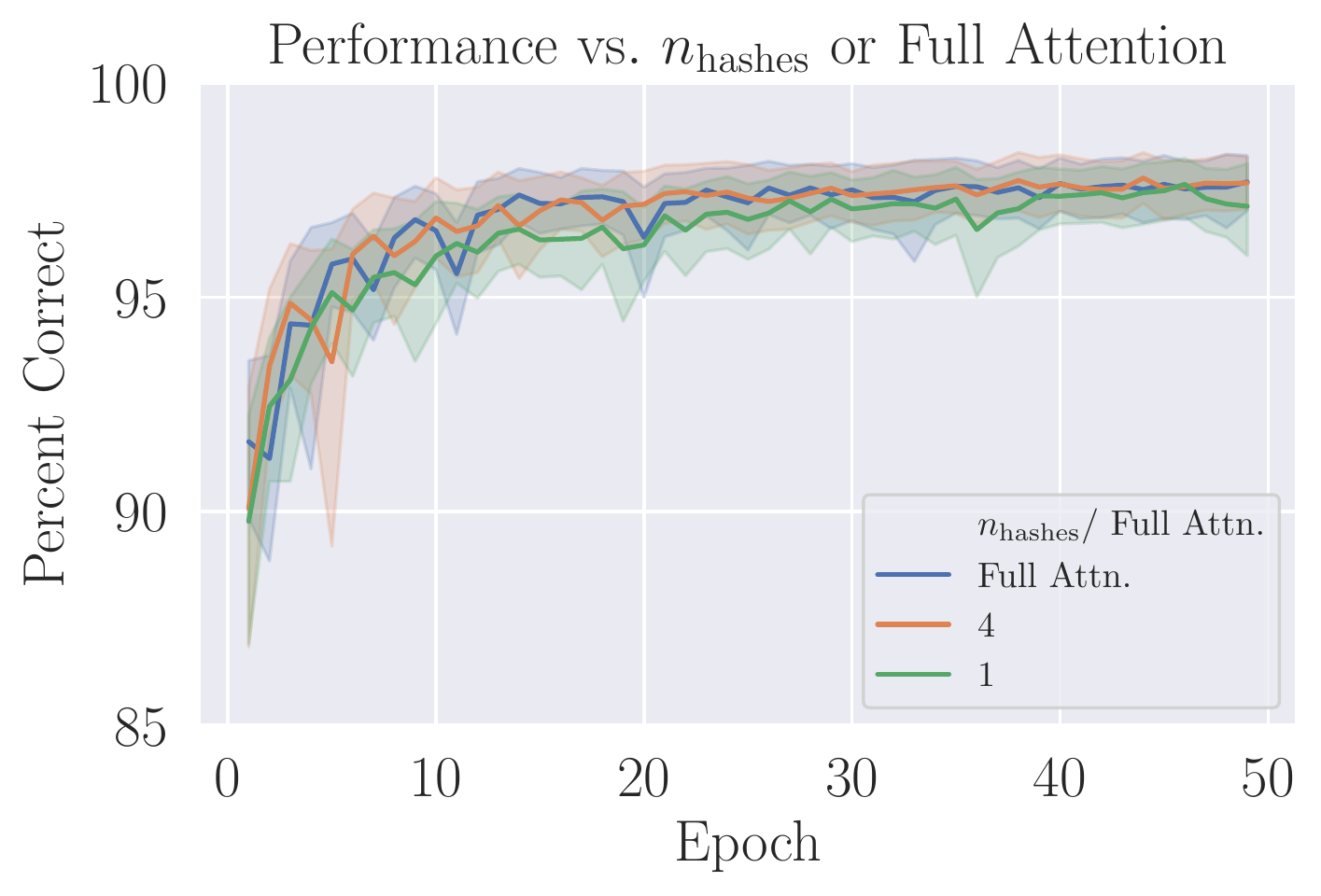}\label{fig:reformer_v_nhashes}}\qquad
\subfloat[][Performance  comparison between reformer and feedforward network.]{\includegraphics[width=.47\textwidth]{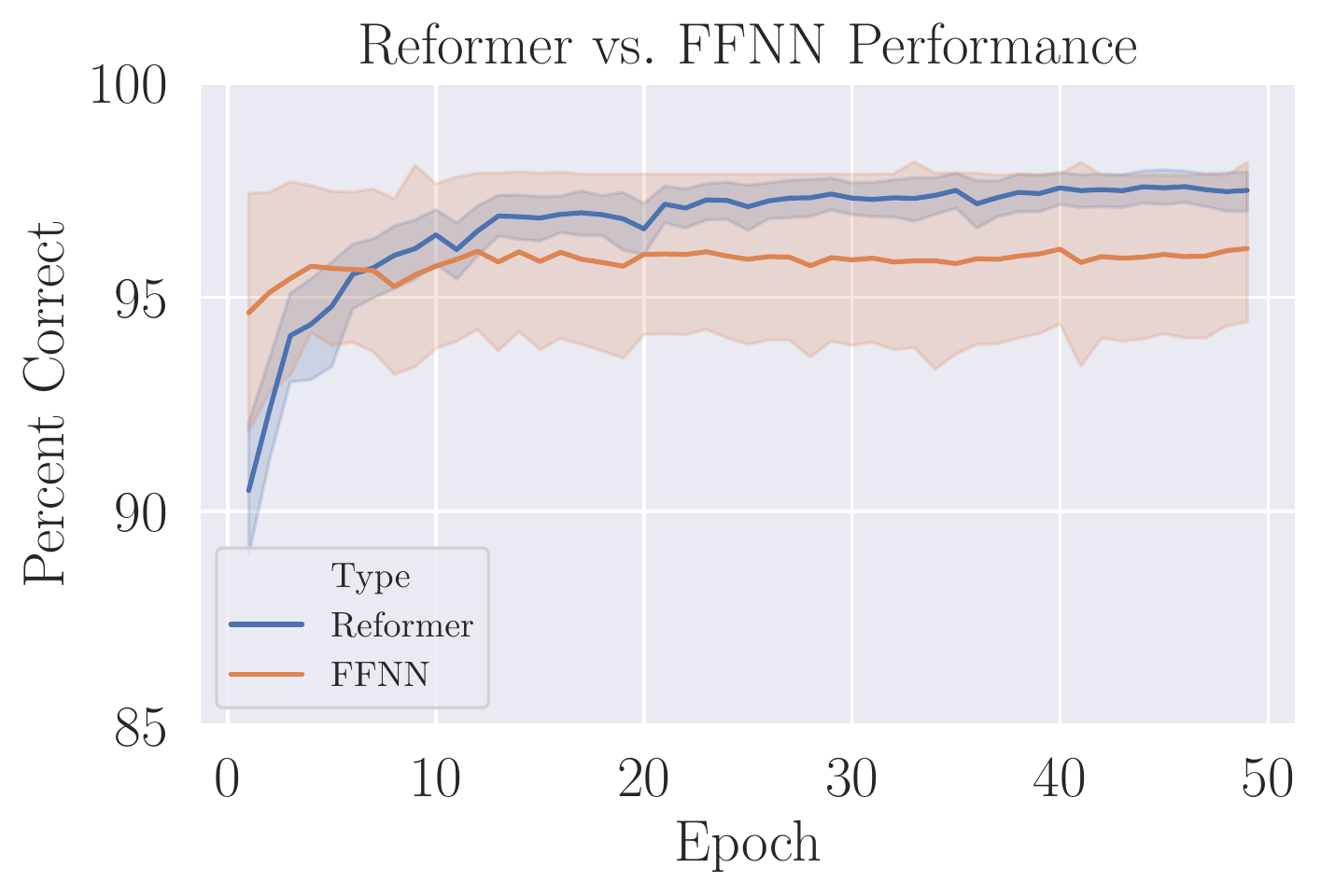}\label{fig:reformer_v_FFNN}}
\newline
\subfloat[][Performance dependence on the braid length. Performance increases with~$N$.]{\includegraphics[width=.47\textwidth]{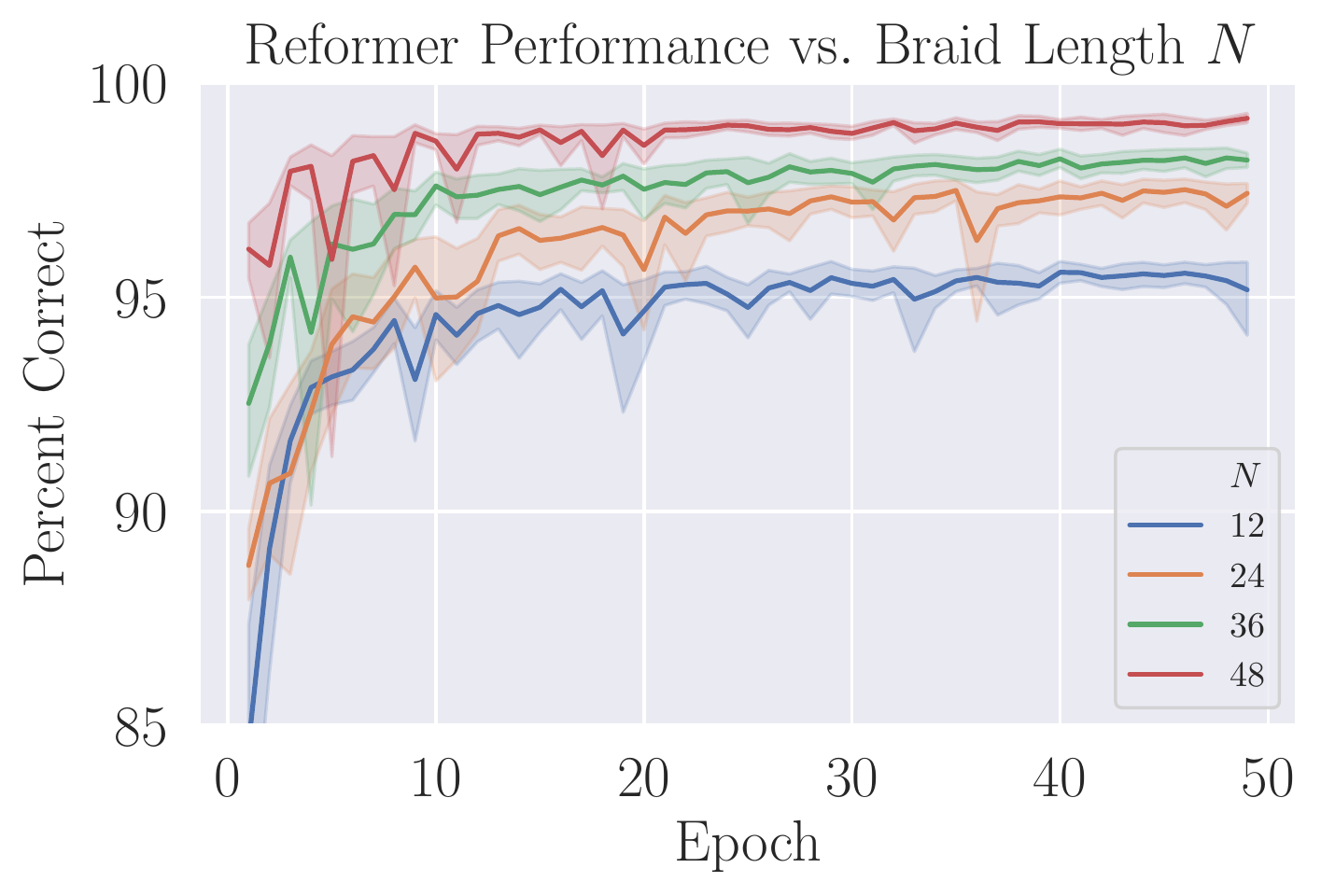}\label{fig:reformer_v_N}}\qquad
\subfloat[][Performance when number of braid letters, rather than number of braid words, is fixed.]{\includegraphics[width=.47\textwidth]{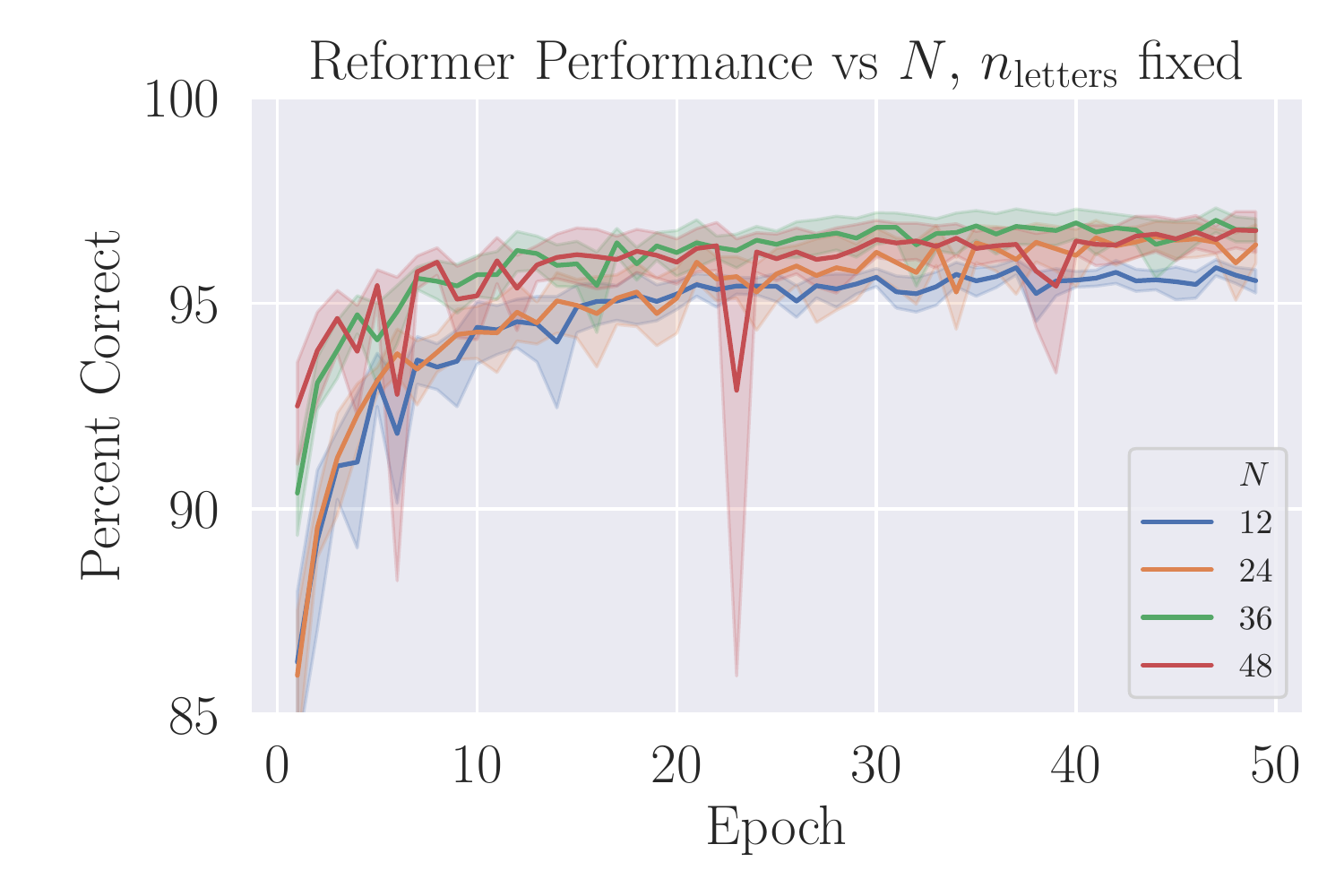}\label{fig:reformer_v_N_fixedNletters}}
\caption{Overview of the performance of the Reformer models for the UNKNOT decision problem. Shaded regions are confidence intervals associated to hyperparameters listed in Table \ref{tab:decision_parameters} that are not displayed on the plots.}
\label{fig:decision}
\end{figure}

\begin{figure}[t]
    \centering
            \includegraphics[width=.48\textwidth]{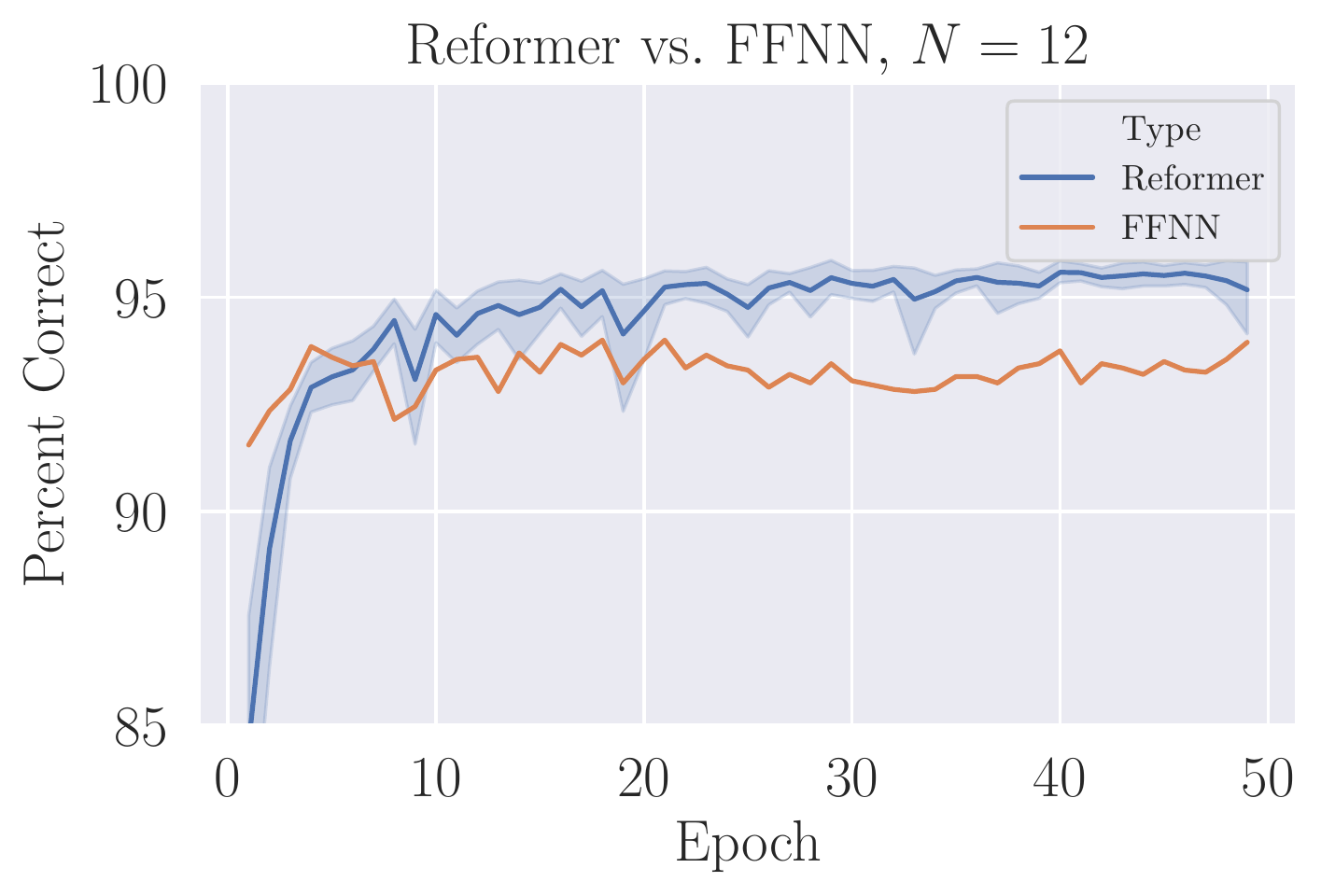}
            \includegraphics[width=.48\textwidth]{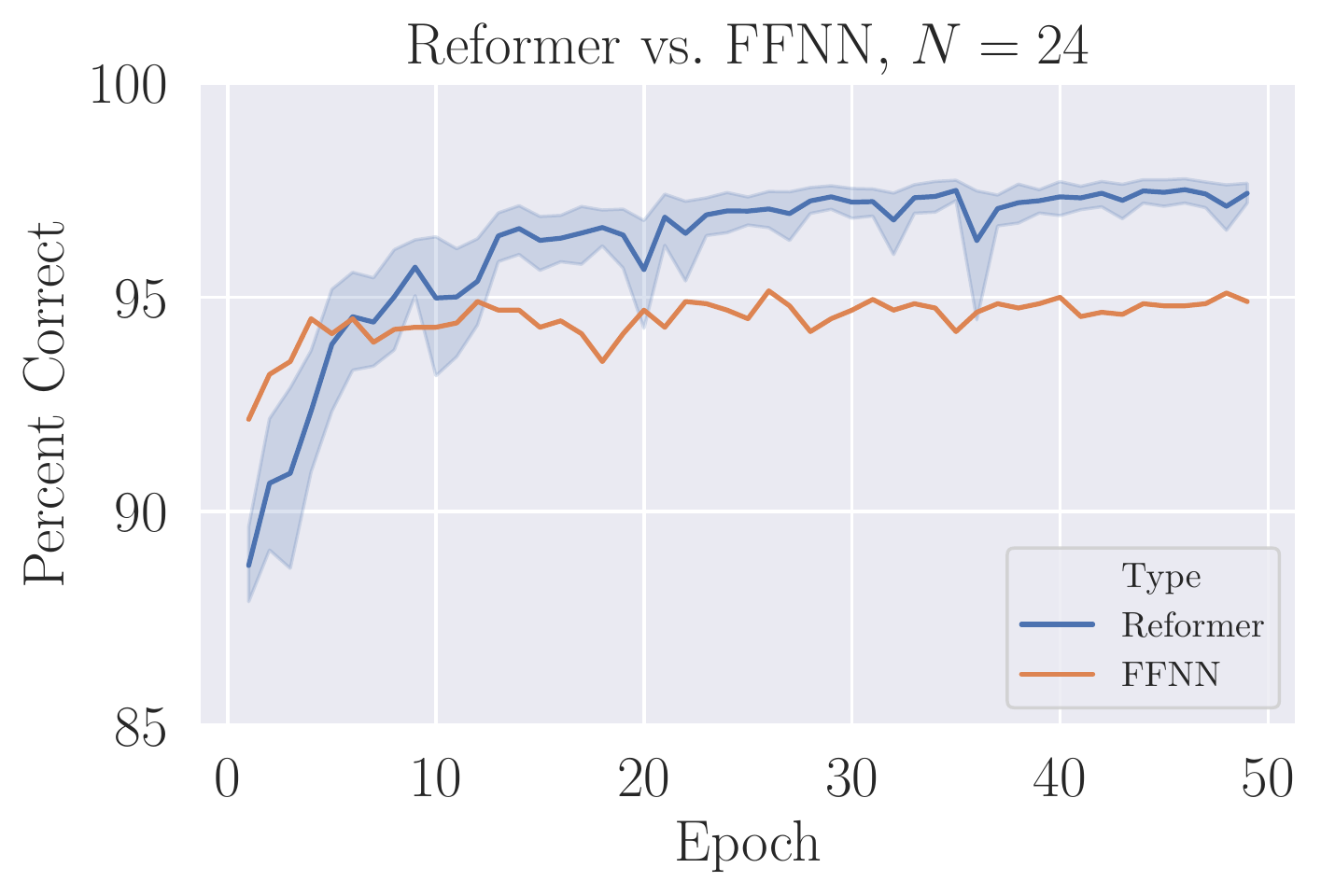}
       \\
            \includegraphics[width=.48\textwidth]{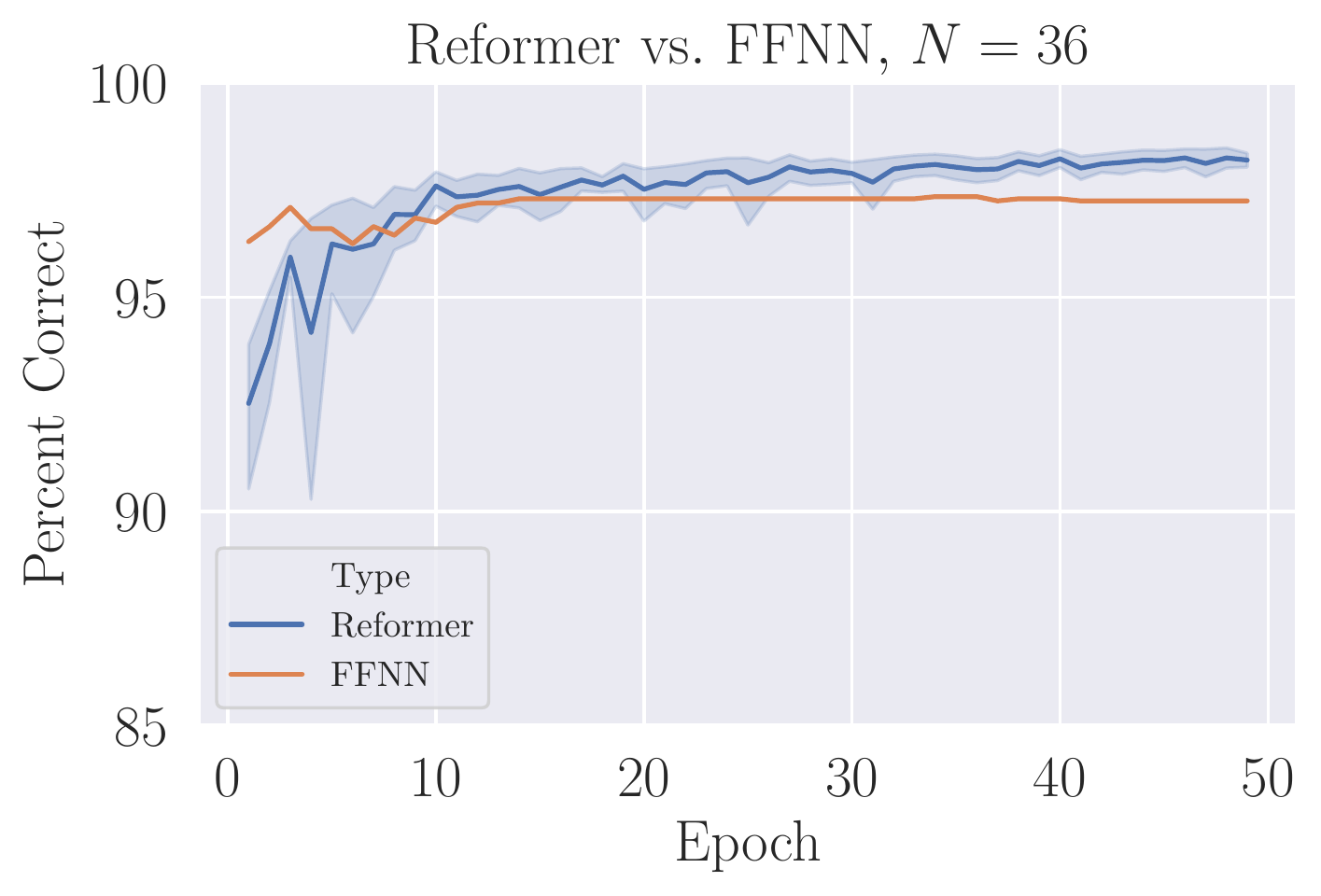}
            \includegraphics[width=.48\textwidth]{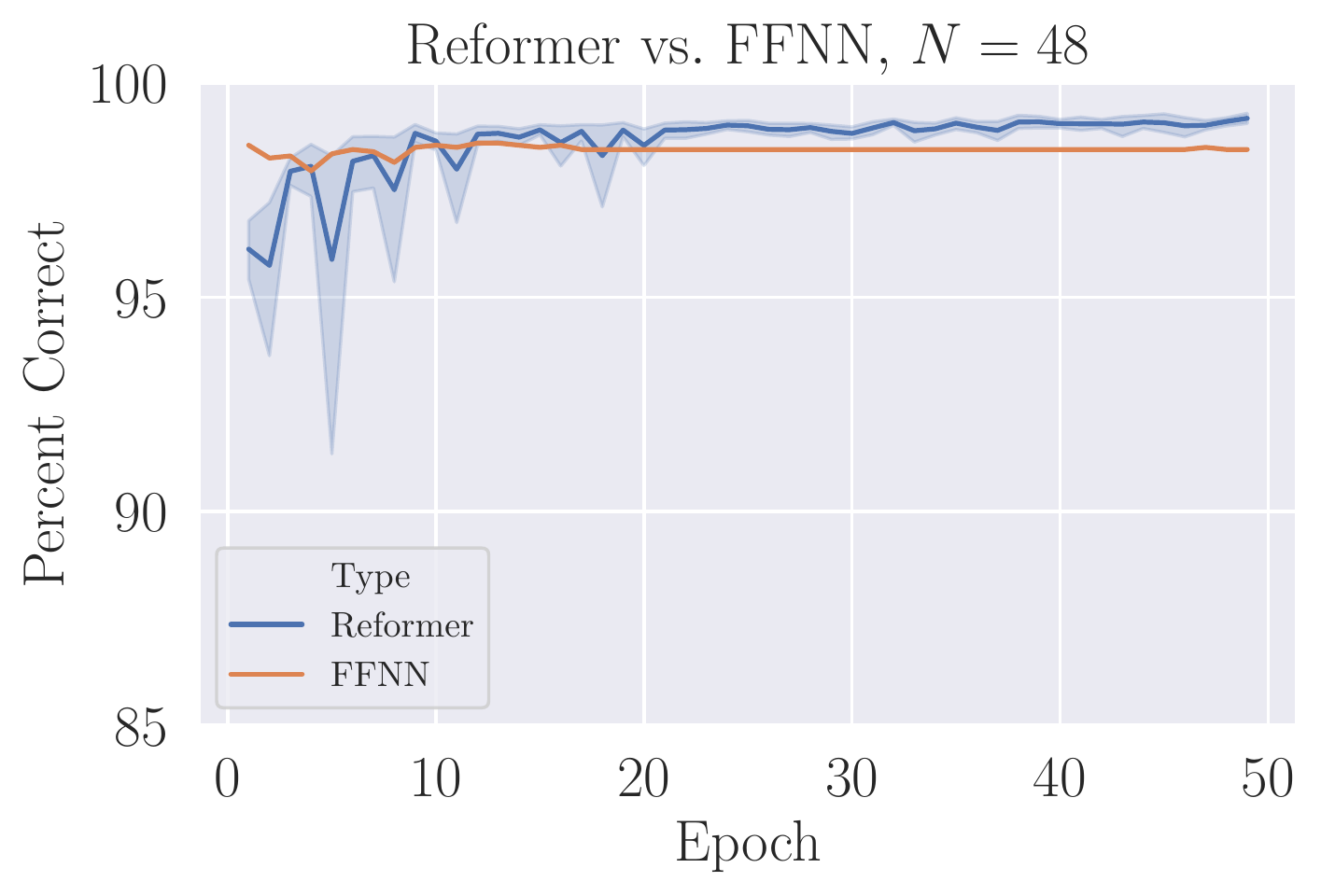}
    \caption{ Performance comparison between Reformer and feedforward network for different braid lengths $N$. Shaded regions are confidence intervals associated to hyperparameters listed in Table \ref{tab:decision_parameters} that are not displayed on the plots.}
    \label{fig:reformer_v_FFNN_fn_of_N}
\end{figure}

The dependence of the performance on $n_\text{hashes}$ or the use of full attention is plotted in Figure~\ref{fig:reformer_v_nhashes}. We find that a single hash performs slightly worse than four hashes, but that four hashes already performs similarly to full attention.

In Figure~\ref{fig:reformer_v_FFNN}, we see that the Reformers and shared-QK Transformer do outperform the feedforward network by a few percent, but the advantage is not as big as one might have expected. This is our first evidence that the UNKNOT problem is not particularly difficult for machine learning; we will see more in Section \ref{sec:unkotting}.  Interestingly, from Figure \ref{fig:reformer_v_FFNN_fn_of_N} we see that the performance gap between the NLP architectures and the feedforward networks is biggest at \emph{small} braid length $N$, and the performance gap is almost negligible at large $N$. Second, quite strikingly we see that performance increases with $N$, which is not what one would naively expect.

To try to better understand the latter point, we ran a slightly different set of experiments. Rather than fixing the total number of braid words to be $10,000$ for each $N\in \{12,24,36,48\}$, we fixed the total number of letters trained on for each $N$. That is, for $N=12$ we take $10,000$ braids, but for $N=48$ we take $2,500$ (and similarly for N=24 and N=36), so that the networks see the same total number of braid letters (generators) for each of the different values of $N$. All the parameters remain as in Table \ref{tab:decision_parameters}, and the performance is plotted in Figure~\ref{fig:reformer_v_N_fixedNletters}. We see that the performance still does increase with $N$, but it is much less drastic than before. While we know that NNs benefit from a larger training set, i.e.\ from seeing more distinct braid words as a whole, this suggests that after a certain threshold, the NN knows what to look for in the global structure, while still benefitting from being exposed to more subpatterns within braid words. It would be interesting to investigate this further by looking at the attention modules of the reformer to see which parts of the braid the NN actually pays attention to at which stage of the decision process, but this is beyond the scope of the current paper.

The first fact, namely that the performance gap gets smaller for larger $N$, is hard to disentangle from the point we have just discussed: Both the feedforward NN and the reformer benefit from the effect discussed above. Hence, as the braid words get longer and the performance increases, the performance gap has to shrink as both networks approach 100\% accuracy.

\subsection{Confident Predictions, Hard Knots, and the Jones Polynomial}

We now study network confidence and its correlation with the Jones polynomial, focusing on the simplest case $N=12$ because computing the Jones polynomial is \#P-hard \cite{jaeger_vertigan_welsh_1990}. Specifically, we train a Reformer on non-trivial knots as well as unknots drawn from the braid priors of Section \ref{sec:generating_data}, using the parameters summarized in Table~\ref{tab:decision_parameters_confidence}.

\begin{table}[t]
	\centering
    \begin{tabular}{|cc|}
    	\hline
        Parameter & Values \\ \hline
        $n_\text{hashes}$ & $ 4$\\
        %Full Attention & \{True, False\} \\ 
        Causal & False \\
        Bucket Size & $6$ \\
        Depth & $10$ \\
        Embedding Dimension & $250$ \\
        Epochs & 250 \\ 
        Optimizer & RMSProp \\
        Learning Rate & $.0001$ \\ \hline
    \end{tabular} 
    \caption{Parameters for the Reformer used to study network confidence.\label{tab:decision_parameters_confidence}}
\end{table}

\begin{figure}[t]
    \centering
        \makebox[\textwidth][c]{
            \includegraphics[width=.48\textwidth]{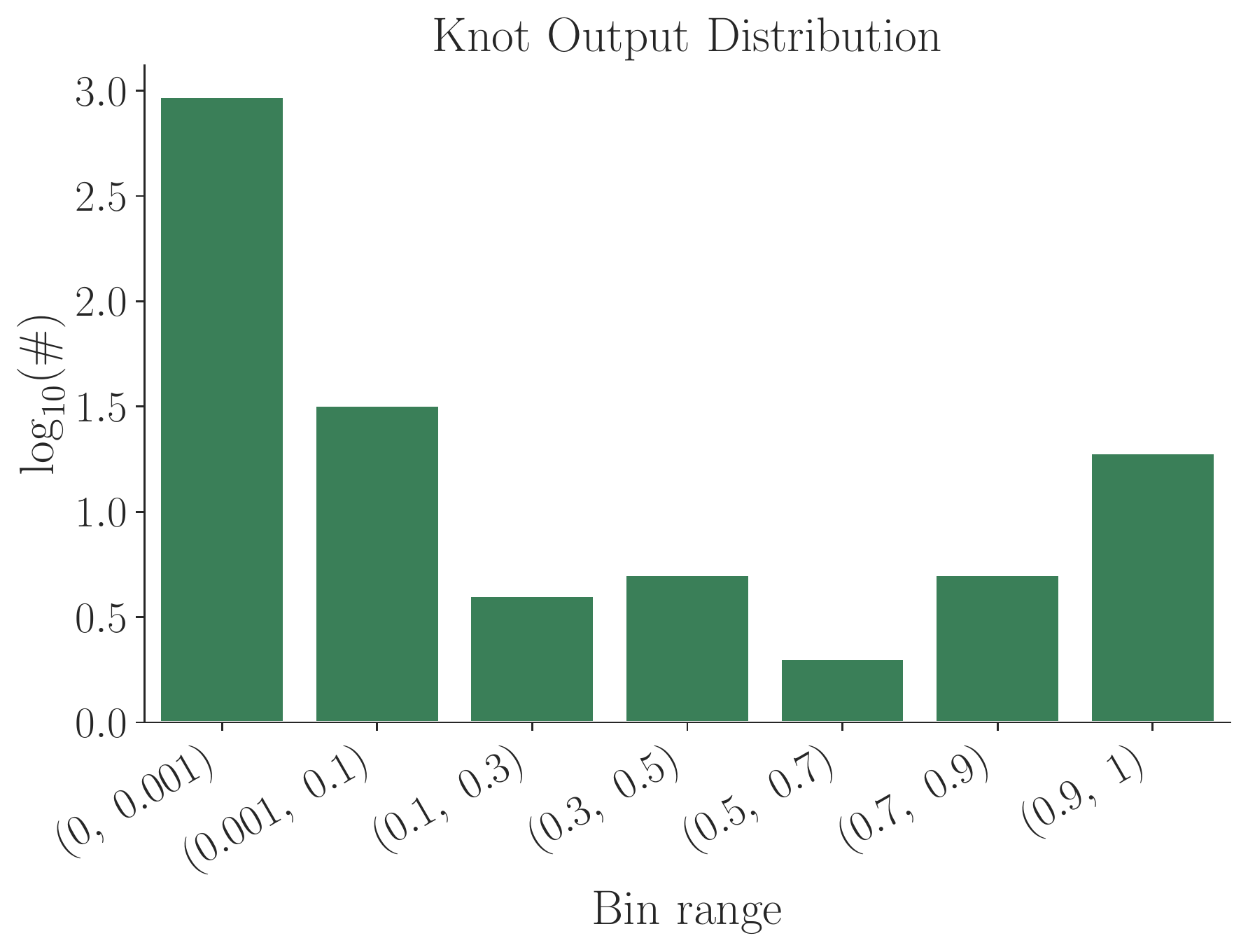}
            \includegraphics[width=.48\textwidth]{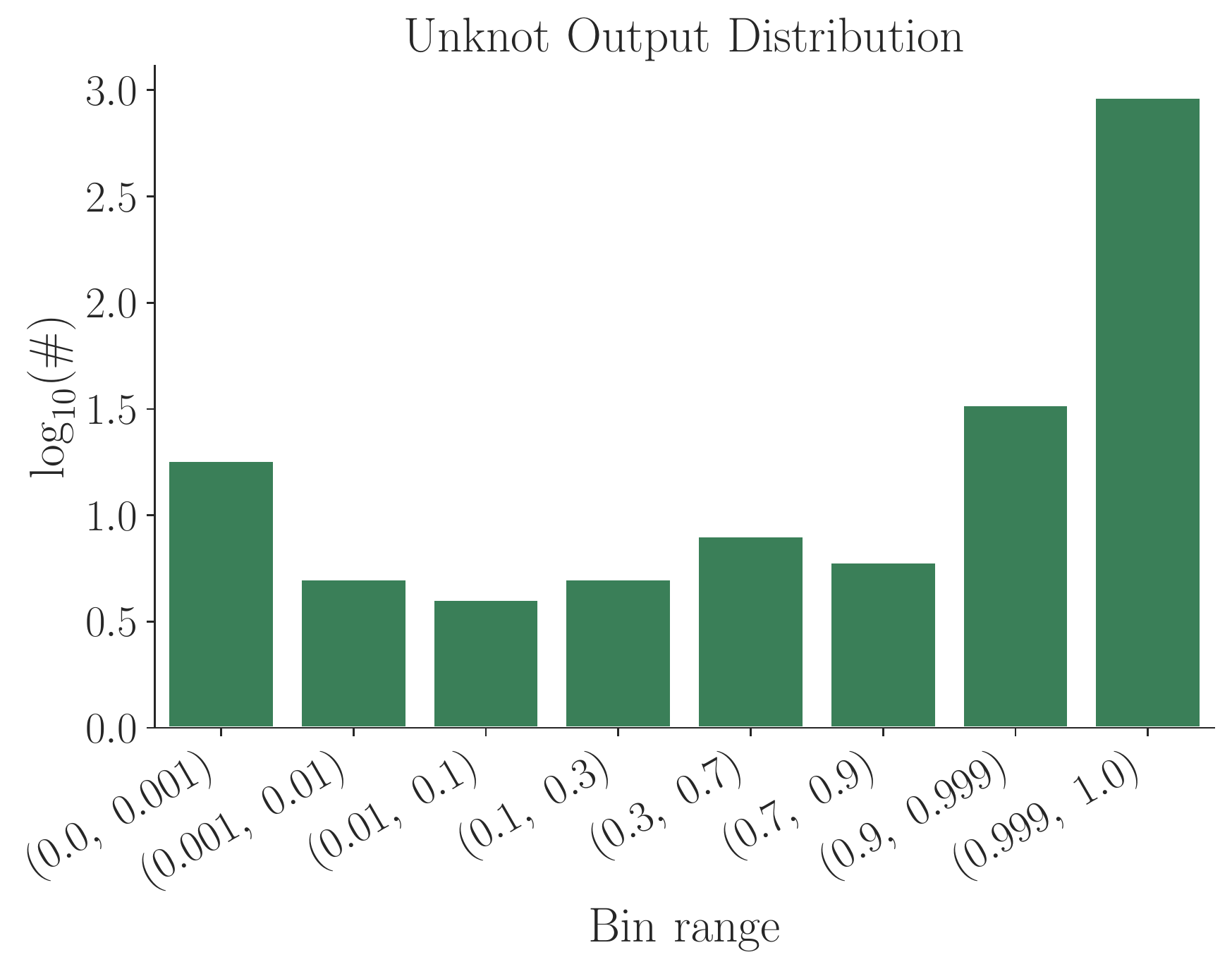}
        }
        \makebox[\textwidth][c]{
            \includegraphics[width=.48\textwidth]{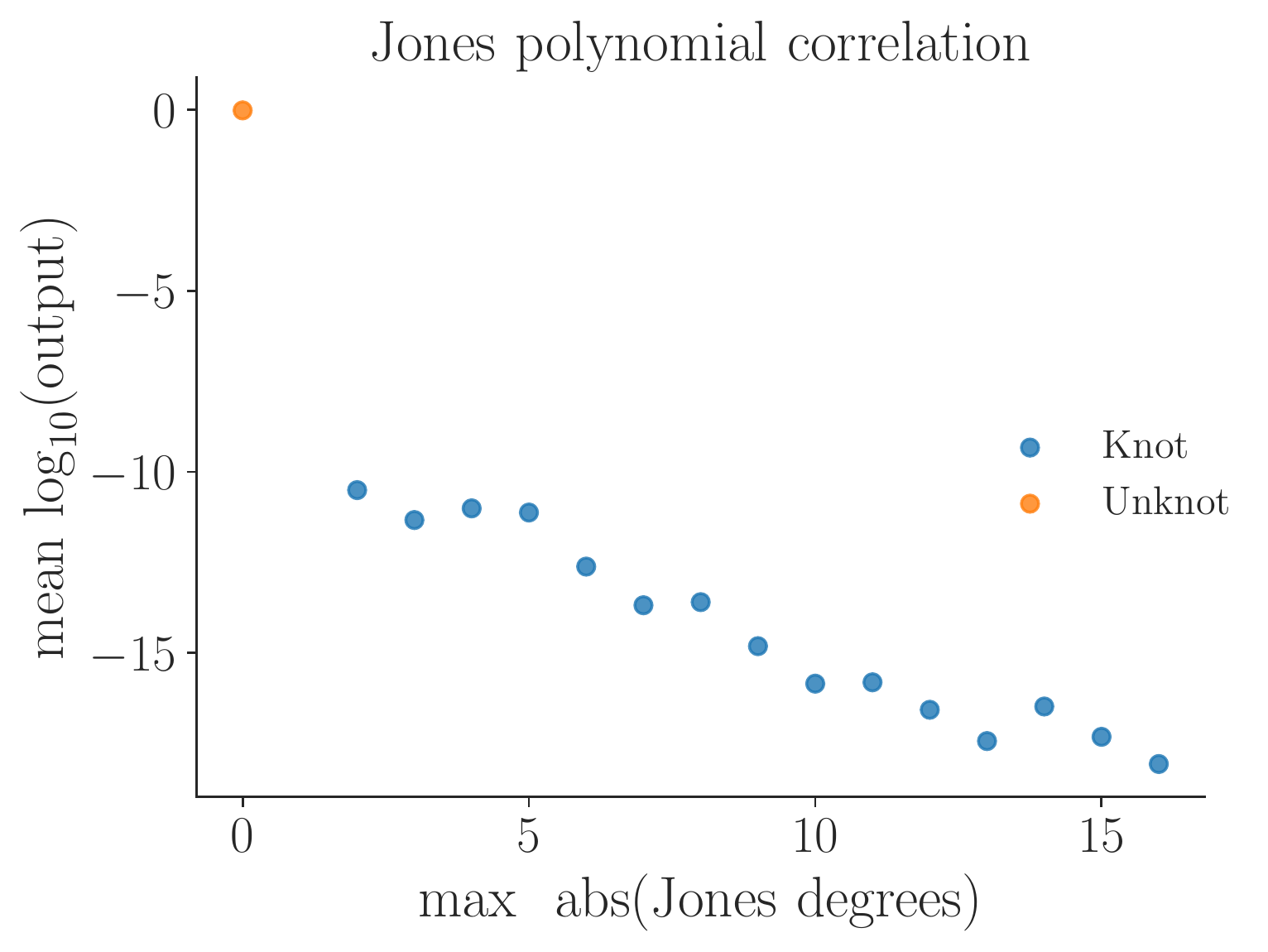}
        }
    \caption{A study of Reformer outputs for $N=12$ knots. \emph{Top Left:} Output distribution for knots. \emph{Top Right:} Output distribution for unknots. 
    \emph{Bottom:}  Correlation between network outputs and the maximum absolute value of Jones polynomial degrees.}
    \label{fig:N12_output_study}
\end{figure}

We test the trained network on $1,000$ non-trivial knots and unknots from a test set. Results are presented in Figure~\ref{fig:N12_output_study}. Since we label non-trivial knots and unknots as $0$ and $1$, respectively, these should be the locations of peaks in the output distributions of well-trained networks. The experimental distributions presented in Figure~\ref{fig:N12_output_study} show this correlation. Indeed, the networks are performing well, with precision >95\% (c.f.\ blue curve in Figure~\ref{fig:reformer_v_N}). We would like to point out the following central observations about the NN predictions for both unknots and non-trivial knots:
\begin{itemize}
\item \textbf{Very high confidence.} Over $900$ of the non-trivial knots (unknots) have outputs within $10^{-3}$ of their target value $0$ ($1$). This shows that for over $90\%$ of the non-trivial knots and unknots, the network is very confident in its prediction.
\item \textbf{Hard knots revealed by small peaks on the wrong label.} For both the non-trivial knot and unknot distributions, we see small peaks at the wrong end of the spectrum, at $1$ for non-trivial knots and $0$ for unknots. While this looks almost negligible on the plots because the peak at the correct values is very large, when one restricts the output distributions to the case of the network making \emph{wrong} predictions, most of the wrong predictions for non-trivial knots (unknots) occur in the bin closest to $1$ ($0$). For these non-trivial knots or unknots, it is not that the network is unsure of its prediction; rather, it is quite sure, but the prediction is wrong. We checked that this was not simply a function of initialization by running the same experiment with $10$ different random initalizations. We found that $22$ of the $1000$ non-trivial test set knots have output $>.95$ for all $10$ runs, and similarly $22$ of the $1000$ test set unknots have output $<.05$ for all $10$ runs. It is therefore natural to conjecture that these examples are fundamentally hard for the NN, i.e.\ they possess some adversarial property that makes the NN predict the wrong answer with high confidence. 

Since knots with $9$ or fewer crossings may be identified by their Jones polynomials, we ran experiments for $N=9$ braids to identify the hard braids. Specifically, we ran five Reformers, each with $4$ hashes, embedding dimension $250$, and depth $10$, on $5000$ $N=9$ braids from our prior with an $80/20$ test-train split. We trained each for $50$ epochs and kept the best model, and accuracy was $\sim 93\%$ on the test set for all five runs. However, of the $1000$ braids (with non-trivial knot closures) in the test set, $30$ of them had outputs above $0.9$ for all five experiments. These hard knots are
\begin{align}
    (3_1,19,242),\,\,\, \,
    (4_1,3,37), \,\,\, \,
    (6_2,1,34), \,\,\, \,
    (6_3,1,30), \,\,\, \,
    (8_{13},1,5), \,\,\,\, 
    (8_{20},4,19),
\end{align}
where the first entry of the tuple is the knot number in the Rolfsen table, the second (third) is the number of ``hard'' (total) instances of the knot in the test set. For knots that are not self-mirror, mirrors are included in the counts. We find that there are no hard knots with $9$ crossings, and that $\sim 2/3$ of the hard knots are trefoils, despite the fact that only $\sim 1/3$ of the knots in the test set are trefoils. While more statistics are necessary to draw a firm conclusion, this preliminary analysis seems to suggest that knots with fewer crossings are more likely to be hard.
A possible explanation for why of all things the simplest non-trivial knot, i.e.\ the trefoil, seems to confuse the networks could be that at fixed length $N$, knots with a smaller number of crossings in the minimal representation (three for the trefoil) contain more crossings that can be undone or disentangled. This is of course also the case for the unknot, where all crossings can be removed. This similarity might cause the NN to being tricked and ``overlooking'' that there is still a non-trivial component left in the braid representation of the trefoil knot after removing all superfluous crossings.

\item \textbf{Network uncertainty.} The network is uncertain when its output is around $0.5$. To set a more precise threshold, let us say that the network is uncertain if $0.3 < \text{output} < 0.7$. As described in the first point, there are very few knots for which a given NN is uncertain to begin with. Moreover, the uncertain knots are not robust to initialization: across the $10$ random initializations just mentioned, we find that there are no unknots or non-trivial knots for which all of the networks are uncertain.
\end{itemize}

We also study correlations between network outputs and the Jones polynomial. Specifically, for all knots in the ensemble we compute the Jones polynomial using \texttt{SageMath} \cite{sagemath} and compute the maximum of the absolute value of the degrees of the monomials. Since the Jones polynomial of the unknot is just the constant polynomial 1 (in the normalization of \texttt{SageMath}), we think of this as serving as a measure of the complexity of the Jones polynomial and hence of the non-triviality of the knot. In the bottom plot of Figure \ref{fig:N12_output_study}, we stratify the knots according to this measure of the degree and plot it against the mean of the network outputs. There is clearly a direct correspondence: on average, the higher the complexity of the Jones polynomial, the more confident the network is that the knot is, in fact, non-trivial. We emphasize that this correlation was learned by the network and \emph{not} put in by hand: the Jones polynomial did not enter anywhere into the training process.

\subsection{Going Up to Go Down: Hard Knots in Dowker-Thistlethwaite Notation}
\label{sec:hardKnots}
Interestingly, some unknots admit a precise notion of ``hardness'':\footnote{We emphasize that this likely has nothing to do with the notion of hardness discussed in the previous section.} there is no sequence of Reidemeister moves that simplifies them without increasing the number of crossings at some point. We provide an example where this phenomenon is illustrated for a simple braid representation of a knot around equation~\eqref{eq:hardBraid} in Section~\ref{sec:RL-environment}, but at this point we do not have a systematic way of constructing examples of such braids. However, 176 examples of such knots with 15 crossings have been constructed in the Dowker-Thistlethwaite representation in~\cite{Sikora1,Sikora2}. It is an interesting question how this property is related to the representation of the knot, which is, however, beyond the scope of this paper. In particular, we do not know whether this property is preserved when moving from the Dowker-Thistlethwaite representation to the braid representation of a knot.Translating Dowker-Thistlethwaite to braids is not straight-forward and might require moving some of the strands, thereby destroying the hardness property.

In order to study these ``hard'' knots nevertheless, we therefore have to use the hard Dowker-Thistlethwaite knots of~\cite{Sikora1,Sikora2}. However, creating a training set in Dowker-Thistlethwaite notation would require new methods of generating random non-trivial knots as well as unknots in this notation. In order to use our existing methods, we generate knots and unknots in braid notation and translate them into Dowker-Thistlethwaite notation. This direction is, in contrast to the other direction, straight-forward. 

With this data set, we trained a Reformer and simple feedforward network on the Dowker-Thistlethwaite representations of length $15$ drawn from our prior on braids. The Reformer was acausal, with bucket size $8$ and $4$ hashes, depth $10$, with $5$ heads, and embedding dimension $250$. Evaluated on a test set of $1000$ non-trivial knots and unknots, the trained Reformer had accuracy $90\%$ on knots and $87.7\%$ on unknots, and the feedforward network performed slightly worse. Moreover, these results are $\sim 5\%$ lower than Reformers or feedforward networks trained on the knots represented as braids, suggesting that braids may be a better representation of knots for machine learning than Dowker-Thistlethwaite notation. 

The reason for doing the analysis on knots in Dowker-Thistlethwaite notation, however, was to test it on the ``hard'' Dowker-Thistlethwaite knots, as just defined. Using the trained Reformer to make predictions for the $176$ known hard Dowker-Thistlethwaite knots from~\cite{}, none of which were explicitly used in the training set, the Reformer achieved only $2.2\%$ accuracy --- it was almost always wrong! Concretely, the mean and standard deviation of outputs for these unknots were $.03$ and $.14$, respectively, demonstrating that the network is quite sure of its wrong predictions. While this performance could be blamed on the out-of-sample prediction of hard unknots, we also observed that including a fraction of the hard unknots does not improve prediction accuracy for the other hard unknots. This gives a second reason, aside from the definition above, to think that these unknots are fundamentally hard.

\section{Unknotting with Reinforcement Learning\label{sec:unkotting}}
Instead of just using the Reformer (or FFNN) as a black box whose output is the likelihood with which the NN thinks the input is or is not the unknot, we also studied unknotting via reinforcement learning. We will discuss in detail the environment, states, actions, and reward we implement for this RL task. For an introduction to RL and these concepts we refer the reader to~\cite{Ruehle:2020jrk,Halverson:2019tkf}. 

\subsection{The RL environment}
\label{sec:RL-environment}
\subsubsection*{State space}
The states of the reinforcement agent are all braid words of a given length $2N_\text{max}$ for starting braids of length $N_\text{max}$ (the reason for the factor of 2 will become apparent once we discuss the actions below). Since the start state will be a braid whose closure is a knot, rather than a link, and this property is preserved under our actions, we are only considering braids with single-component closures. For such braids, a braid word of length $N$ has at most $2n$ generators and is thus an element of (at most) Br$_{n+1}$. Therefore, there are (at most)
\begin{align}
\label{eq:numStatesRL}
N_\text{states} = 1+\sum_{n=1}^{2N_\text{max}} (2n)^n
\end{align}
possible states. Out of these, the only terminal state is the state corresponding to the empty braid word in Br$_1$. Here, the subtlety with the destabilization move discussed around~\eqref{eq:destabilization}: Even if we start with a braid word $w\in\text{Br}_{N_\text{max}}$ corresponding to a one-component knot, after reducing its length to some $N<N_\text{max}$ the braid is taken to be an element of $\text{Br}_{N}$ rather than $\text{Br}_{N_\text{max}}$ to preserve the single-component property.\footnote{Of course, our algorithm can be run on multi-component links by first identifying all components (which is easy) and then running the algorithm on each component individually. This will simplify each component as much as possible (but won't contain information on e.g.\ the linking number).}

\subsubsection*{Reward function}
The purpose of using RL is that we want to find an equivalent braid representation of any input knot with as short a braid word as possible. The reason for this is two-fold:
\begin{itemize}
\item When we want to use NNs to analyze braids that represent knots (not necessarily just for the unknot question addressed in this paper), the (input dimension of the) NN can be smaller if the input word is shorter.
\item Since the unknot is represented by an empty braid word in Br$_1$, we can use this to detect whether a knot is the unknot.
\end{itemize}
This makes it natural to use the negative length of the current braid word as a reward (or rather punishment): The fact that shorter braid words are punished less strongly means that the agent will attempt to minimize the length of the braid word. Moreover, since each move receives a punishment, the agent is incentivized to reduce the length of the braid word as fast as possible. As discussed in Section~\ref{sec:hardKnots} and exemplified around~\eqref{eq:hardBraid} below, there do exist knots for which the braid word has to become longer before it can be simplified. It is hence important that the agents maximize their (long-term) return rather than just their (short-term) reward. As we shall see next, the action space necessarily contains illegal actions for some states. We punish such illegal actions with a negative reward of $4N_\text{max}$.

Let us make one further comment on the length of the braid word: there is a fast algorithm due to Dehornoy~\cite{Dehornoy:1997aaa} that solves the word problem for braids. He defines a notion of a reduced braid, which is unique within each equivalence class of braids. This means that two braids with the same reduced braid word can be transformed into each other using the braid relations~\eqref{braid-relations}. Since the reduced braid word for the empty braid is the empty braid word, this gives a sufficient criterion for any braid to describe the unknot: if its reduced braid word is empty, then there exist a sequence of braid relations that will turn the knot (given as a braid) into the unknot. We find that this sufficient criterion is extremely weak for our unknots. Essentially none of the braids representing the braid word have an empty reduced form. Since we use Markov moves together with braid relations to generate the unknots, and since Markov moves do change the braid, this result is not too surprising. But it begs the question of whether a better measure for triviality of the braid would be the reduced braid word rather than just the length of the braid word. Our experiments clearly show that this is not the case; if we base the rewards on the reduced braid word obtained from Dehornoy's algorithm, the agent learns slower and performs worse. This might be due to the case that the reduced word length changes more erratically when braid-altering actions are performed than the non-reduced braid word length. This erratic change in reward might make it harder for the agent to learn.

\subsubsection*{Action space}
Remember that performing the Markov moves depicted in Figure~\ref{fig-markov} change the braid but not the braid closure, i.e., not the knot. However, just using Markov moves does not guarantee that the unknot which corresponds to the empty braid word can be reached. Indeed, any given braid configuration of the unknot will be reducible to the empty braid word by using Markov moves together with the braid relations~\eqref{braid-relations}. As a simple example, consider the braid word
\begin{align}
w=[1,2,1,-2]\,.
\end{align}
Just Markov moves alone will not simplify this to an empty braid word corresponding to the unknot. However, if we use the first braid relation, this braid can be seen to be equivalent to the braid
\begin{align}
w'=[2,1,2,-2]\,.
\end{align}
Now, removing the consecutive inverses $[2,-2]$ will lead to the braid word $[2,1]$,\footnote{Note that we could equally well have removed the generator 2 at the fist position together with the generator -2 at the last position, since we are interested in the closure of the braid. This would have left us with the braid word $[1,2]$, which describes an inequivalent braid, but the same knot.} which, after two destabilization moves, collapses to the empty braid word.

This means that we need in principle 4 different kinds of actions for the agent: The two types of braid relations~\eqref{braid-relations}, as well as the two Markov moves. Let us count the number of possible actions: In order to carry out the two braid relations, we need to specify the position in the braid at which the relations are to be used, adding $N_\text{max}$ actions each. Markov move 1, i.e.\ conjugation, consists of two consecutive actions: multiplying by a generator and its inverse on either side of the braid word, and subsequently simplifying the braid. This adds 2 (composite) actions. Markov move 2, if allowed, either adds or removes a strand from the braid, adding another 2 actions. 

However, there is another subtlety here: If we perform conjugation as a composite action of adding consecutive inverses on the closure of the braid and removing a different set of consecutive inverse operators, we prescribe an order in which these operations are to be carried out. Fixing this order means that some braids of unknots cannot be simplified to the empty braid word anymore. As an example, consider $w=[-1,2,1,-2]$. This cannot be simplified by applying conjugation (i.e.\ cyclic shifts), braid relations, and destabilization moves. However, if we (i) insert consecutive inverses at the second to last position, (ii) use braid relation 1 on positions 2-4, (iii) remove consecutive inverses at positions 1 and 2, (iv) conjugate by $-2$ on the left and 2 on the right and remove consecutive inverses at the ends of the braid, and (v) performing two destabilization moves, the braid word $w$ becomes the empty braid word:
\begin{align}
\label{eq:hardBraid}
\begin{split}
w=[-1,2,1,-2] &\xrightarrow{\text{(i)}} [-1,2,1,2,-2,-2] \xrightarrow{\text{(ii)}} [-1,1,2,1,-2,-2] \xrightarrow{\text{(iii)}} [2,1,-2,-2]\\
& \xrightarrow{\text{(iv)}} [1,-2]\xrightarrow{\text{(v)}} \emptyset
\end{split}
\end{align}
This is an example of the previously mentioned fact that for some knots the length of the braid word (i.e.\ the number of crossings) has to be increased before it can be decreased and illustrates that we need to allow the insertion of consecutive inverse operators at some point in the braid, without specifying a priori the order in which they are removed, or whether one applies other operations, such as using the first braid relation and does not remove the generators at all. Naively, this requires adding the possibility of inserting $2N_\text{max}$ generators $\sigma_i^{\pm1}\sigma_i^{\mp1}$ at any of the $N_\text{max}$ positions of the braid, thus adding $2N_\text{max}^2$ actions.

In total, this set of actions would have 
\begin{align}
N_\text{actions}=2N_\text{max}^2+2N_\text{max}+4 
\end{align}
actions, which can be several thousands. In our experience, RL works particularly well for small action spaces, while the state space can be very large.\footnote{For small state spaces, the problem can often be brute-forced and RL is not necessary.} In the case at hand, the action spaces can become quite sizable. Moreover, they contain many illegal actions, since most braid relations can only be applied at very few positions in the knot. This large number of illegal actions means that the braid is often not changed by an action, which in turn requires a very long exploration phase by the agent in order to realize which actions are valid based on which input states. To counter this, we do not use the set of actions described above but rather introduce a different set of high-level actions that does not grow as fast with $N_\text{max}$.

We have tried several agents with different types of composite actions, but we will only discuss the one that works best here. First, we add cyclic shifts to the left and right in the braid word, thus including Markov moves of type 1. This way, we have covered the Markov moves of type 1 with 2 actions. Next, we need to address equivalences of the braid. Here, we can potentially save a huge number of actions, if we find a more efficient way to describe the insertion of consecutive inverses (which introduced $2N_\text{max}^2$ actions) and for the braid relations (which introduced $2N_\text{max}$ actions).

For the consecutive inverses, we can get away with only $N_\text{max}$ actions in the following way: Since conjugation, i.e.\ cyclic shifts, are already included in the actions, it does not matter where the inverse operators are inserted, so we insert them at the beginning and the end of the braid word. This reduces the number of actions by a factor of $N_\text{max}$, since we need not specify the position anymore. Theoretically, this reduction comes at the cost of performing up to $\lceil N_\text{max}/2\rceil$ cyclic shifts after adding the inverse pair of generators in order to move them to any desired position in the braid word. In practice, we only care about the braid closure anyways, so this shift will never be necessary. Hence, it seems as if we need to add $2N_\text{max}$ actions, $N_\text{max}$ that send $w\to\sigma_i w \sigma_i^{-1}$ and another $N_\text{max}$ that send $w\to\sigma_i^{-1} w \sigma_i$. However, it is enough to include only the first $N_\text{max}$: The second $N_\text{max}$ can be obtained from the former by performing a cyclic right-shift, commuting the pair of inverse operators, and performing a cyclic left-shift,
\begin{align}
\sigma_i w \sigma_i^{-1} \to \sigma_i^{-1}\sigma_i w \to \sigma_i\sigma_i^{-1} w \to \sigma_i^{-1} w\sigma_i\,.
\end{align}
We have thus reduced the number of actions coming from insertion of consecutive inverses from $2N_\text{max}^2$ to $N_\text{max}$. However, note that inserting consecutive inverses increases the length of the braid word by 2 (this is also true for the original $2N_\text{max}^2$ actions). Since we use NNs to approximate the state and action value functions, and since we need to fix the input dimension of these NNs, we allow a maximum intermediate length of $2N_\text{max}$ for a braid word with original length $N_\text{max}$.

In addition, we bundle several simplifications (i.e.\ removing trivial link components, relabeling braid generators, performing destabilization moves, and removing inverses of a pair of operators that are either consecutive or separated by other braid generators with which the generator and its inverse commute) into one action called S\textsc{mart}C\textsc{ollapse}. These operations are repeated until the braid does not simplify further. Note that we did not add a stabilization move. While we do not know whether there are situations where one needs to perform a stabilization in order to eventually get to the empty braid word, we observe empirically that the performance of the agents went slightly down when adding this extra action. Of course, longer training or better hyperparameter tuning, plus maybe allowing for longer intermediate braid words, should overcome this drop in performance. However, from the drop in performance we do not see evidence that, for the knots in our database, stabilization contributes significantly to simplifying braid words of unknots to the empty braid word. If in doubt, one could of course add this extra action at the cost of a few percent accuracy or finding better hyperparameters.

Next, let us address reducing the $2N_\text{max}$ actions coming from the braid relations. This can be achieved as follows:
\begin{enumerate}
\item Start at the beginning of the braid word
\item Use a braid relation at the first possibility
\item Use conjugations (i.e.\ cyclic shifts) to move the position where the braid relation has been applied to the end of the braid word
\end{enumerate}
This reduces the number of actions from $2N_\text{max}$ to 2. The price to pay is that if the network wants to perform a braid relation at a specific position, it might have to perform, in the worst case, $N_\text{max}$ additional cyclic shifts in order to move the position where the braid relation is to be applied far enough to the left such that it is the first occurrence. 

Note that one could be tempted to not do the cyclic shift in the third step, but instead remember where the last action was performed and start from that position the next time. However, this would break the Markov property (since the next action on a state would then depend on how the state was reached) and this would make the problem not (necessarily) amendable to be solved with RL.

To summarize, this leaves us with
\begin{align}
N_\text{actions}=5+N_\text{max} 
\end{align}
actions on a braid word $[i_1,i_2,\ldots,i_k]=w\in\text{Br}_{n+1}$, $i_j\in[-n,n]$:
\begin{itemize}
\item \textsc{SmartCollapse} (see Algorithm~\ref{alg:SmartCollapse}): removes twists, performs destabilization, removes inverses.
\item shift left (conjugation + remove inverses):\[w\xrightarrow{a_{n+1}}(-i_1)\circ w \circ (i_1) = [i_2,\ldots,i_k,i_1]\].
\item shift right (conjugation + remove inverses):\[w\xrightarrow{a_{n+2}} (i_k) \circ w \circ (-i_k) = [i_k,i_1,i_2,\ldots,i_{k-1}]\].
\item braid relation 1 and shift right: let $m$ be the position where the braid relation can be applied and $s=[i_{m+1},i_{m+2},i_{m+3}]$ be the three-letter substring to which it is applied, yielding $s'$. Then:\[[i_1,i_2,\ldots,i_m,s,i_{m+4},\ldots, i_k]\xrightarrow{a_{n+4}}  [i_{m+4},\ldots, i_k,i_1,i_2,\ldots,i_m,s']\].
\item braid relation 2 and shift right: let $m$ be the position where the braid relation can be applied and $s=[i_{m+1},i_{m+2}]$ be the two-letter substring to which it is applied, yielding $s'$. Then:\[[i_1,i_2,\ldots,i_m,s,i_{m+3},\ldots, i_k]\xrightarrow{a_{n+5}}  [i_{m+3},\ldots, i_k,i_1,i_2,\ldots,i_m,s']\].
\item $n$ Markov moves of type 1, i.e.\ conjugations by an arbitrary generator $i_m\in[-n,n]$, do not remove inverses:\[w\xrightarrow{a_m}(i_m)\circ w \circ (-i_m) = [i_m,i_1,i_2,\ldots,i_k,-i_m]\].
\end{itemize}

Note that the only action that simplifies (i.e.\ reduces the length of) a given braid word is \smartcollapse. All other actions serve the purpose of applying braid relations and Markov moves that will eventually allow \smartcollapse to simplify the braid word. Also note that now the only illegal actions are the ones that try to add a pair of inverse generators to a braid word that is already of length $2N_\text{max}$. In order to get to such a braid word, the agent has to perform $N_\text{max}$ of such actions before encountering an illegal move for the first time. We do not know whether there exist knots whose braids require increasing the length (i.e. the number of crossings) by more than a factor of two at some intermediate step, before the knot can be simplified and collapsed to the unknot. However, especially for larger $N_\text{max}$, this will incur quite a sizable punishment over a large number of steps, which makes it unlikely that the agent would follow this policy unless we choose a discount factor of 1 (or very close to 1). While the \smartcollapse action might not (be able to) change the braid, we do not consider this an illegal action.

\subsection{The RL algorithm}
Many RL algorithms have been developed to solve an MDP. We have tried:
\begin{itemize}
\item A3C (Asynchronous Advantage Actor-Critic)~\cite{Mnih:2016A3C} and the synchronous version A2C, with feedforward and reformer NNs.
\item DQN (Deep Q-Networks)~\cite{Mnih:2013DQN}, with and without dueling, with feedforward and recursive (LSTM) NNs.
\item PPO (Proximal Policy Optimization)~\cite{Schulman:2017PPO}, with feedforward NNs.
\item TRPO (Trust-Region Policy Optimization)~\cite{Schulman:2015TRPO} with GAE (Generalized Advantage Estimation), with feedforward and reformer NNs.
\end{itemize}
We used different libraries based on Tensorflow/Keras~\cite{Plappert2016:kerasrl}, Pytorch~\cite{Kostrikov:PytorchA3C}, and ChainerRL~\cite{Fujita:2019chainerrl}. We performed thorough (but not excessively extensive) box searches for the hyperparameters and tried varying the NN architecture. Based on these experiments we found that PPO and DQN performed worse. Moreover, DQN showed oscillatory behavior in the reward (and the actual performance). A2C and A3C performed approximately on the same level (with A3C slightly outperforming A2C), and both were much better than DQN and PPO. By far the best algorithm was TRPO. We tried using conventional FFNN as well as the Reformer NN with TRPO. While the reformer seemed to perform slightly better, it trains much slower. Since the results with a classic FFNN were already very good, we ended up using the latter in a ChainerRL implementation of the algorithm with the hyperparameters of~\cite{Henderson:2017reinforcement}.\footnote{TRPO promises to require little tuning of hyperparameters. For us, the parameters of~\cite{Henderson:2017reinforcement} (which are standard in ChainerRL) worked well and we did not attempt to tune them further.} We train the agents for $5\times10^6$ steps on a standard GPU, which takes around $24$ hours.

Since we ended up using TRPO, let us briefly explain the ideas behind this algorithm. As is the case for example in A3C, TRPO uses a NN to approximate the state value and action value functions. However, while A3C uses gradient descent to minimize the mean square error for the loss function and uses this as a baseline for the policy updates (which are also performed using gradients), TRPO follows a different approach. The problems TRPO tries to address are the following: First, since gradient-based algorithms are by definition insensitive to curvature, this can lead to problems for strongly curved loss landscapes. Second, the step size should be a function of the curvature, making smaller updates in strongly curved regions and larger updates in flatter regions. Third, TRPO can guarantee policy updates that improve the policy and will eventually lead to the optimal policy. 

To address the first point, TRPO uses standard SGD for the value function updates, but a second-order update (based on conjugate gradients and a line search) for the policy function updates. The second point is addressed by introducing trust regions. The algorithm determines, for a given size of the trust region, the maximum step size that it wants to explore and then finds the optimal point within this trust region. Concerning the third point, instead of optimizing the policy function, TRPO optimizes a surrogate function (based on the KL divergence between the old and the new policy) that approximates the expected reward (computed from the current policy) locally. The authors of~\cite{Schulman:2015TRPO} show that if this surrogate function bounds the expected reward from below, this is guaranteed to lead to an improvement.

\subsection{Results}
In order to evaluate the performance of the RL agents, we run them on braid words that represent unknots. The reason is that for these we know that there exists a series of moves that reduce the braid word to the empty braid word. For non-trivial knots, we do not know what the simplest braid word is (this is the whole point of using ML for this task), hence we cannot judge how well the algorithm performs on non-trivial knots. 

The agents receive knots of length $N_\text{max}\in\{12,24,36,48,72,96\}$. Since, as discussed above, inserting a pair of inverse generators can increase the length of the braid, we introduce a maximum length of $2N_\text{max}$ as input size for the NN and apply zero-padding. Note that zero is not a valid generator in our convention; since we use $\pm i$ to encode the generator $\sigma_i^{\pm1}$, we start at $i=1$. For that reason, 0 is a good pad value. We define accuracy as the fraction of unknots for which the RL agents found a series of moves that reduced the input braid word to the empty braid word. We also benchmark the trained agents against a random walker (RW), which does not follow any sophisticated policy but draws the next action from the set of all actions with a flat prior. We set an upper bound of 500 actions and check performance on the same 10,000 unknots for all algorithms. 

\begin{figure}[t]
    \centering
    \includegraphics[width=.99\textwidth]{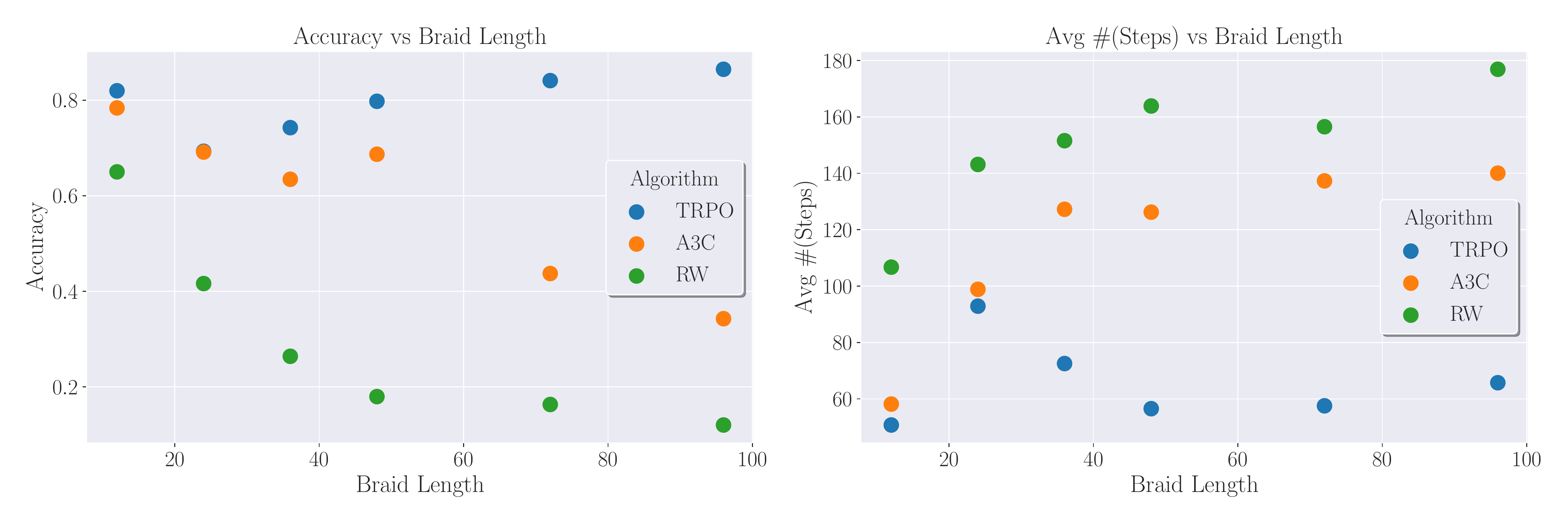}
    \caption{Performance comparison of the TRPO, A3C and RW algorithms. Left: Fraction of unknots whose braid words could be reduced to the empty braid word as a function of initial braid word length $N$. Right: Average number of actions necessary to reduce the input braid word to the empty braid word as a function of $N$.}
    \label{fig:RLResults}
\end{figure}

We present the results of the RL runs in Figure~\ref{fig:RLResults}. To keep the plot manageable, we only show the two best-performing RL agents (TRPO and A3C), both with FFNNs, as well as the performance of the random walker. Let us start by discussing the plot on the left in Figure~\ref{fig:RLResults}. For TRPO, we find that the fraction of braids encoding the unknot that can be fully reduced to the empty braid word is around 85\% for all unknots, more or less irrespective of their length and number of generators (within the interval we tested). While results for A3C are comparable for smaller braids (up to $N_\text{max}=48$), the accuracy drops significantly for $N_\text{max}=72$ and $N_\text{max}=96$ from around 70\% to 34\%. Since there seems to be no fundamental obstruction and A3C is working well for smaller $N_\text{max}$, we expect that more complex NNs, longer exploration phases, a discount factor closer to $1$, and possibly tweaks to other hyperparameters could lead to better performance also at larger $N_\text{max}$. However, we did not attempt this since TRPO already performed at a constant, higher rate. Both RL runs are to be contrasted with the performance of the random walker. This shows a sharp drop in accuracies, from around 64\% accuracy for $N_\text{max}=12$ to around 10\% for $N_\text{max}=96$. This illustrates that the strategy (policy) learned by the agents significantly outperforms brute-force searches, especially for larger $N_\text{max}$.

We also plot the average number of actions the agent takes until obtaining the empty braid word on the right in Figure~\ref{fig:RLResults}. For $N_\text{max}=12$, the trained TRPO agent takes an average of around 50 steps to fully reduce a braid word of the unknot to the empty word. In contrast, the random walker needs twice as many steps until it stumbles upon a solution. For this low $N_\text{max}$, the A3C agent still performs very similarly to the TRPO agent. For larger $N_\text{max}$, the average number of actions needed by the A3C agent and the RW increases. However, the average number of steps needed by the RW is by a factor of 1.3 to 1.8 larger as compared to the A3C agent. Interestingly, the number of steps is almost constant for the TRPO agent and a factor of 2 or more below the random walker. Note that this is only averaging the number of actions over those knots that could actually be reduced to the unknot, which is around 8,800 for TRPO, 3,400 for A3C, and 1,000 for RW at $N=96$. The results suggests that if the number of steps would be bound to be around 60 instead of 500, the accuracy of the random walker would already be very low even for $N_\text{max}=12$. Since the computational complexity of the actions (especially of S\textsc{mart}C\textsc{ollapse}) grows with the length of the braid word, the smaller the number of actions, the faster the algorithm performs. Conversely, this means that a brute-force approach to the unknotting problem becomes unfeasible for larger $N_\text{max}$.

\subsection{Actions taken to unknot}

\begin{figure}[t]
    \centering
    \includegraphics[width=.9\textwidth]{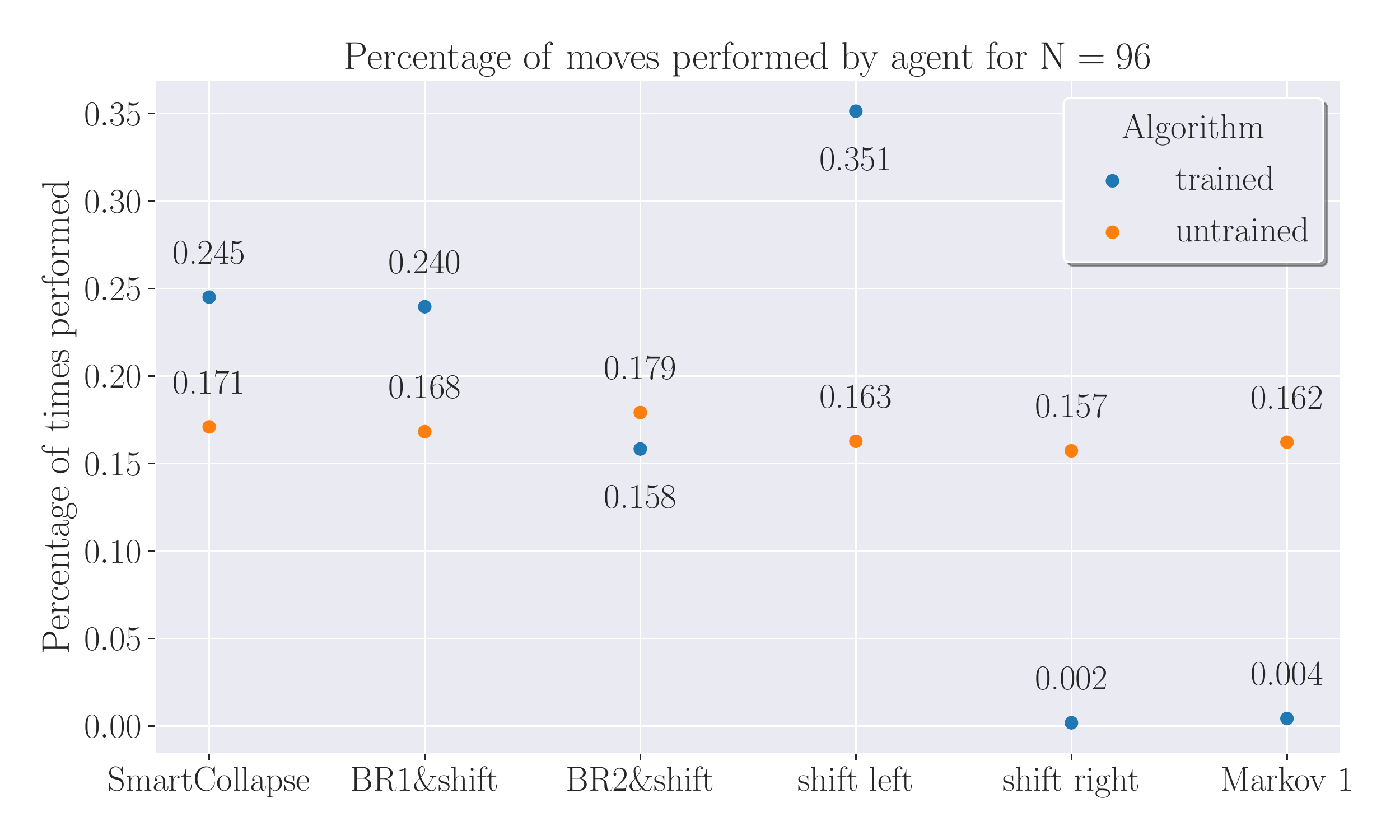}
    \caption{Comparison of the moves that are picked by an untrained agent vs a trained agent.}
    \label{fig:RLmoves}
\end{figure}

In reinforcement learning, the learning process of the AI agent is represented by the flow of the state-dependent distribution on action space, i.e., a change in policy.
Therefore, as a simple attempt in trying to understand what the agent has learned, we look at the distribution of actions that the agent performs before and after training, as summarized in Figure~\ref{fig:RLmoves}. We look at the largest, most complicated braids at $N=96$. The actions are those summarized at the end of Section~\ref{sec:RL-environment}. In the plot, we have counted all $n$ actions that perform Markov move 1 in one category labeled ``Markov 1'' and divided by $N$. 

The first thing to observe is that the untrained agent performs every actions equally often. This is to be expected, since it is just choosing actions randomly from a flat prior. The fact that around 10\% of the unknots represented by braid words of length 96 can be reduced to the empty braid word by just performing actions drawn from our base set of actions with a flat prior (see Figure~\ref{fig:RLResults}) sets a baseline against which we should benchmark the performance.

For the trained agent, we see that shift left and \smartcollapse are the most frequent actions. This makes sense. First, \smartcollapse is the only action that can actually reduce the length of the braid word (and thus decrease the punishment); all other actions either leave the length unchanged or increases it by two. Second, the large asymmetry between shifting left and right (the former is the most frequent action, the latter the least frequent action and essentially negligible) is owed to the fact that applying braid relations 1 and 2 are followed by many right shifts, which move the part that the braid relation has acted upon to the end of the braid word. Hence, the agent gets many right shifts for free, which it might have to counter with shifting left in order to rearrange the braid such that the part where the agent wants to apply a braid relation to next is the first possible occurrence. 

Next, we observe that braid relation 1 is used around 10\% points more often than braid relation 2. This is explained by the fact that \smartcollapse will remove non-consecutive inverses if the generators between a generator and its inverse commute with the generator, in other words $\smartcollapse(w \sigma_i w' \sigma_i^{-1}) =ww'$ if $\sigma_i$ commutes with all generators in $w'$. Hence, the agent gets again many commutation moves for free in smart collapse and hence needs to perform these less frequently.

Finally, we find that Markov move 1 is performed rather infrequently. This is interesting and tells us that the braid relations seem far more important than Markov move $1$ when it comes to simplifying our braids. Note that these increase the braid length and with it the punishment, and hence the agent is also reluctant to perform these. The fact that the agent performs them nevertheless (and even more frequently than shift rights, which do not change the punishment) illustrates that the agent is indeed maximizing its long-term return and the fact that some braids need to become longer before they can be simplified as illustrated in an example in Section~\ref{sec:RL-environment}.

\section{Conclusion \label{sec:conclusion}}

In this paper we have proposed studying knot theory using techniques from natural language processing (NLP), as is natural since any knot may be represented by a braid word. The statement that any topologically equivalent knots are related by a sequence of Reidemeister moves gives rise to a natural action on the space of braid words, comprised of Markov moves and braid relations. From the NLP perspective, then, the problem of knot equivalence becomes a question of whether two words (more specifically, their closures) are equivalent in the language. Transformers and Reformers are natural architectures for such studies.

The entirety of the paper focused on a fundamental problem in knot theory, the UNKNOT problem, which asks whether a given knot is trivial or non-trivial. This is equivalent to the question of whether a given representation of knot may be continuously deformed to the simple circle without ripping or tearing, which requires the existence of a sequence of Reidemeister moves that performs the simplifications, or alternatively an equivalent sequence of moves on a braid representative. 
Much of our interest in this problem stems from its relation to the smooth four-dimensional Poincar\'e conjecture, a fundamental open problem in geometric topology. From the perspective of computational complexity, the UNKNOT problem resides in $\text{NP}\cap \text{co-NP}$, though if one bounds the number of Reidemeister moves from above it is NP-hard. Further information on both our motivation and complexity issues are discussed in the main text.

In Section \ref{sec:generating_data}, we gave a detailed description of the induced prior for our method of randomly generating braids (and consequently knots) based on drawing braid generators from a flat prior and performing a sequence of Markov moves (again chosen randomly with flat prior). Letting $N$ be the number of letters in the braidword, $N=9$ braids have knot closures that, when simplified, have $9$ or fewer crossings. Such knots happen to be uniquely identifiable by their Jones polynomials, which allowed us to compare the distribution of knots drawn from our prior to the flat prior on knots with $9$ or fewer crossings. The latter induces a distribution on the number of crossings that increases exponentially, whereas our prior induces a much flatter distribution.

In Section \ref{sec:decision}, we studied the UNKNOT decision problem, aiming at solving the binary classification of whether or not a given knot is the unknot. We performed systematic experiments using Reformers, Shared-QK Transformers, and feedforward neural networks, finding that the NLP architectures outperformed feedforward networks, but only by a few percent, perhaps suggesting that the UNKNOT problem is still relatively easy as an ML problem; accuracy above $90\%$ was easy to achieve in all cases, and for the NLP architectures accuracy in the mid- to high- $90$s was achieved. Interestingly, we saw the counter-intuitive result that performance clearly \emph{increased} with increasing $N$, where $N$ is the number of letters in the braid word. These experiments had a fixed number of braid words, which meant that the neural network saw a larger total number of letters for increasing $N$. Fixing instead the number of total letters seen by the networks, we found the increase in performance with increasing $N$ is less drastic. Still, this result surprises us, and a better conceptual understanding may provide a useful perspective on the UNKNOT problem. For $N=12$ knots, we also found that the certainty with which the networks predicted the non-triviality of a knot was directly correlated with the maximimum of the absolute value of the Jones polynomial degrees. We emphasize that this correlation was learned, and no information about the Jones polynomial was put in by hand. 

A notion of a ``hard'' knot also emerged from our analysis: in some cases, the network applies the \emph{wrong} knot-vs-unknot label to a knot, and it is quite sure of its wrong prediction. For instance, in our convention, non-trivial knots were labelled with a $0$, and in some rare cases the network would make a prediction of, e.g.,  $\sim .999$, on a non-trivial knot, indicating that the network is quite sure that it is an unknot. To test whether this was an accident, we repeated the experiments with ten randomly initialized neural networks and found that many of the incorrectly labelled non-trivial knots were incorrectly labelled in all of the experiments. This suggested that some knots may be fundamentally ``harder'' than others, or possess an adversarial property that tricks the NN into assigning the wrong label with high confidence. Since knots with $9$ or fewer crossings may be identified by their Jones polynomial, we were able to identify ``hard'' knots associated to $N=9$ braids. We found that no knots with $9$ crossings were hard, and that $\sim 2/3$ of the hard knots were trefoils, despite the fact that only $\sim 1/3$ of the knots in the test set were trefoils. This suggests that knots with less than $N$ crossings may be more likely to be hard knots. We speculate that a possible reason why small knots seem to be harder for the NN is because at fixed length $N$, these knots contain many superfluous crossings which can be removed, simplifying the knot tremendously. This might trick the network into thinking that the knot can be completely simplifed to the unknot, overlooking the non-trivial trefoil hiding within the mostly trivial braid word. In that case, it would be interesting whether this adversarial property is robust against changing the description of the knot from a braid closure to e.g.\ Dowker-Thistlethwaite, Gauss codes, etc. To study this a bit more, we looked at the NN performance when using the Dowker-Thistlethwaite notation. Overall, we found that the performance of the networks were slightly worse than for those trained with braids. Moreover, we used the 176 unknots of~\cite{Sikora1,Sikora2}, which require increasing the number of crossings before one is able to reduce the knot to the trivial unknot with no crossings. We find that indeed knots with this property are hard for the NN, and the NNs tend to consistently misidentify such knots as non-trivial.

In Section \ref{sec:unkotting}, we studied the unknotting problem using reinforcement learning. The idea is to train an agent that, given a knot (encoded as a braid), can find a sequence of moves that simplifies the braid as much as possible. We define simplicity via the length of the braid word. Since every unknot can be presented as the empty braid word in Br$_1$, this allows us to identify unknots using this agent: if the agent can find a sequence of actions that turns a braid word into the empty braid word, the given knot is provably the unknot. This has the advantage that the result can be verified for unknots, in contrast to the black-box models of section~\ref{sec:decision}, which can only assign a probability to a knot being the unknot.

While using braid relations and Markov moves would be enough to reduce any braid representing the unknot to the empty braid word, we use a different set of more high-level actions. The reason is that the number of actions can become quite large, which makes the agent difficult to train. Using trust region policy optimization (TRPO) and the base set of actions, we find that our trained agent can identify a set of moves that reduces a starting braid word to the empty braid word in over 80 percent of the cases even for knots with up to 96 crossings (i.e.\ length of the braid word). This beats the next best RL algorithm (A3C) by more than a factor of 2, and a brute-force algorithm that performs random actions by a factor of almost 8.

Interestingly, we find that the number of actions necessary to obtain the empty braid word from a starting braid word that represents the unknot is constant over the range of braid words we consider (from length 12 to 96). Moreover, we observe that, for our randomly generated unknots, using braid relations is much more important to simplify the knot than using Markov move 1.

\bigskip

\noindent \textbf{Acknowledgments.} We thank Peter Battaglia, Kyle Cranmer, Michael Freedman, Mark Hughes, Ciprian Manolescu, Alex Radovic, Danilo Rezende, and Adam Sikora for useful discussions. The work of S.G.\ is supported by the U.S. Department of Energy, Office of Science, Office of High Energy Physics, under Award No.\ DE-SC0011632, and by the National Science Foundation under Grant No.\ NSF DMS 1664227. J.H.\ is supported by NSF CAREER grant PHY-1848089. The work of P.S. is supported by the TEAM programme of the Foundation for Polish Science co-financed by the European Union under the European Regional Development Fund (POIR.04.04.00-00-5C55/17-00).

\appendix

\section{Algorithms \label{sec:app_algos}}

\begin{algorithm}[H]
    \caption{\textsc{RandomMarkovMove}}
    \label{alg:RandomMarkovMove}
\begin{algorithmic}[]
    \State $i\sim \mathcal{U}(\{0,1\})$
    \If{$i=0$} \Comment{Conjugation Markov move}
        \State $j \sim \mathcal{U}(\{1,\dots,\text{max}(\text{abs}(B)\})$.
        \State $k \sim \mathcal{U}(\{0,1\})$.
        \State $B \gets [(-1)^{k} \, j] + B + [(-1)^{k+1} \, j]$
    \Else \Comment{New strand Markov move}
        \State $k \sim \mathcal{U}(\{0,1\})$.
        \State $B \gets B + [(-1)^{k} \, \text{max}(\text{abs}(B))]$
    \EndIf
    \State \Return $B$.
\end{algorithmic}
\end{algorithm}

\begin{algorithm}[H]
    \caption{\braidrelationone:}\label{alg:BraidRelation1}
\begin{algorithmic}[]
    \Require Braid $B$, int start, bool take\_closure.
    \State $i\gets$ start
    \While{$i$ mod length($B$)$\;\neq\;$ start - 1}
        \If{not takeClosure and $i$ $\stackrel?=$ length($B$)} \Comment{Reached end of Braid}
            \State \Return $B$
        \EndIf
        \State $[p1, p2, p3]$ = $[i, i+1, i+2]$ mod length($B$)
		\If{[B[$p1$], B[$p2$], B[$p3$]] $\stackrel?=$ $[\pm k,\pm(k+1),\pm(k)]$ for some $k\in\mathbb{Z}$}        
            \State [B[$p1$], B[$p2$], B[$p3$]] $\gets  [\pm (k+1),\pm(k),\pm(k+1)]$
         \ElsIf{[B[$p1$], B[$p2$], B[$p3$]] $\stackrel?=$ $[\pm(k+1),\pm(k),\pm(k+1)]$ for some $k\in\mathbb{Z}$}       
            \State [B[$p1$], B[$p2$], B[$p3$]] $\gets  [\pm (k),\pm(k+1),\pm(k)]$
        \EndIf
    \EndWhile
    \State \Return $B$.
\end{algorithmic}
\end{algorithm}

\begin{algorithm}[H]
    \caption{\textsc{SmartCollapse}: method to reduce braid length.}
    \label{alg:SmartCollapse}
\begin{algorithmic}[]
    \Require Braid $B$.
    \State New braid $B' \gets$ empty braid word.
    \While{$B'\neq B$ as braid words}
        \State $B' \gets B$.
        \State $B\gets$ \textsc{RemoveConsecutiveInverses}($B$).
        \State $B\gets$ \textsc{RemoveFreeStrands}($B$).
        \State $B\gets$ \textsc{Destabilize}($B$).
        \State $B\gets$ \textsc{RemoveNonconsecutiveInverses}($B$).
    \EndWhile
    \State \Return $B$.
\end{algorithmic}
\end{algorithm}

\begin{algorithm}[H]
    \caption{K\textsc{notify}: turn braid representative of a link into a knot. Iteratively weaves together two strands not in the same link component until the braid closure is a knot.}\label{alg:Knotify}
\begin{algorithmic}[]
    \Require Braid $B$.
    \State $CSL \gets$ component strand list; list of list of strands in each component of $B$. 
    \While{$B\neq[]$ and $|CSL|\neq 1$}
        \State strands $\gets \{1,\dots,$ max(abs($B$))$+1\}$
        \For{ $CS\in CSL$}
            \For{strand $\in CS$} 
                \State up, down $\gets$ strand$+1$, strand$-1$
            \State $i \sim \mathcal{U}(\{0,1\})$.
            \If{up $\notin$ component and up $\in$ strands}
                \State $B\gets B + [(-1)^i\, \text{strand}]$ \Comment{Weave together strand, strand $+1$.}
                \State \textbf{break} twice.
            \ElsIf{down $\notin$ component and down $\in$ strands}
                \State $B\gets B + [(-1)^i\, (\text{strand}-1)]$ \Comment{Weave together strand, strand $-1$.}
                \State \textbf{break} twice.
            \EndIf
            \EndFor
        \EndFor
        \State $CSL \gets$  component strand list; list of list of strands in each component of $B$. 
    \EndWhile
    \State \Return $B$.
\end{algorithmic}
\end{algorithm}

\begin{algorithm}[H]
    \caption{R\textsc{andom}U\textsc{nknot}: generate random unknot representative.}
    \label{alg:RandomUnknot}
\begin{algorithmic}[]
\Require $n_\text{letters}, M\in \mathbb{Z}$.
    \State Braid $B \gets$ empty braid word.
    \While{$|B| \neq n_\text{letters}$}
            \If{$|B| > n_\text{letters}$}
                \State $B \gets$ empty braid word.
            \EndIf
            \For{$k\in \{1,\dots,M\}$}
                \State $B$ $\gets$ \textsc{RandomMarkovMove}($B$).
                \If{$|B|-1\geq 0$}
                    \State $B\gets$ \textsc{BraidRelation}2($B$, start position$\sim\mathcal{U}(\{1,\dots,|B|\})$).
                \EndIf

            \EndFor
            \State $B \gets$ \textsc{SmartCollapse}($B$).
    \EndWhile
    \State \Return $B$.
\end{algorithmic}
\end{algorithm}

\begin{algorithm}[H]
    \caption{\textsc{RandomKnot}: generate random non-trivial knot representative.}
    \label{alg:RandomKnot}
\begin{algorithmic}[]
\Require $n_\text{letters}, n_\text{strands},M\in \mathbb{Z}$.
    \State Braid $B \gets$ empty braid word $[]$.
    \While{$|B| \neq n_\text{letters}$}
            \If{$|B| > n_\text{letters}$}
                \State $B \gets$ empty braid word.
            \EndIf
            \While{$|B| < n_\text{letters}$}
                \State $i \sim \mathcal{U}(\{0,1\})$.
                \State $j \sim \mathcal{U}(\{0,\dots, n_\text{strands}-1\})$.
                \State $B \gets B + [(-1)^i\, j]$
            \EndWhile
            \State $B\gets$ \textsc{Knotify}($B$)
            \If{$B\neq []$} \Comment{Knotify sometimes yields an empty word.}
                \For{$k\in \{1,\dots,M\}$}
                    \State $B$ $\gets$ \textsc{RandomMarkovMove}($B$).
                    \State $B\gets$ \textsc{BraidRelation}2($B$, start position$\sim\mathcal{U}(\{1,\dots,|B|\})$).
                \EndFor
                \State $B \gets$ \textsc{SmartCollapse}($B$).
            \EndIf
    \EndWhile
    \State \Return $B$.
\end{algorithmic}
\end{algorithm}

\clearpage
\section{Knot or not? A game for children. \label{sec:game}}

Every child needs to be introduced to low-dimensional topology at the earliest possible age. We have developed a game to
help you in your quest. 

In \emph{Knot or not?}, your child will develop aptitude for determining whether or not a given
knot diagram is topologically trivial, i.e., is it a non-trivial knot or the unknot? Mastery of the game will help your children contribute  to household
chores, including (but not limited to) increasing parental health by unknotting their headphones before a jog, or increasing a sibling's well-being by unknotting their hair without resorting to
the scissor trick.

\begin{figure}
    \makebox[\textwidth][c]{
        \includegraphics[scale=.4]{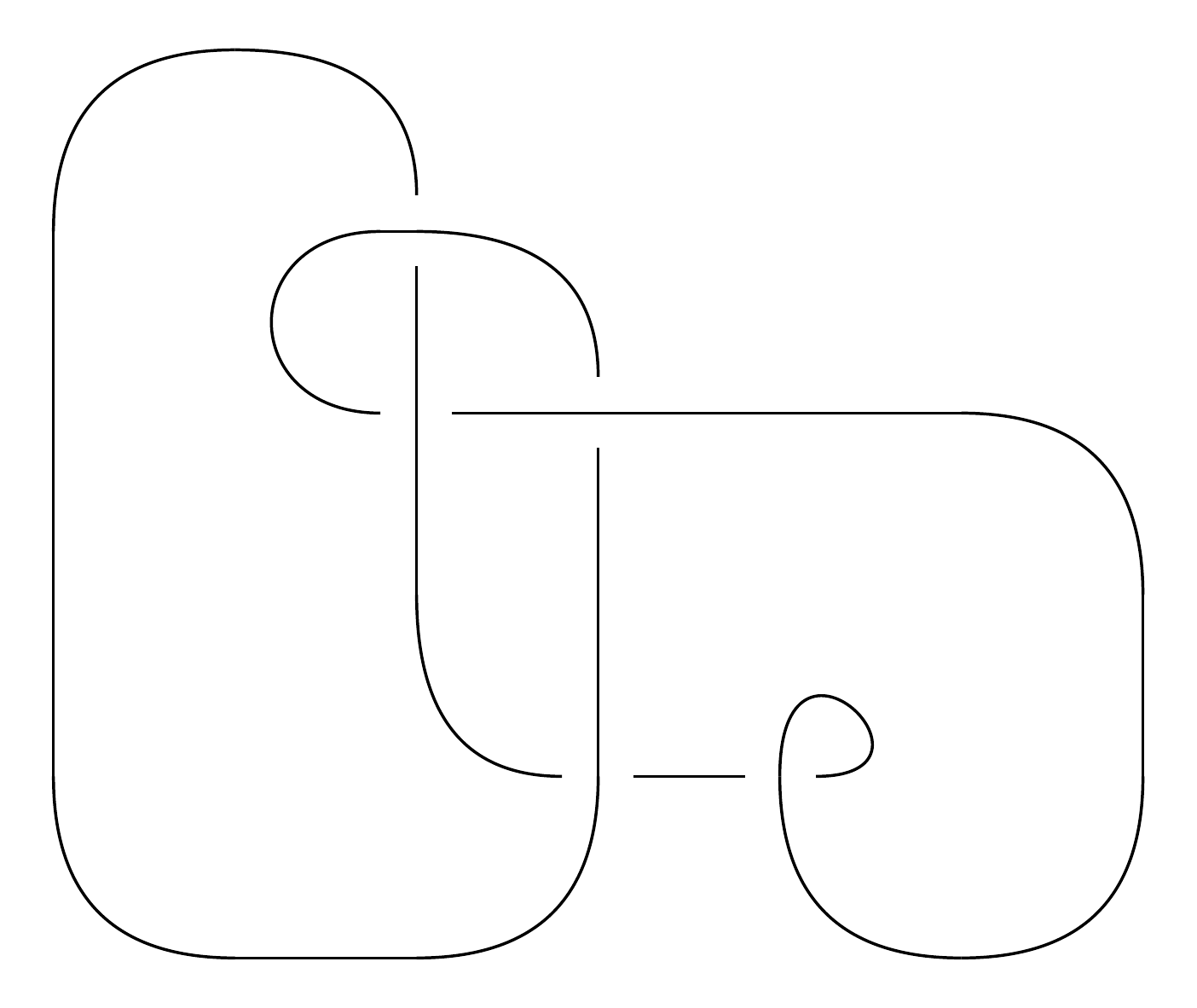}
        \includegraphics[scale=.4]{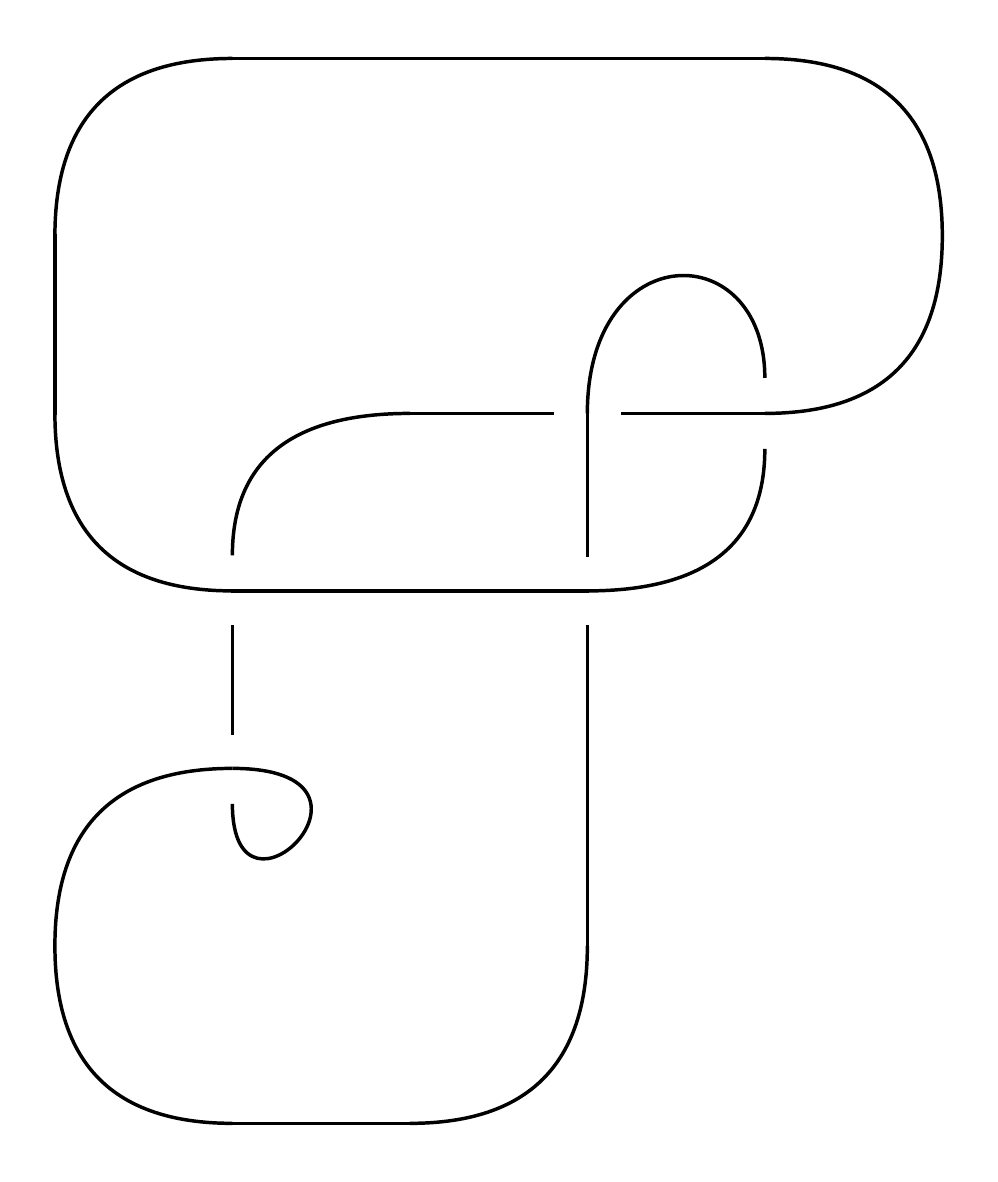}
        \includegraphics[scale=.4]{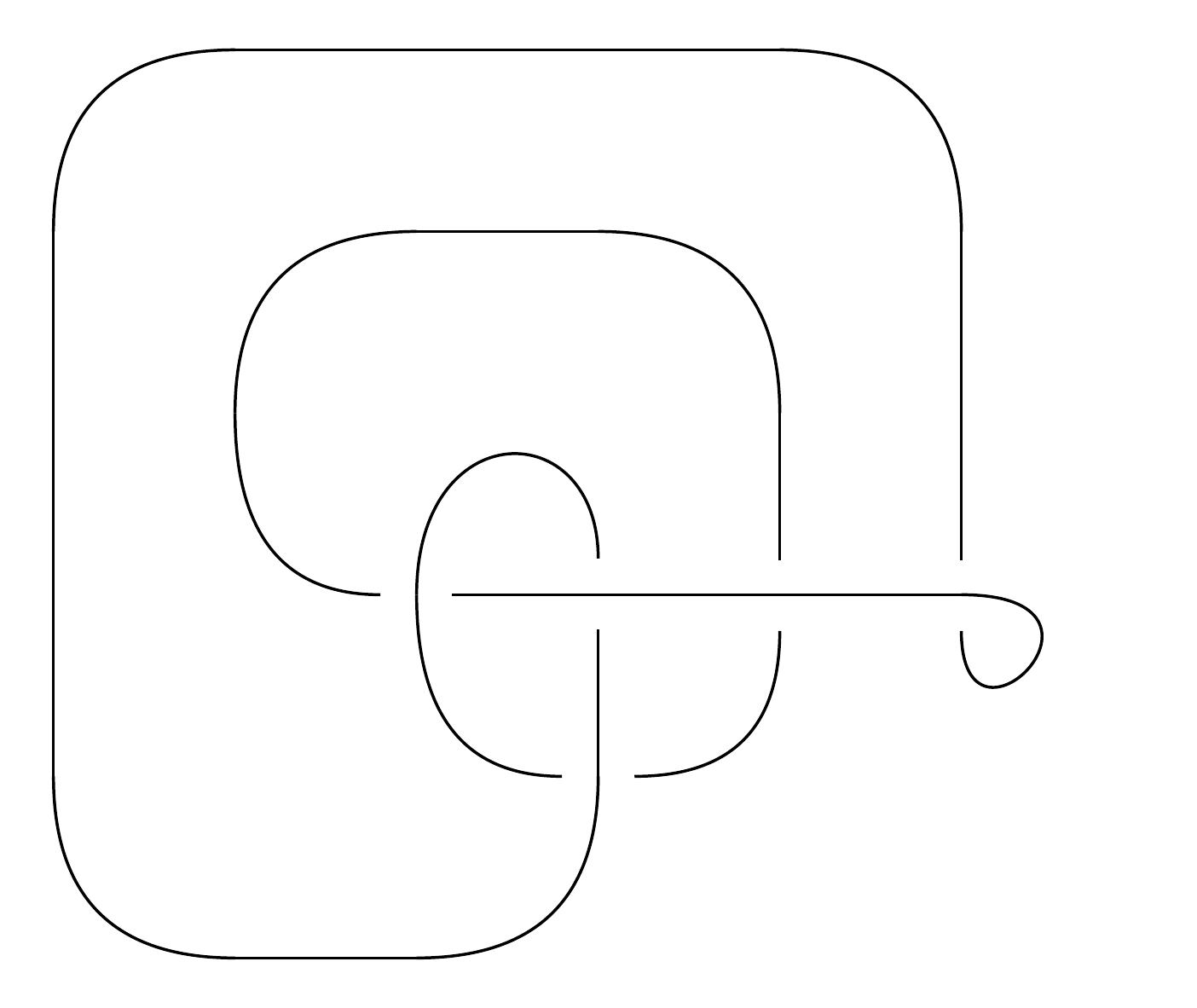}
    } 
    \makebox[\textwidth][c]{
        \includegraphics[scale=.4]{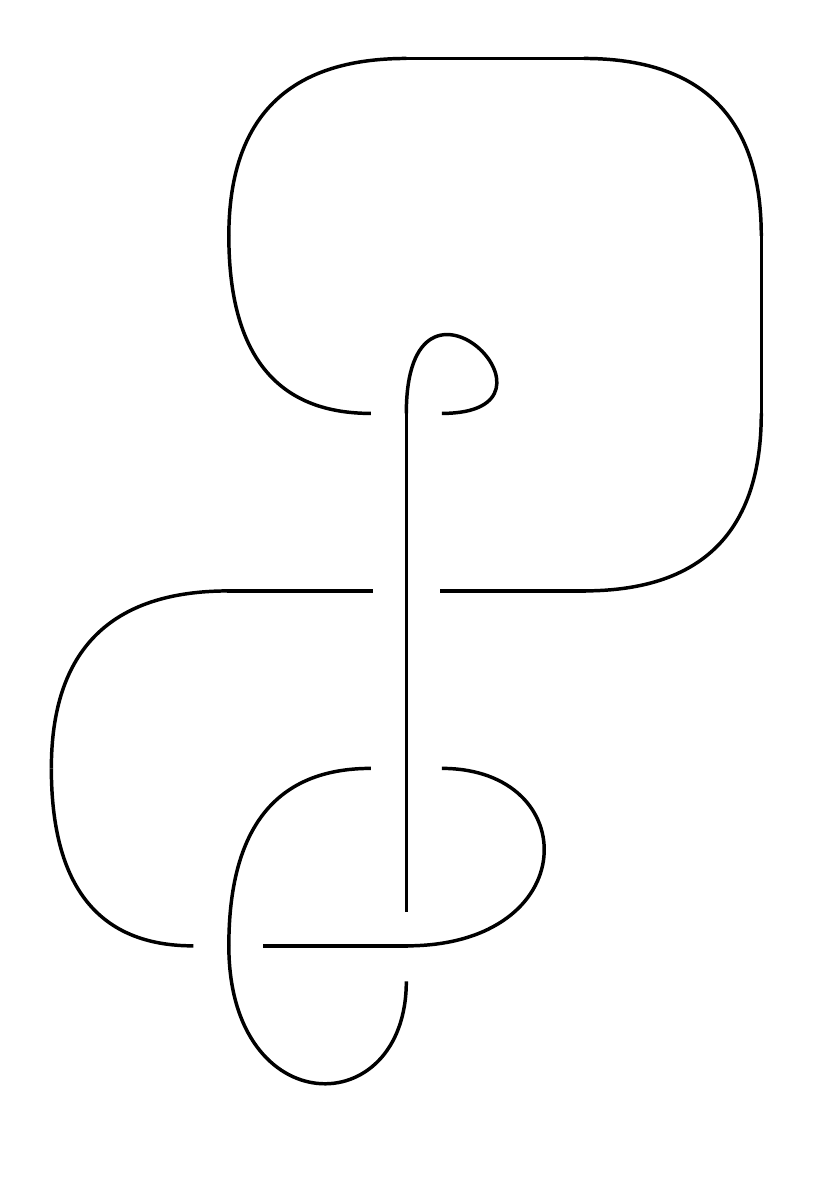}
        \includegraphics[scale=.4]{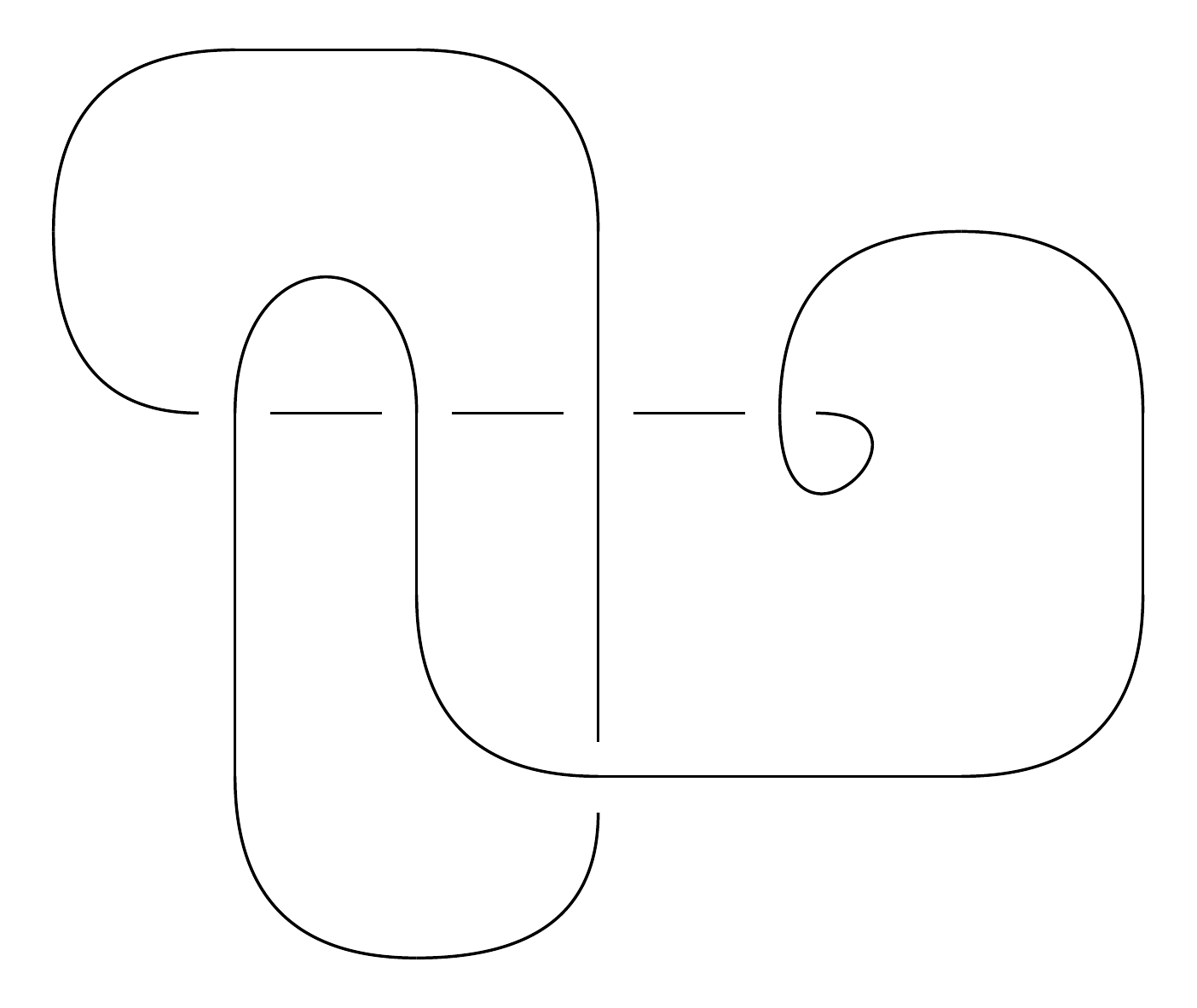}
        \includegraphics[scale=.4]{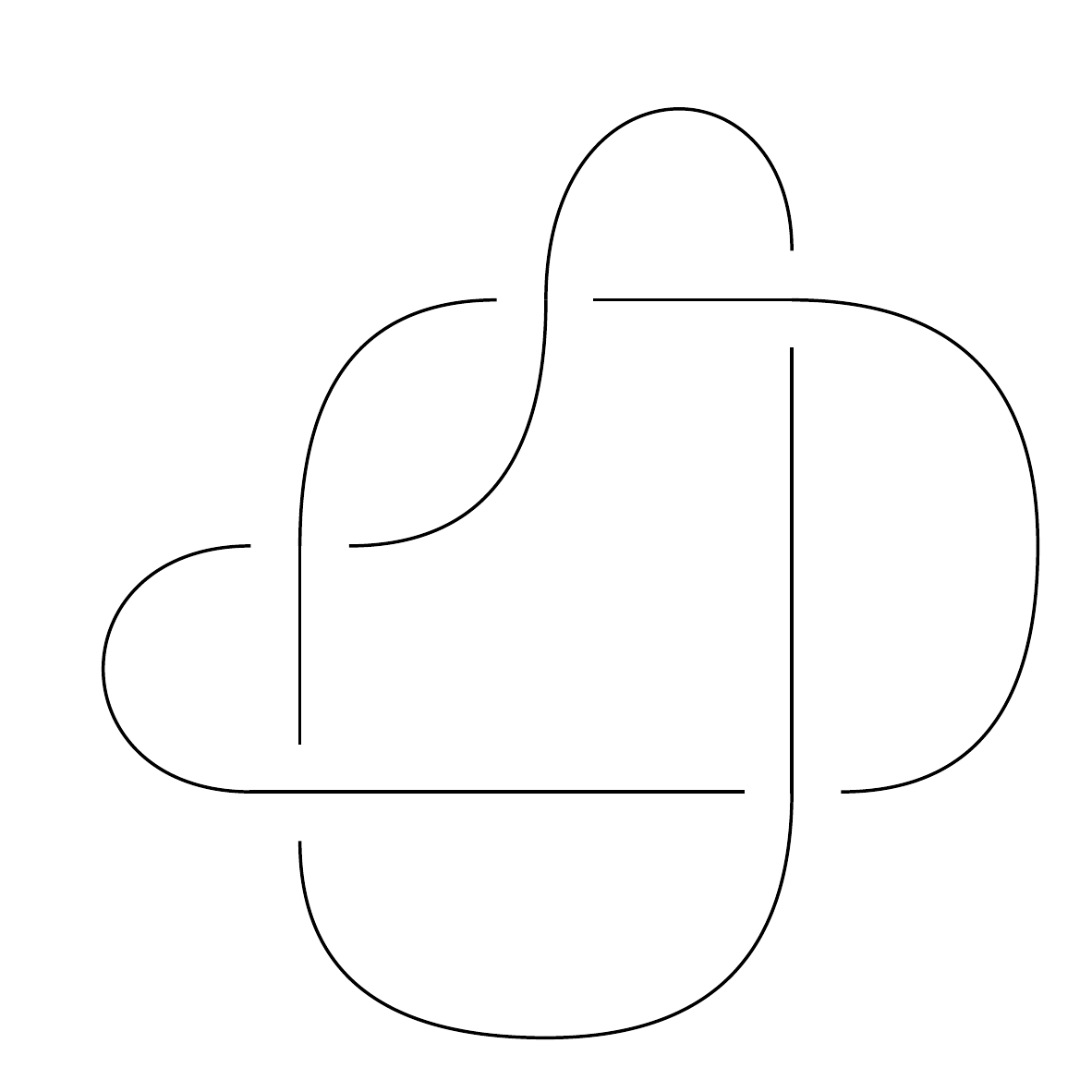}
    } 
    \makebox[\textwidth][c]{
        \includegraphics[scale=.4]{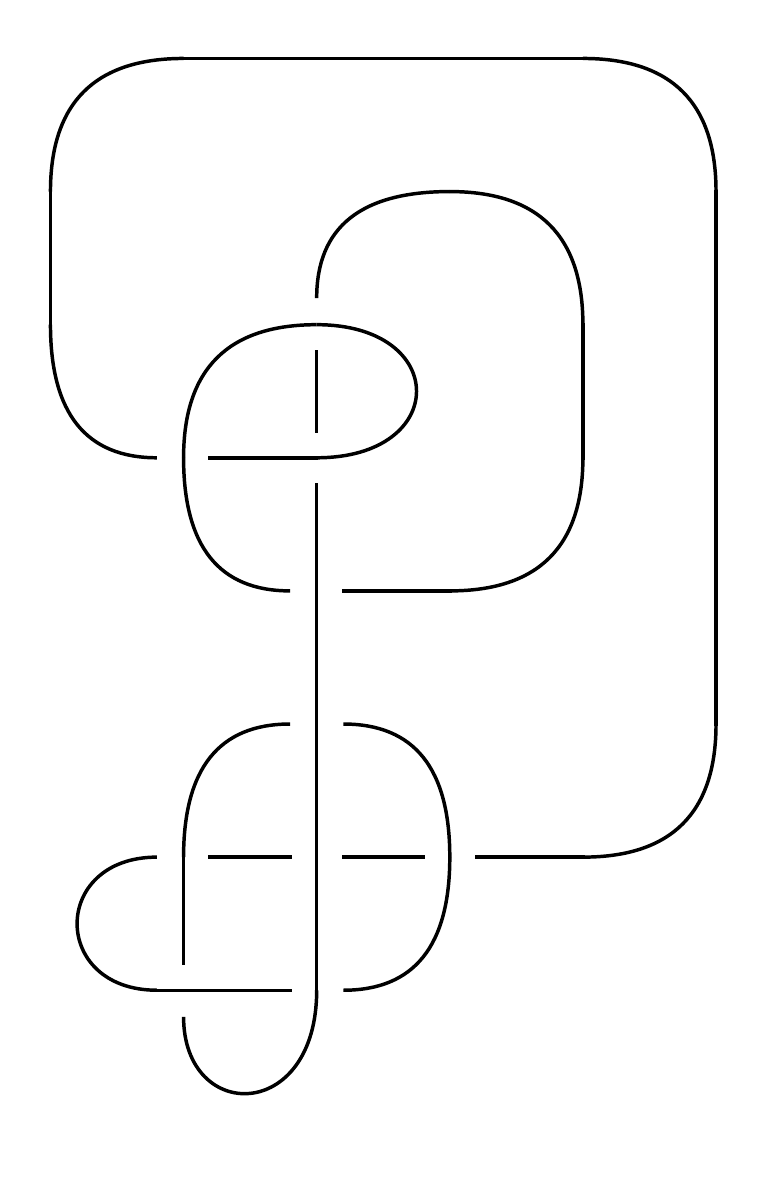}
        \includegraphics[scale=.4]{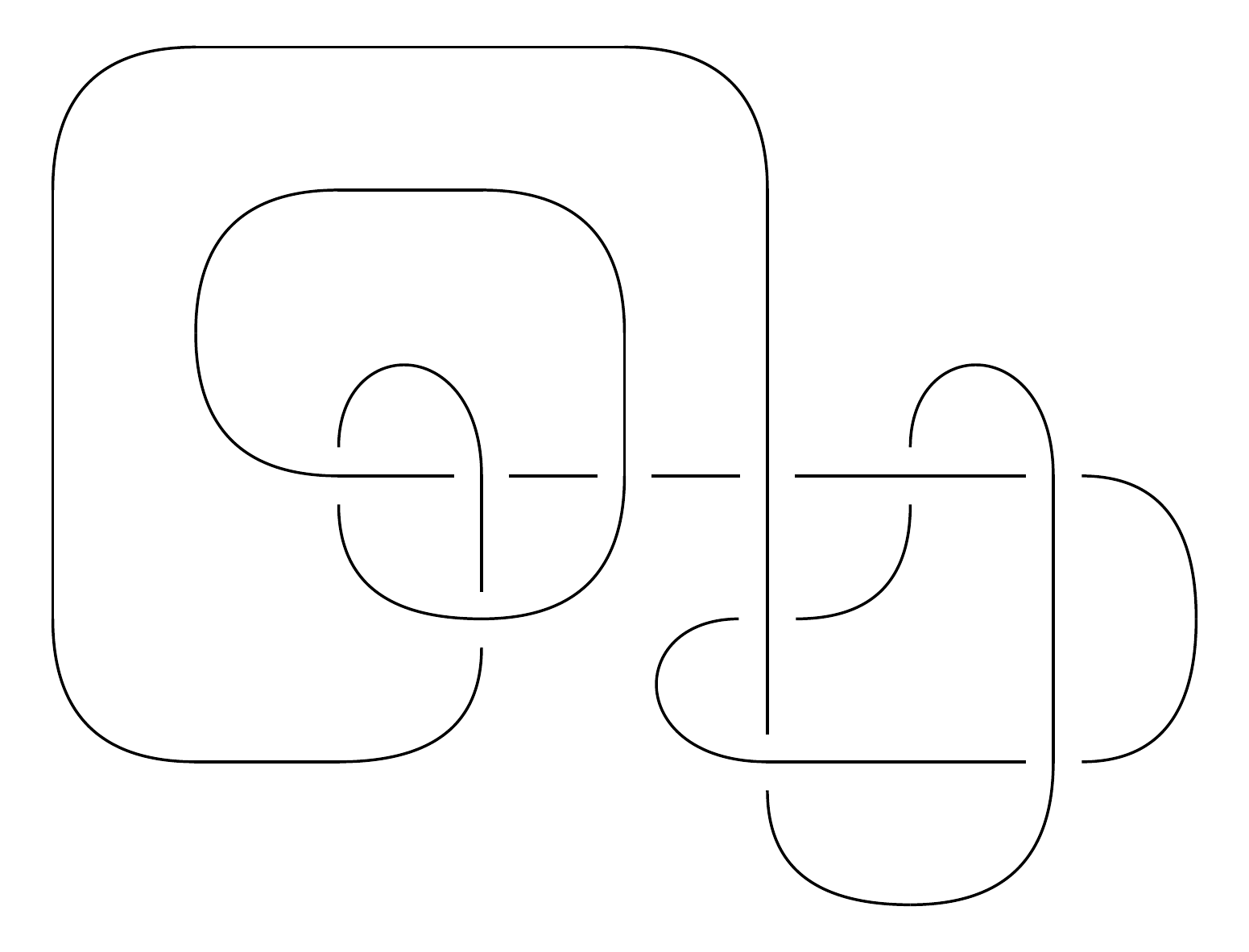}
        \includegraphics[scale=.4]{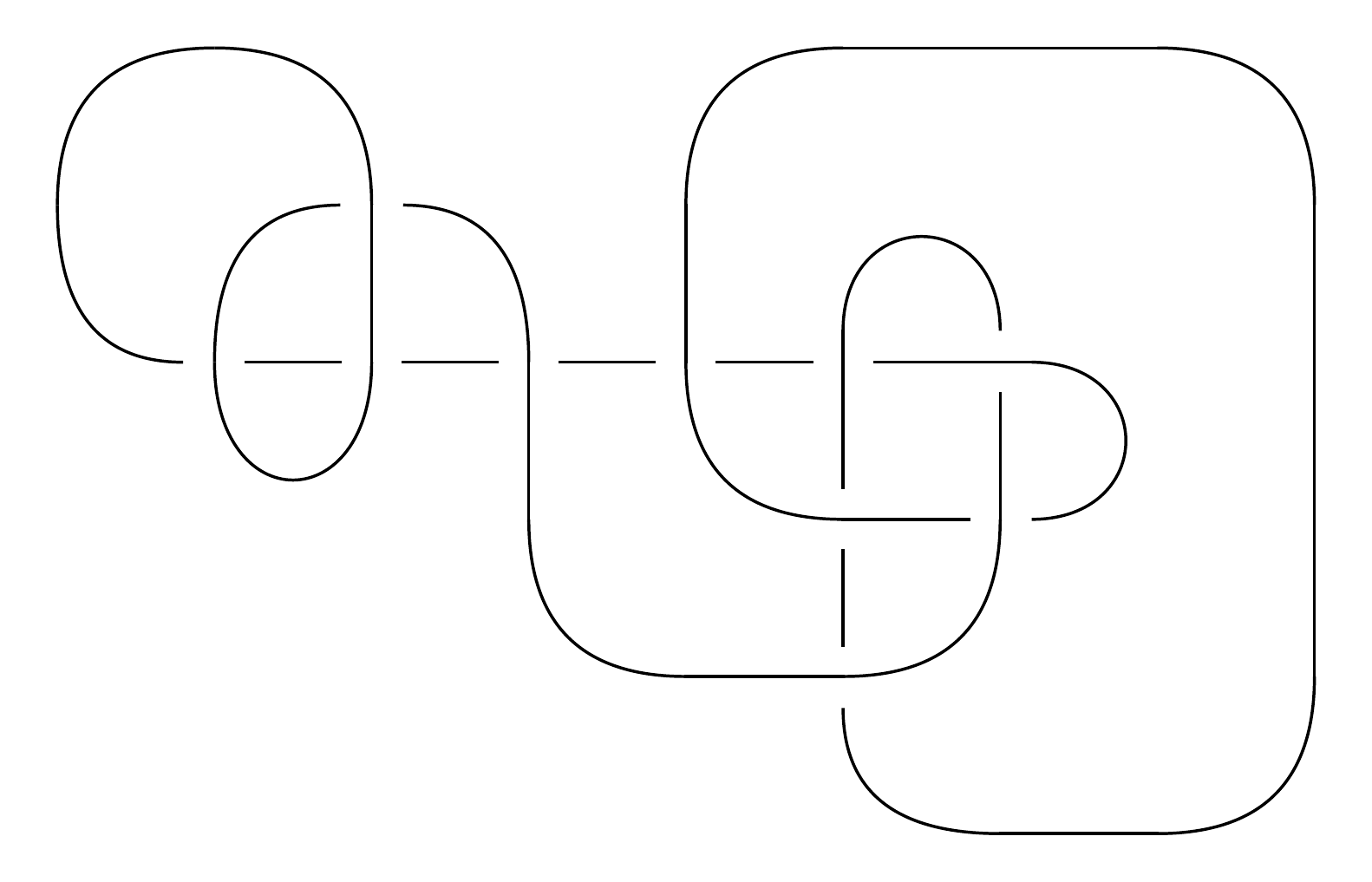}
    }
    \makebox[\textwidth][c]{
        \includegraphics[scale=.4]{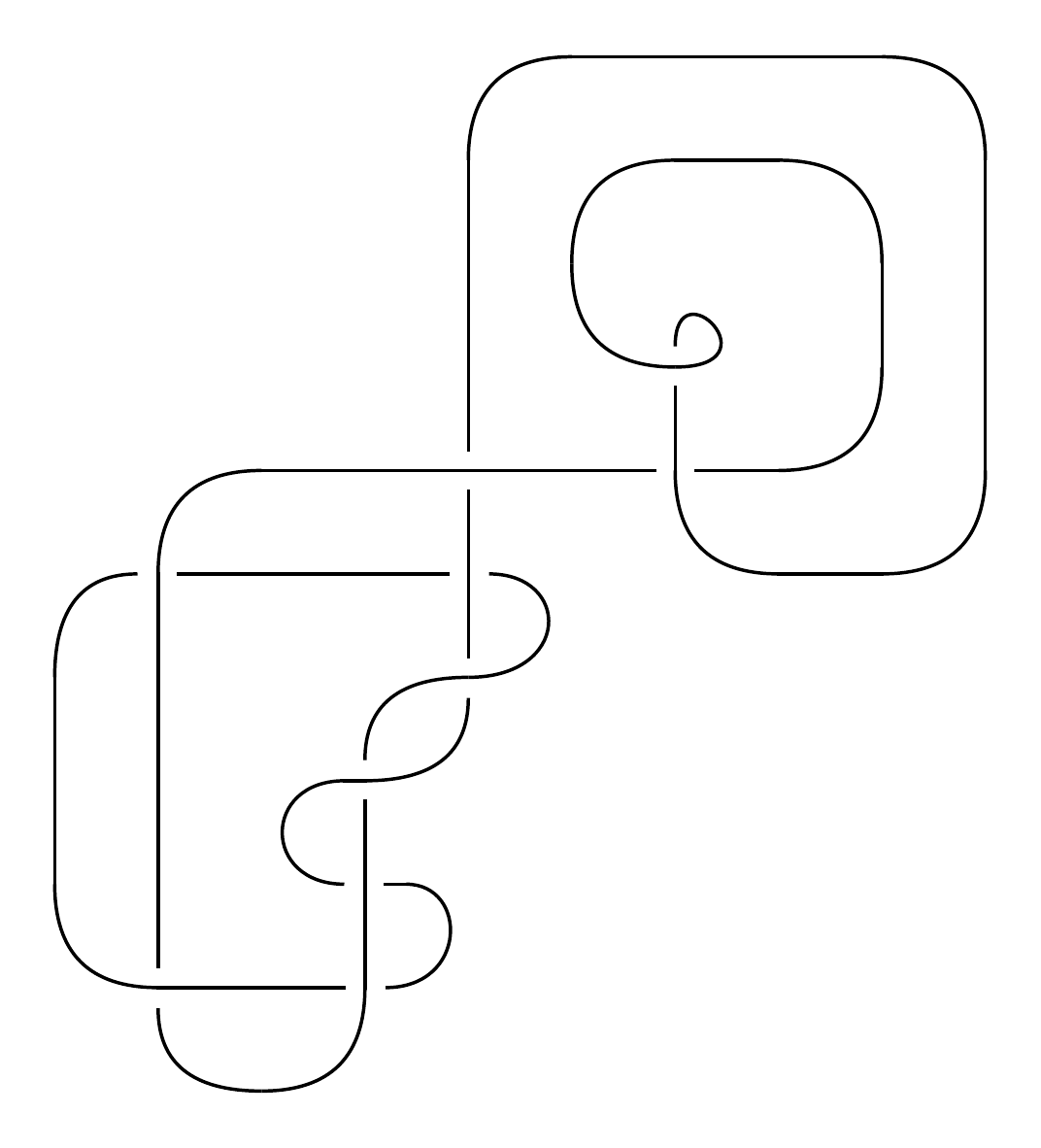}
        \includegraphics[scale=.4]{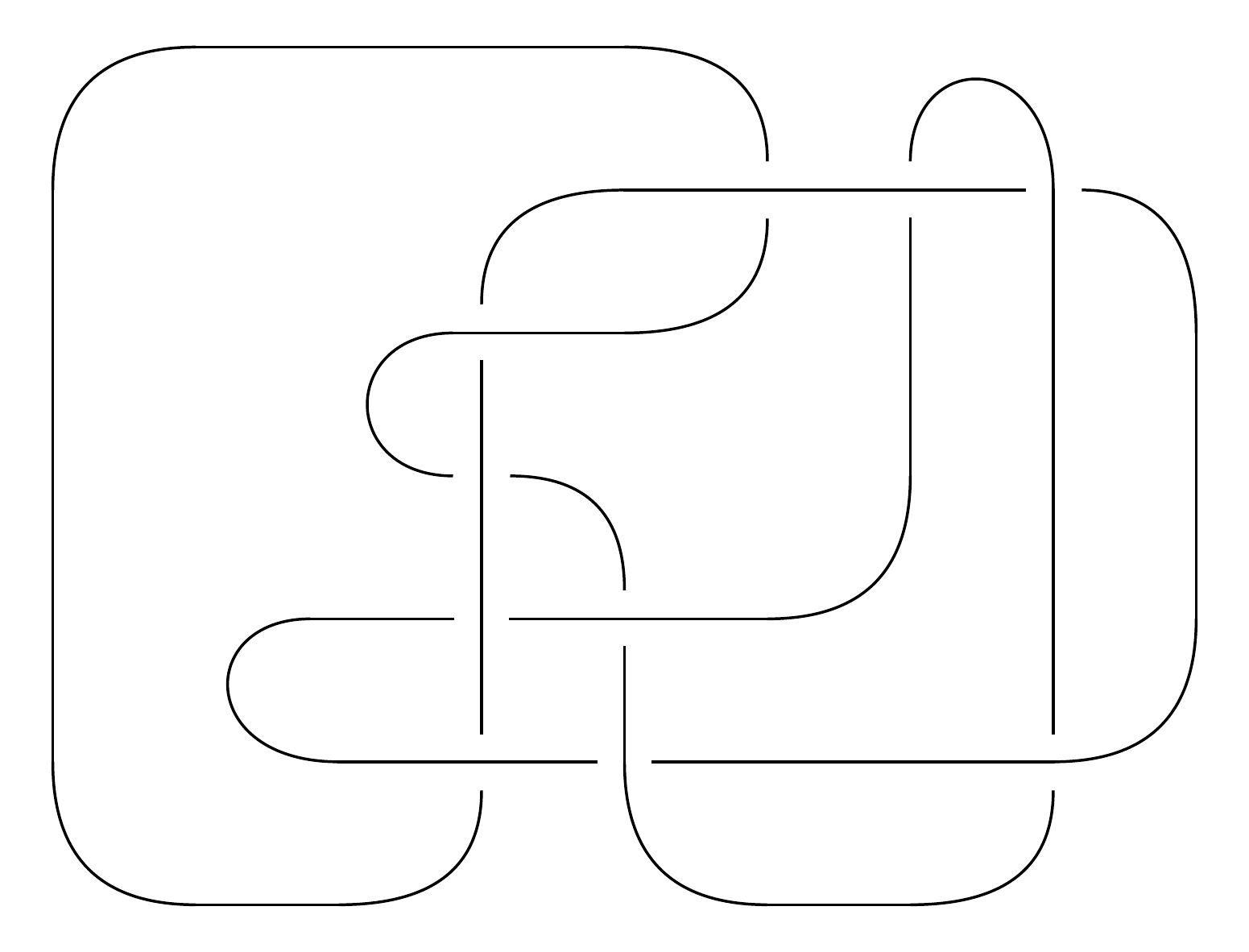}
        \includegraphics[scale=.4]{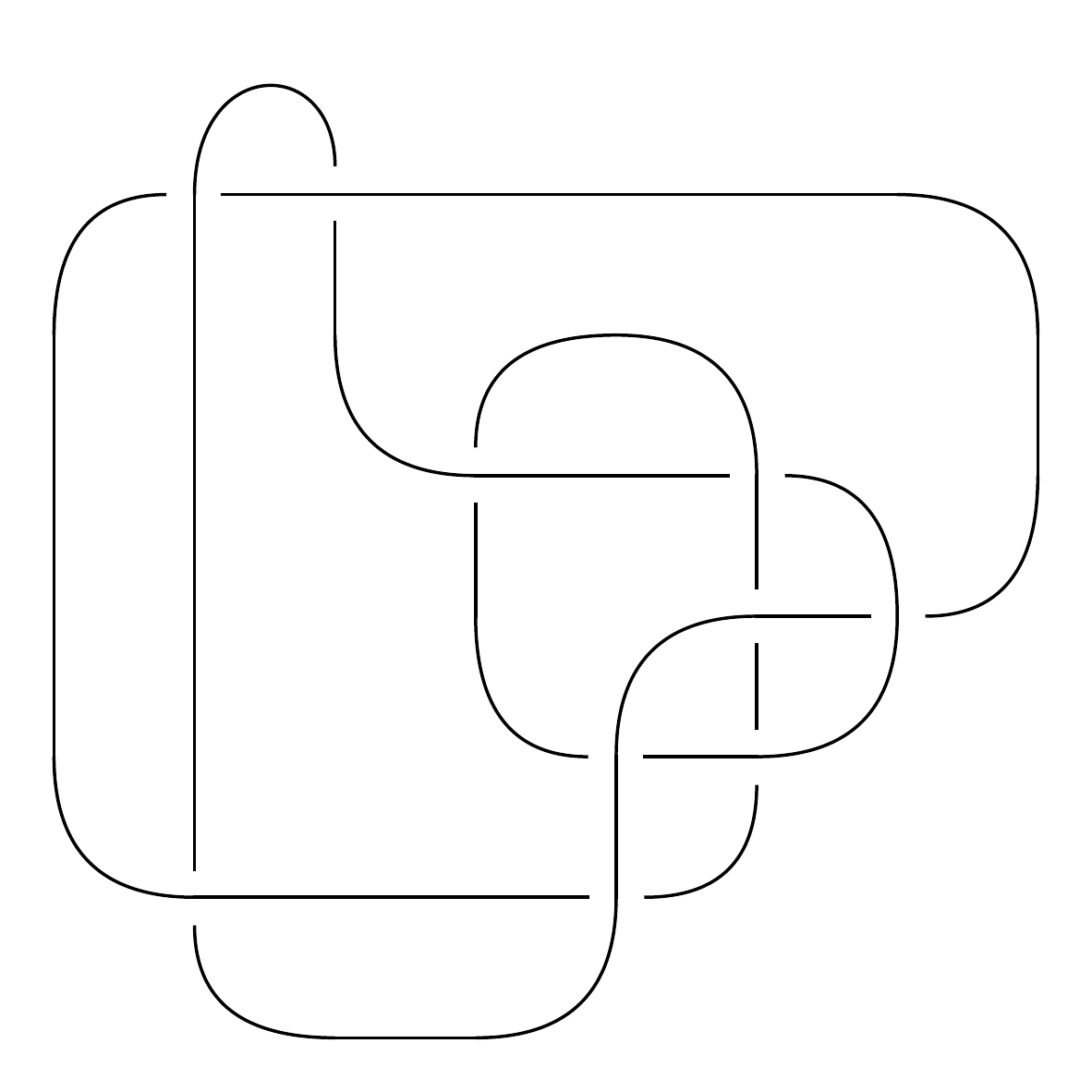}
    }
    \caption{Knot or not? Five and ten crossing in rows $1$-$2$ and $3$-$4$, respectively.}
    \label{fig:knotornot1}
\end{figure}

\begin{figure}
    \makebox[\textwidth][c]{
        \includegraphics[scale=.4]{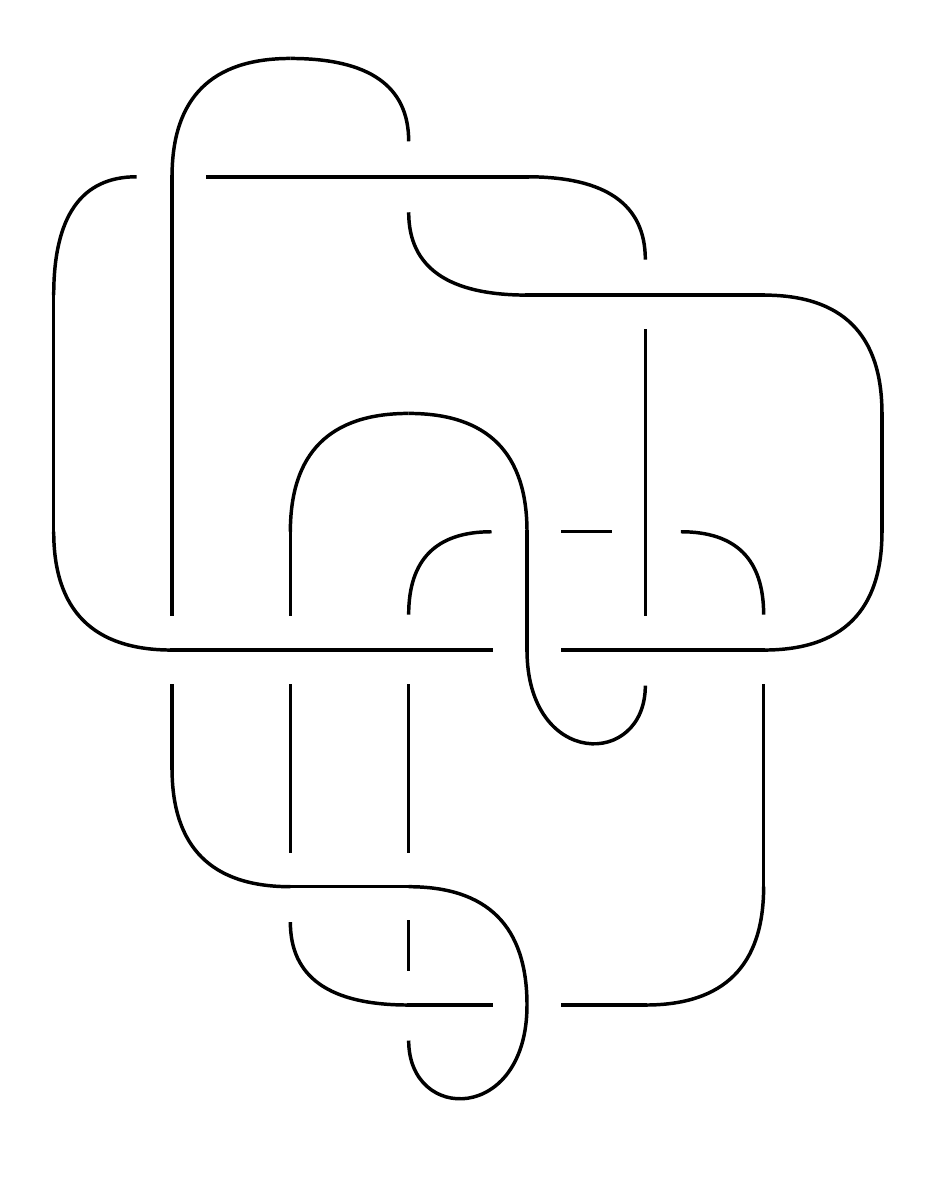}
        \includegraphics[scale=.4]{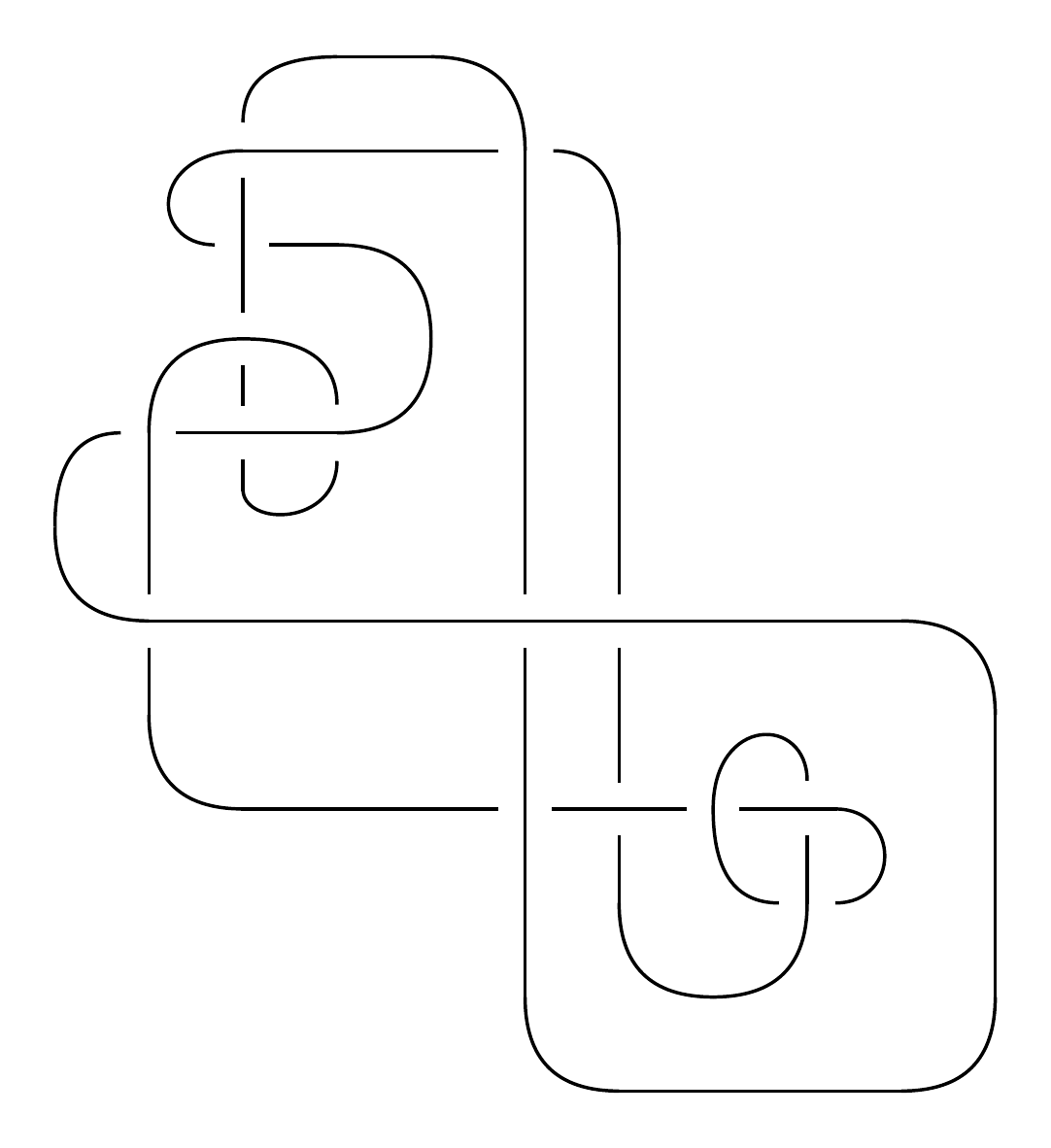}
        \includegraphics[scale=.4]{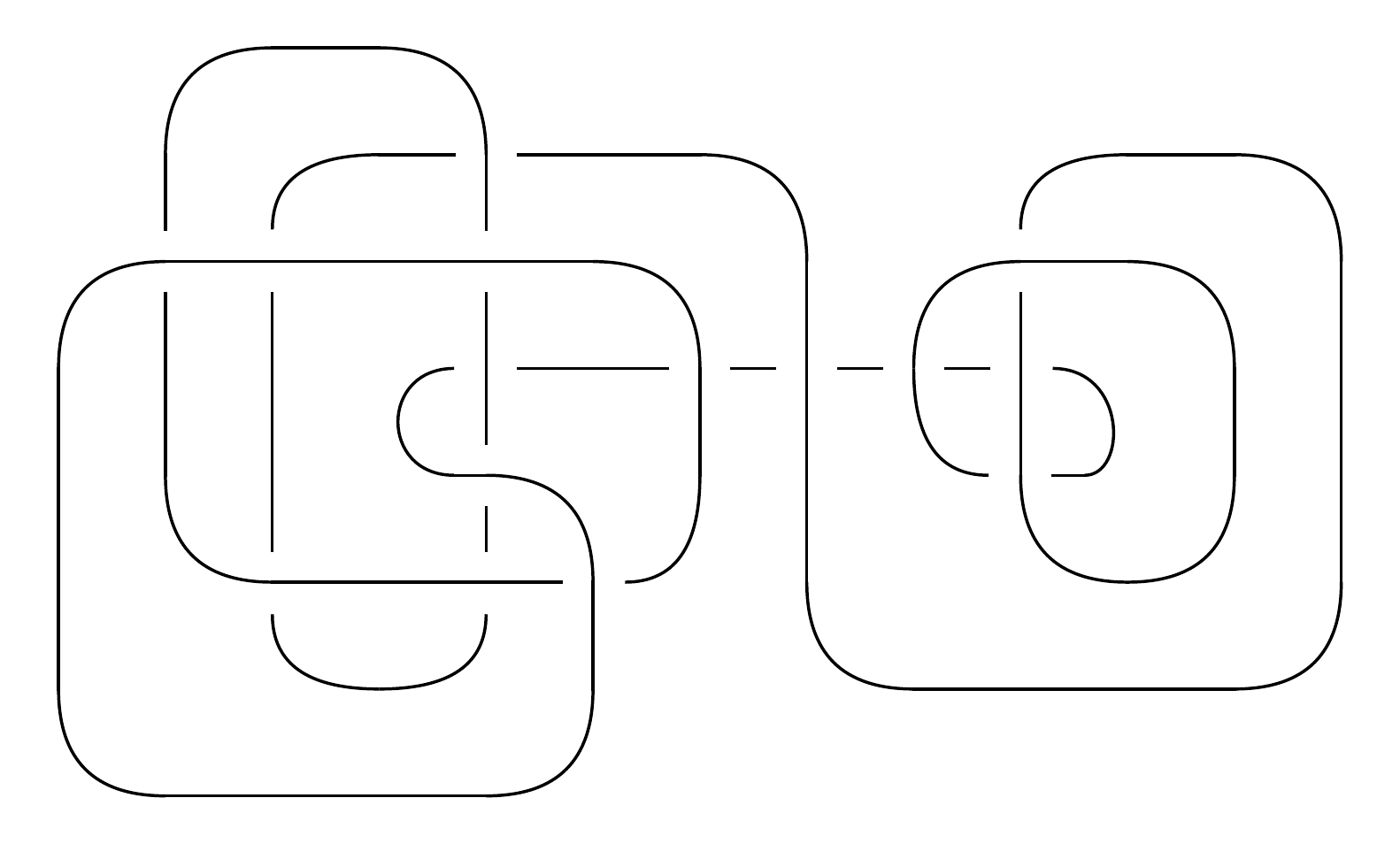}
    } 
    \makebox[\textwidth][c]{
        \includegraphics[scale=.4]{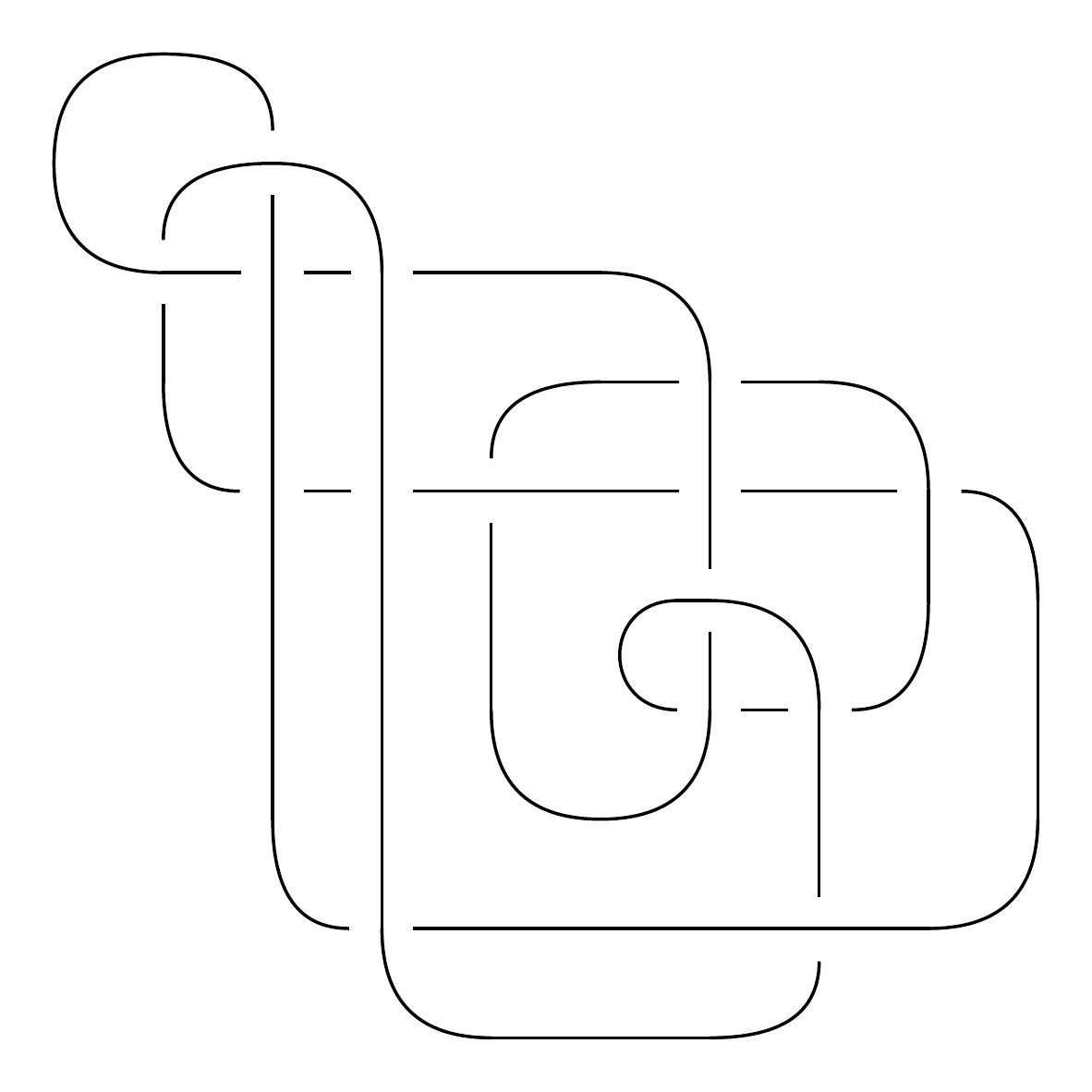}
        \includegraphics[scale=.4]{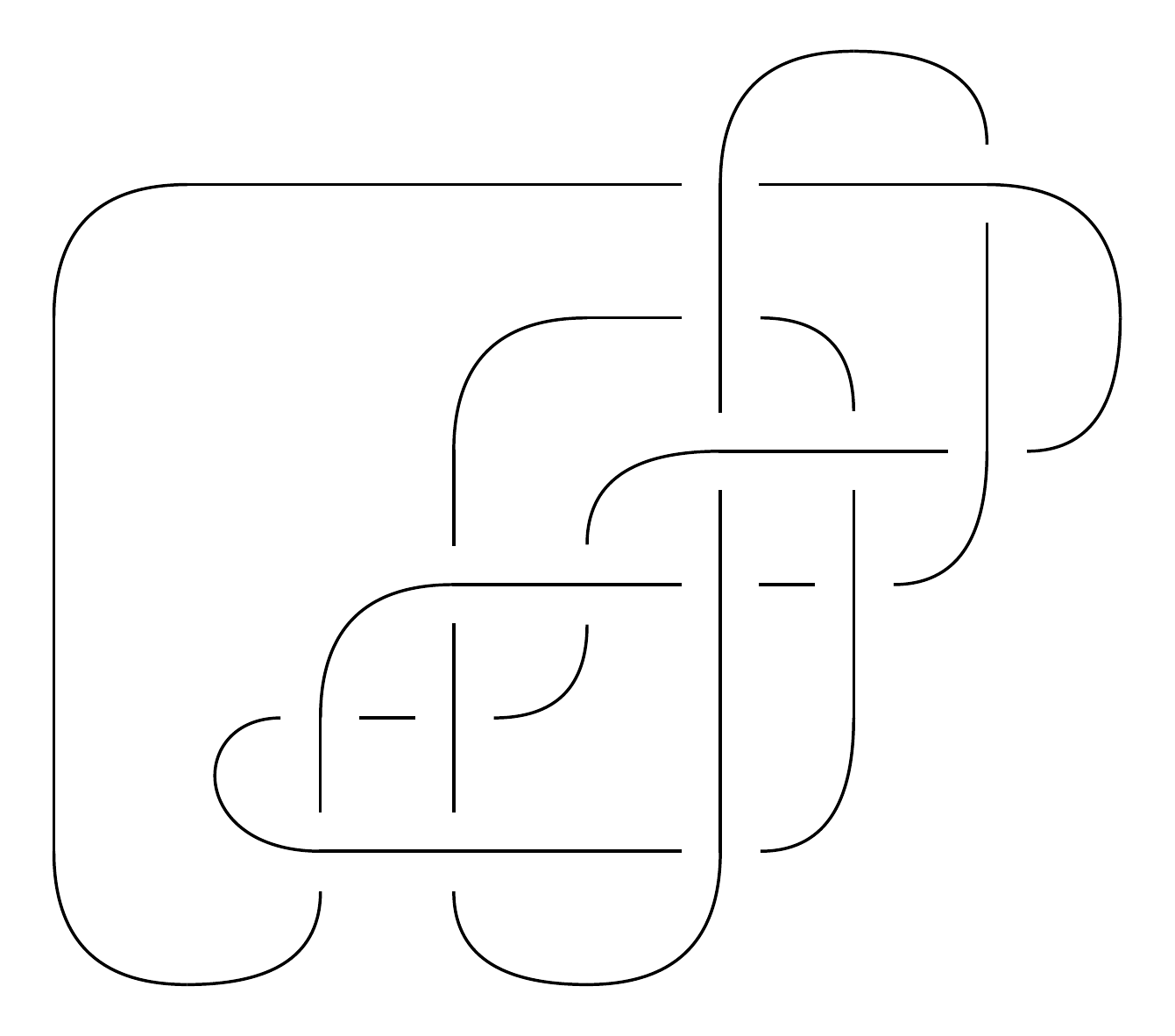}
        \includegraphics[scale=.4]{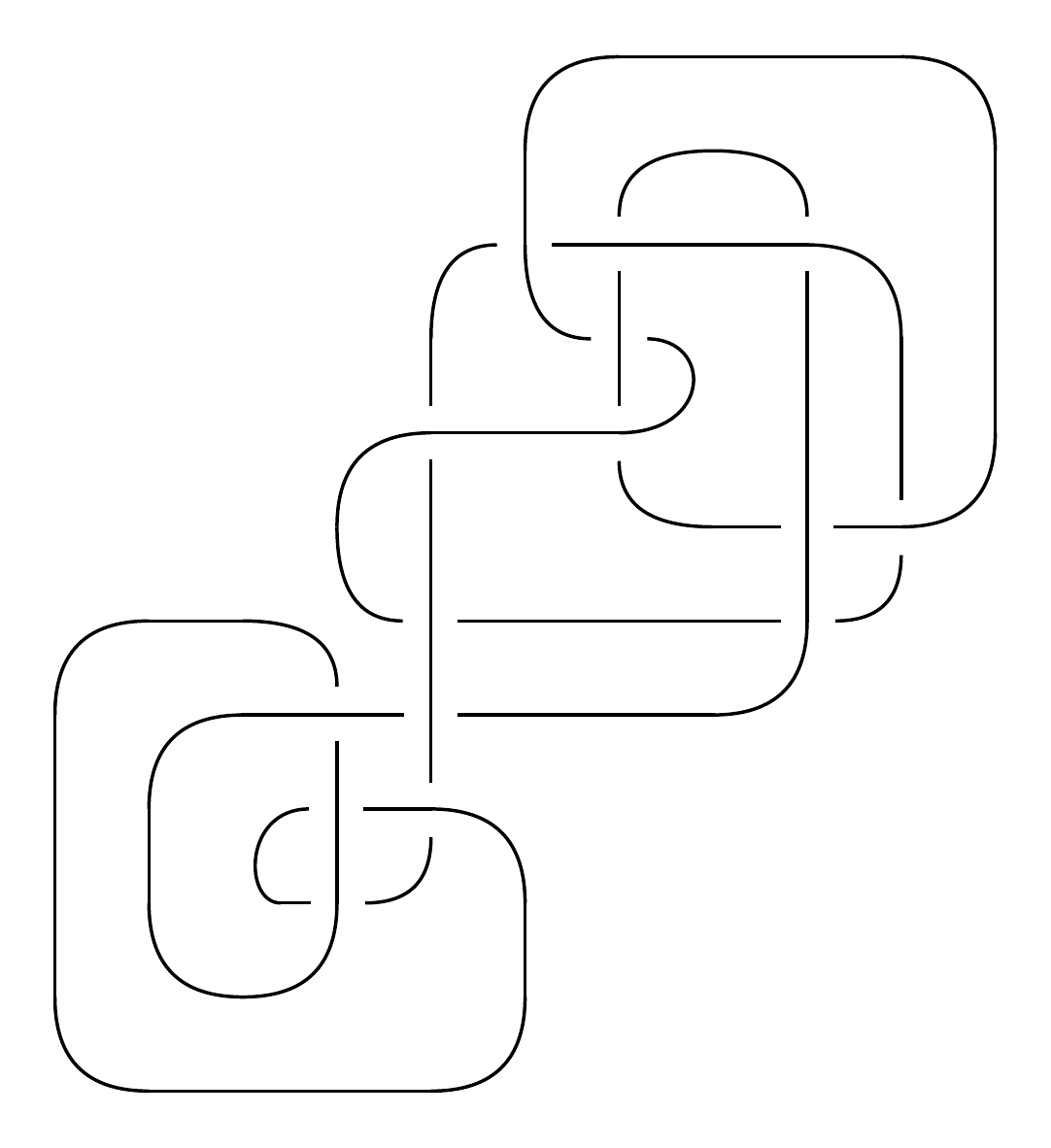}
    } 
    \makebox[\textwidth][c]{
        \includegraphics[scale=.4]{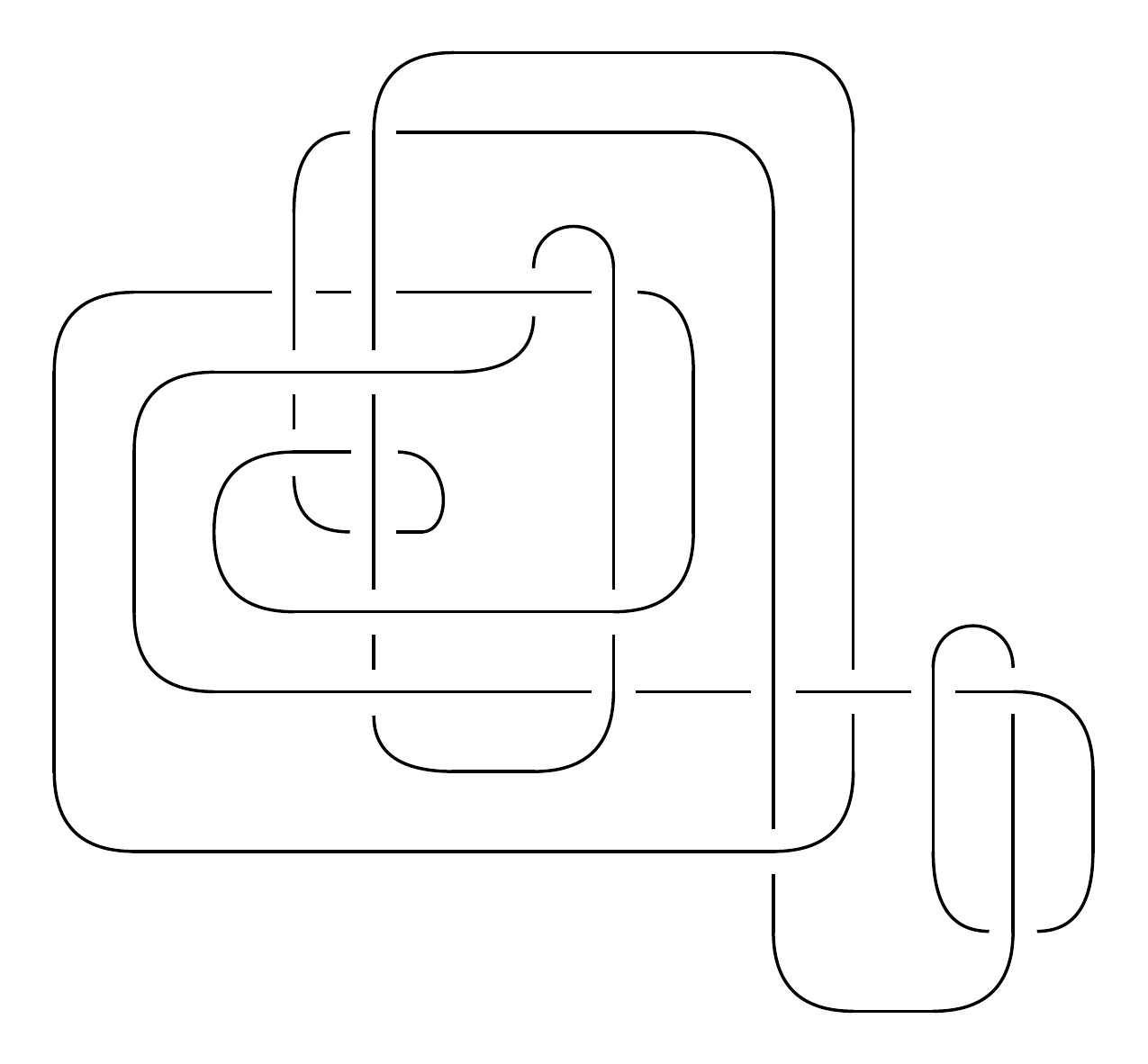}
        \includegraphics[scale=.4]{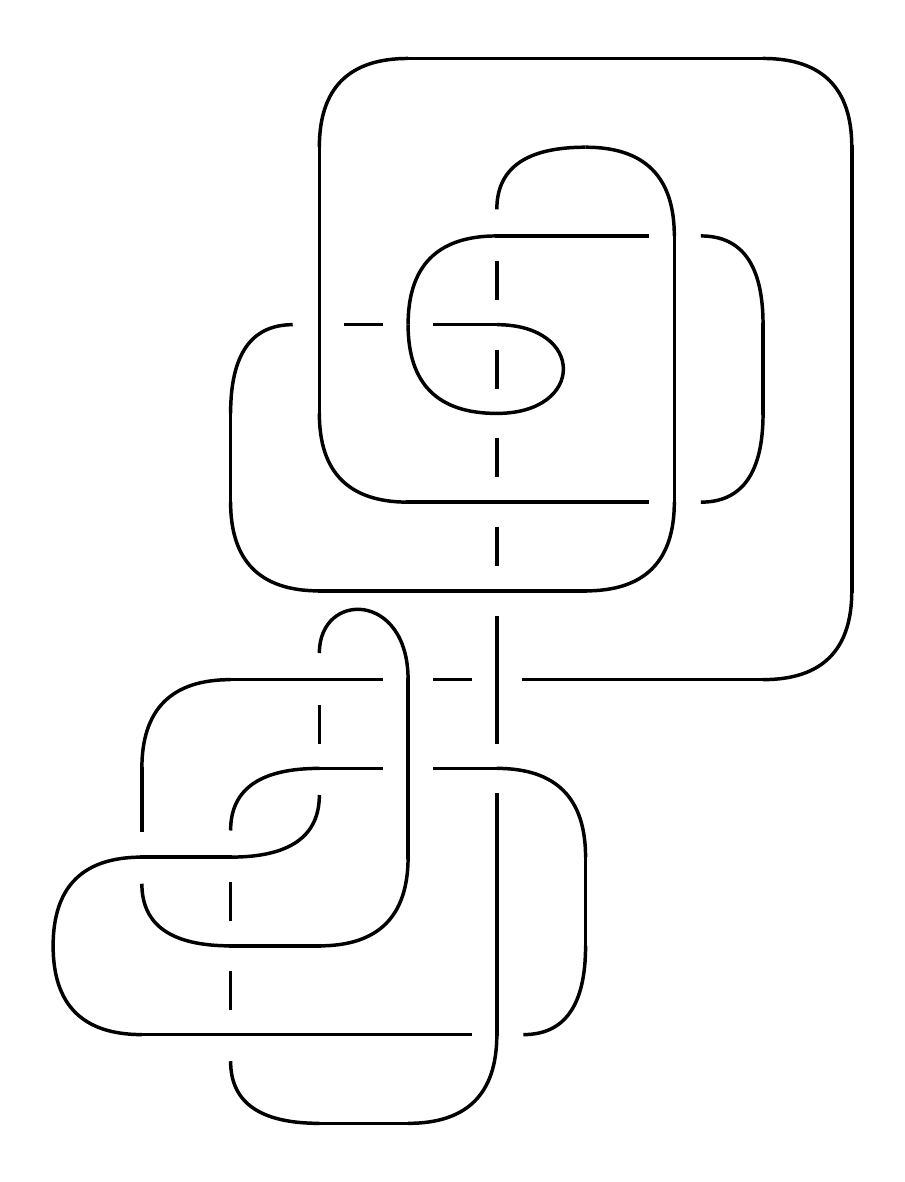}
        \includegraphics[scale=.4]{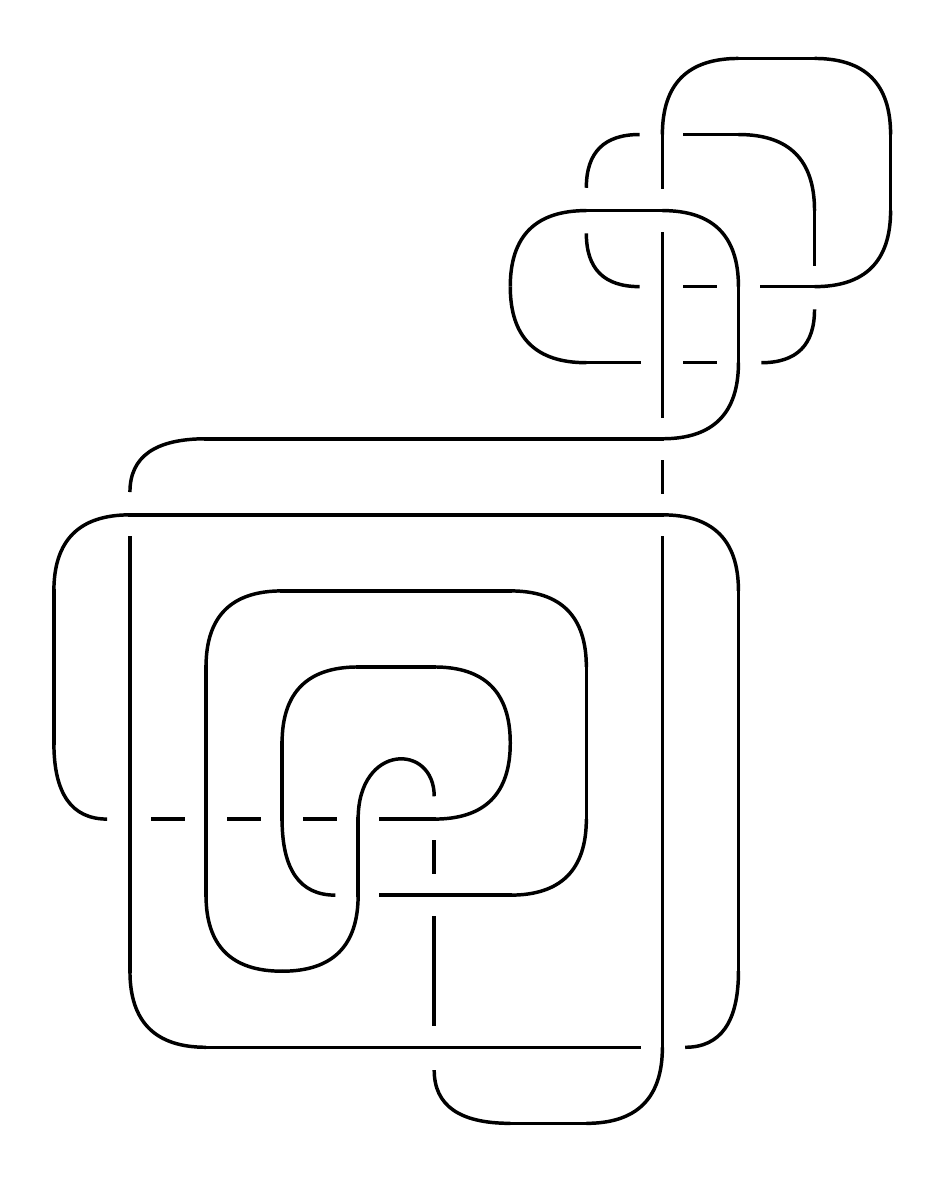}
    }
    \makebox[\textwidth][c]{
        \includegraphics[scale=.4]{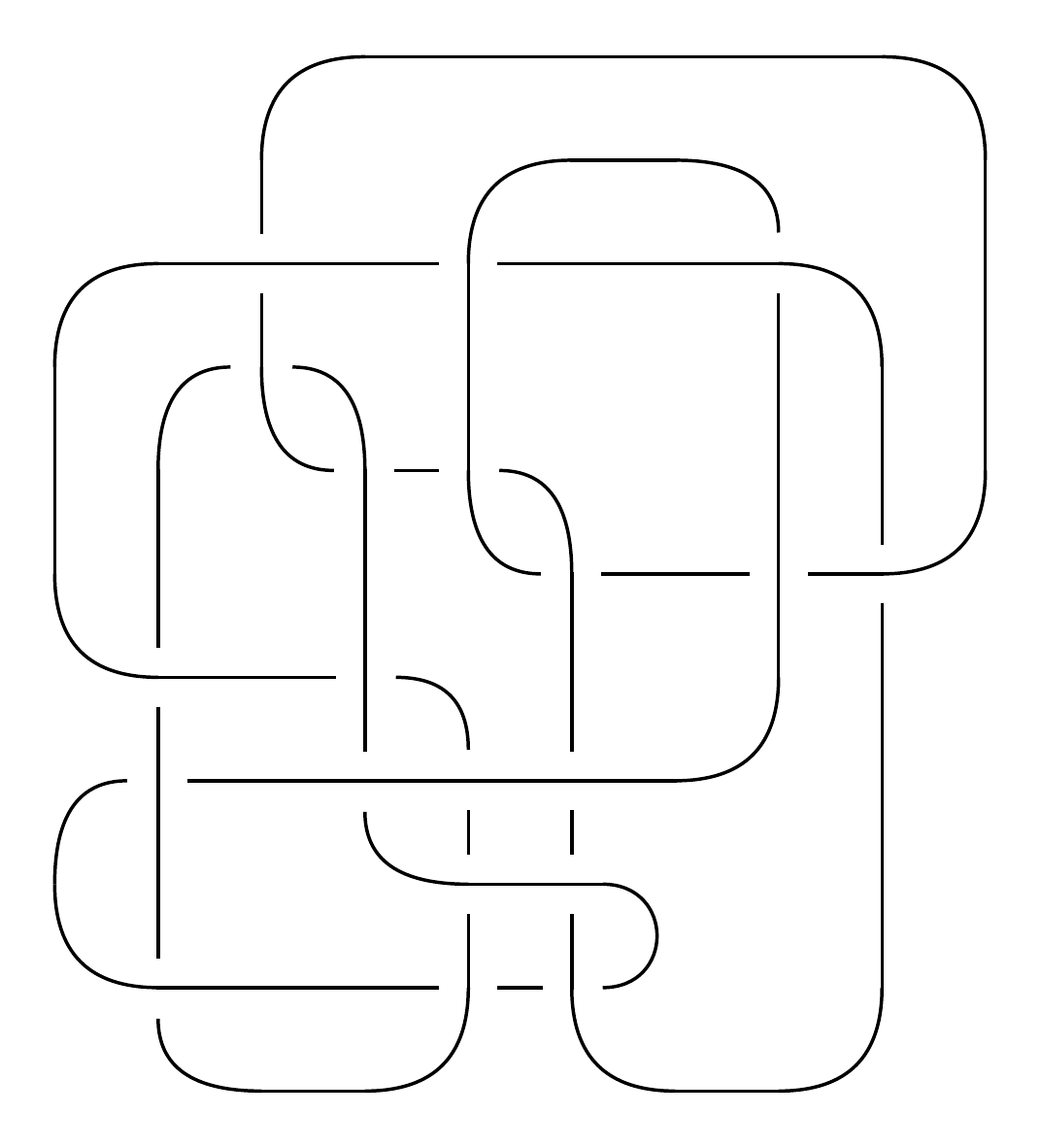}
        \includegraphics[scale=.4]{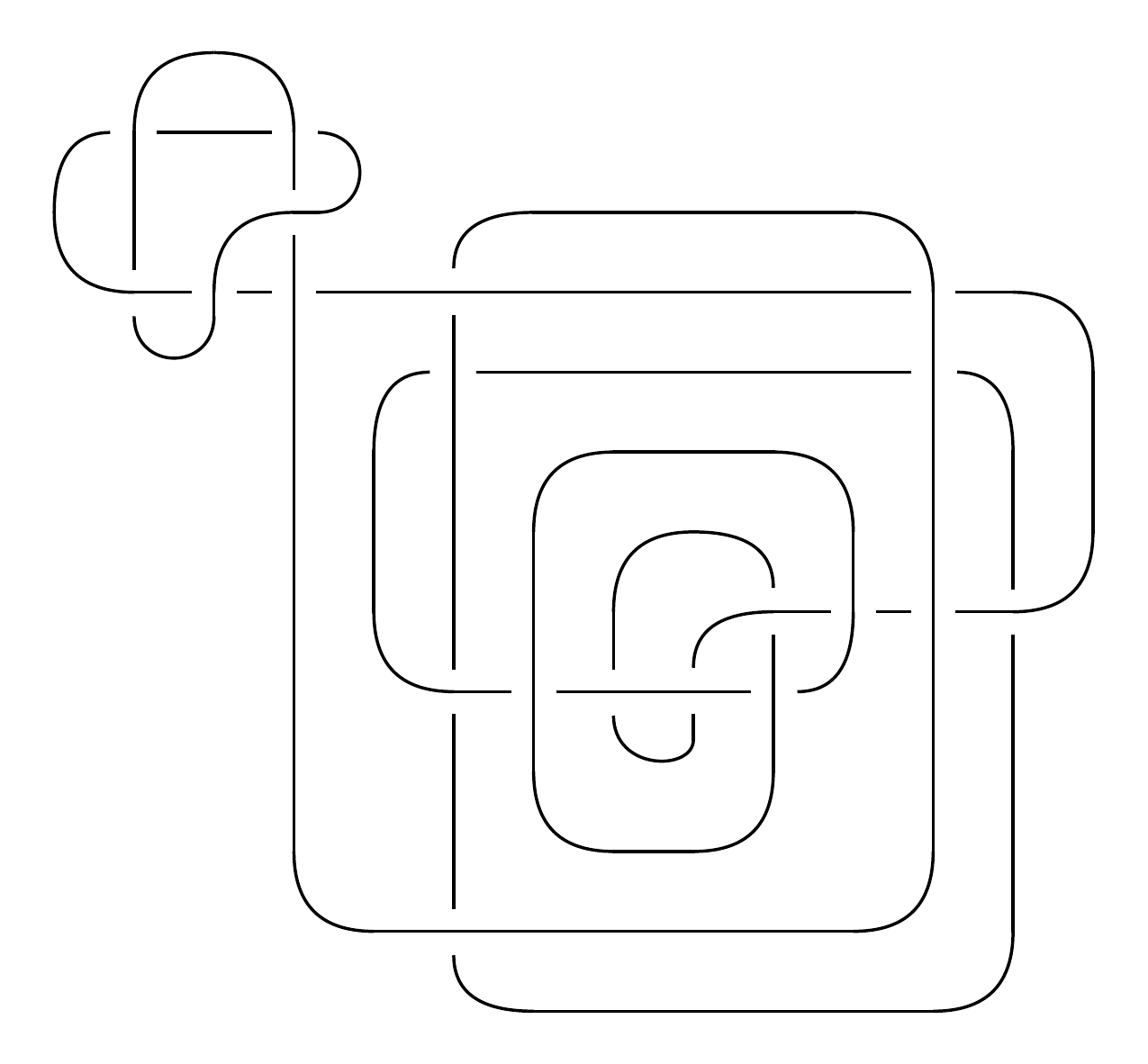}
        \includegraphics[scale=.4]{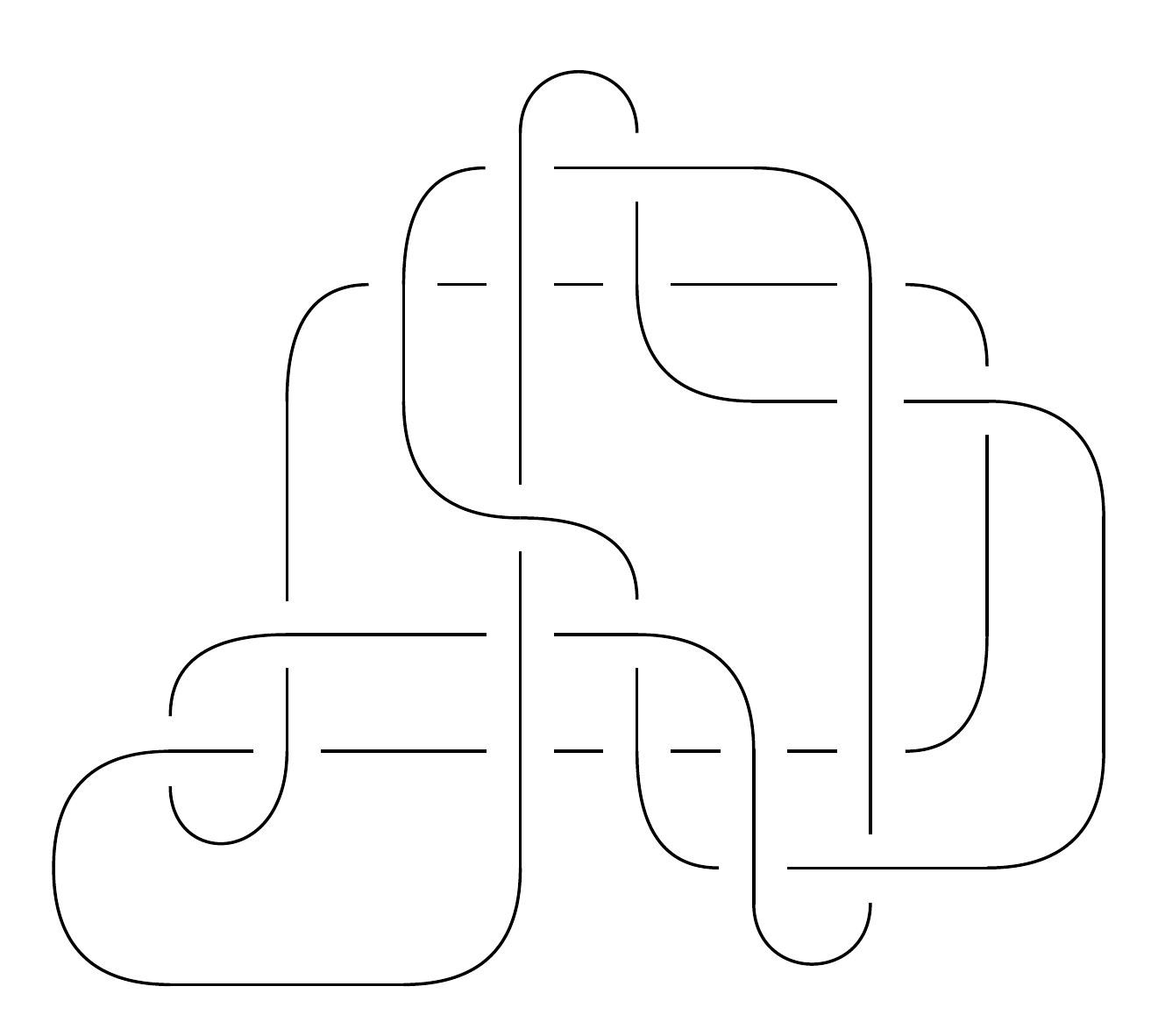}
    }
    \caption{Knot or not? Fifteen and twenty crossing in rows $1$-$2$ and $3$-$4$, respectively.}
    \label{fig:knotornot2}
\end{figure}

\begin{figure}
    \makebox[\textwidth][c]{
        \includegraphics[scale=.4]{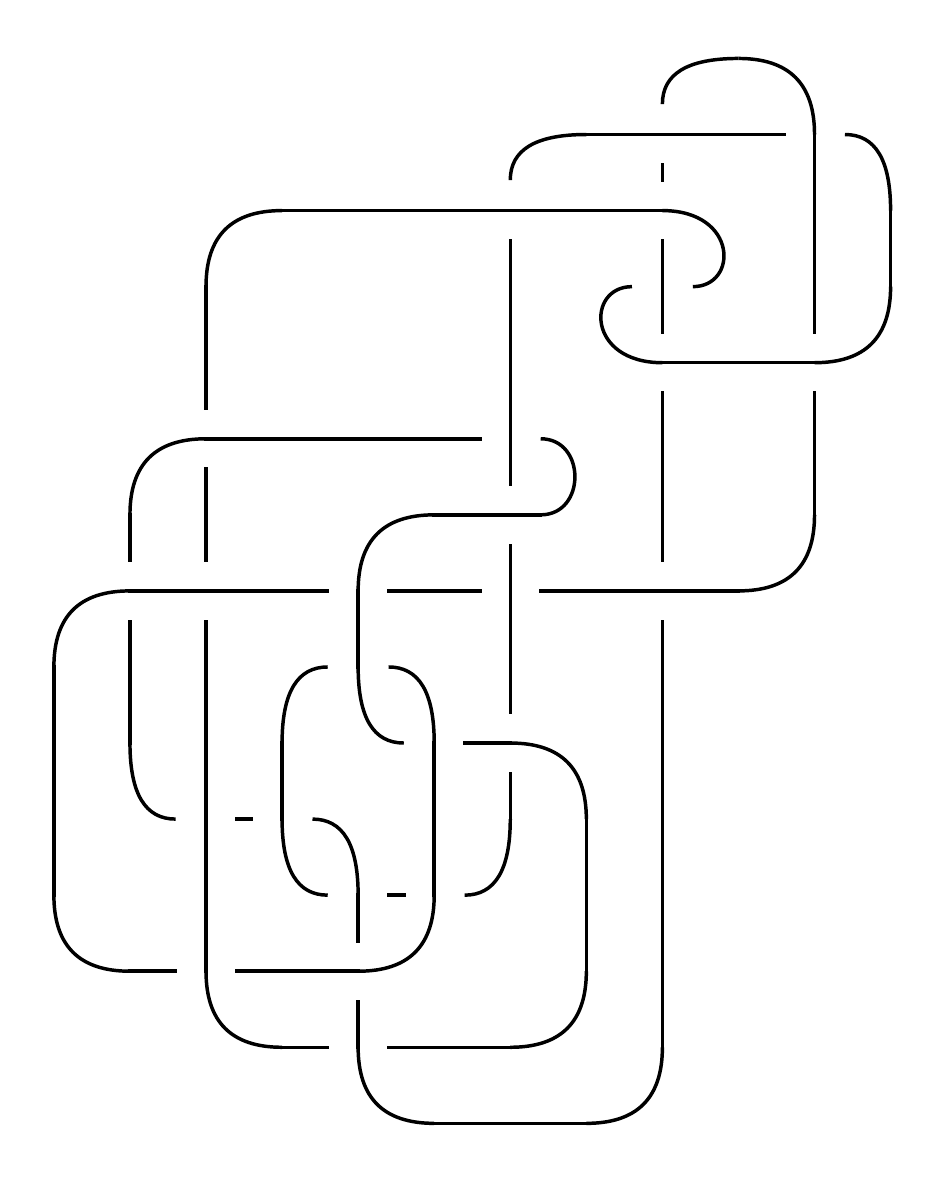}
        \includegraphics[scale=.4]{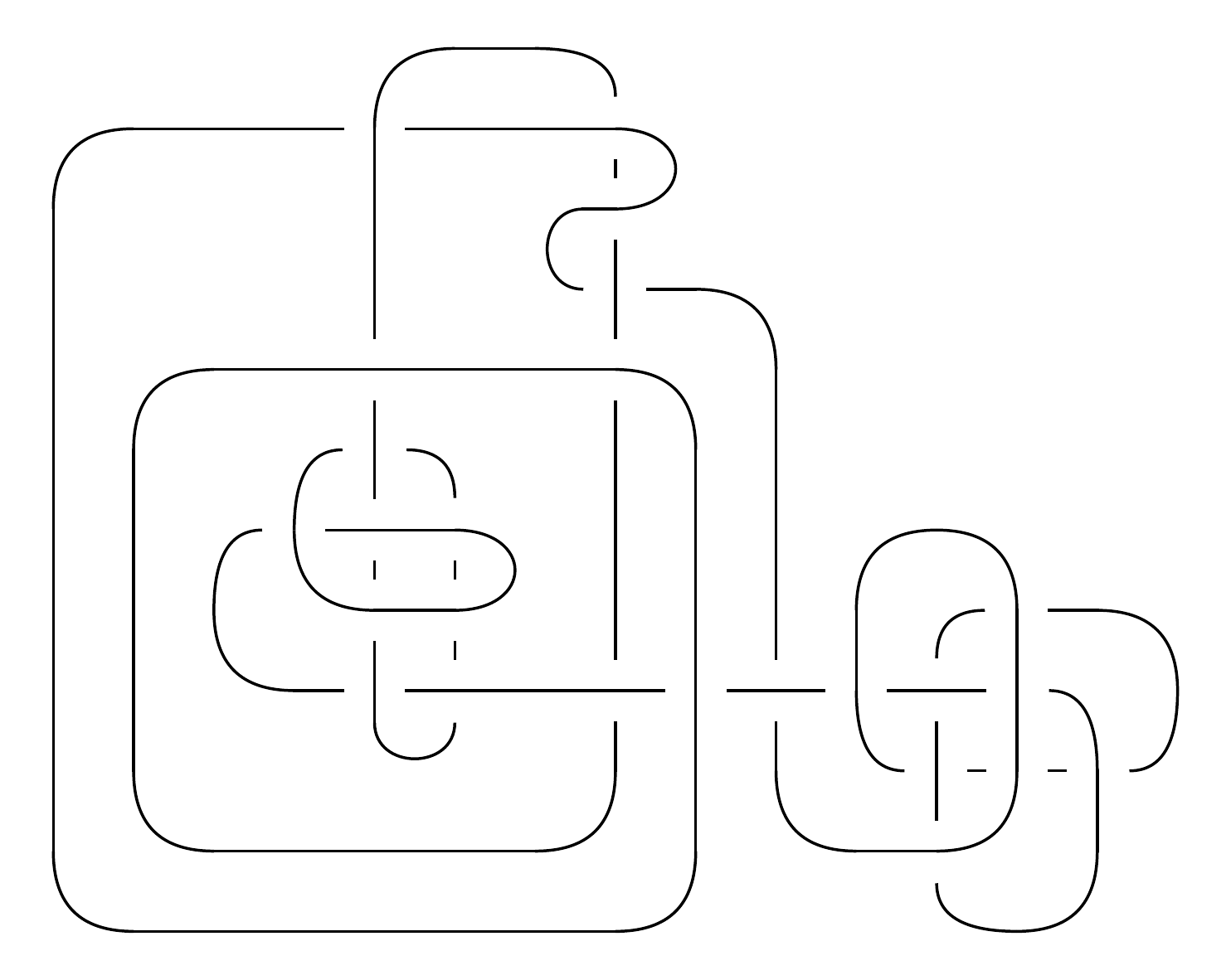}
        \includegraphics[scale=.4]{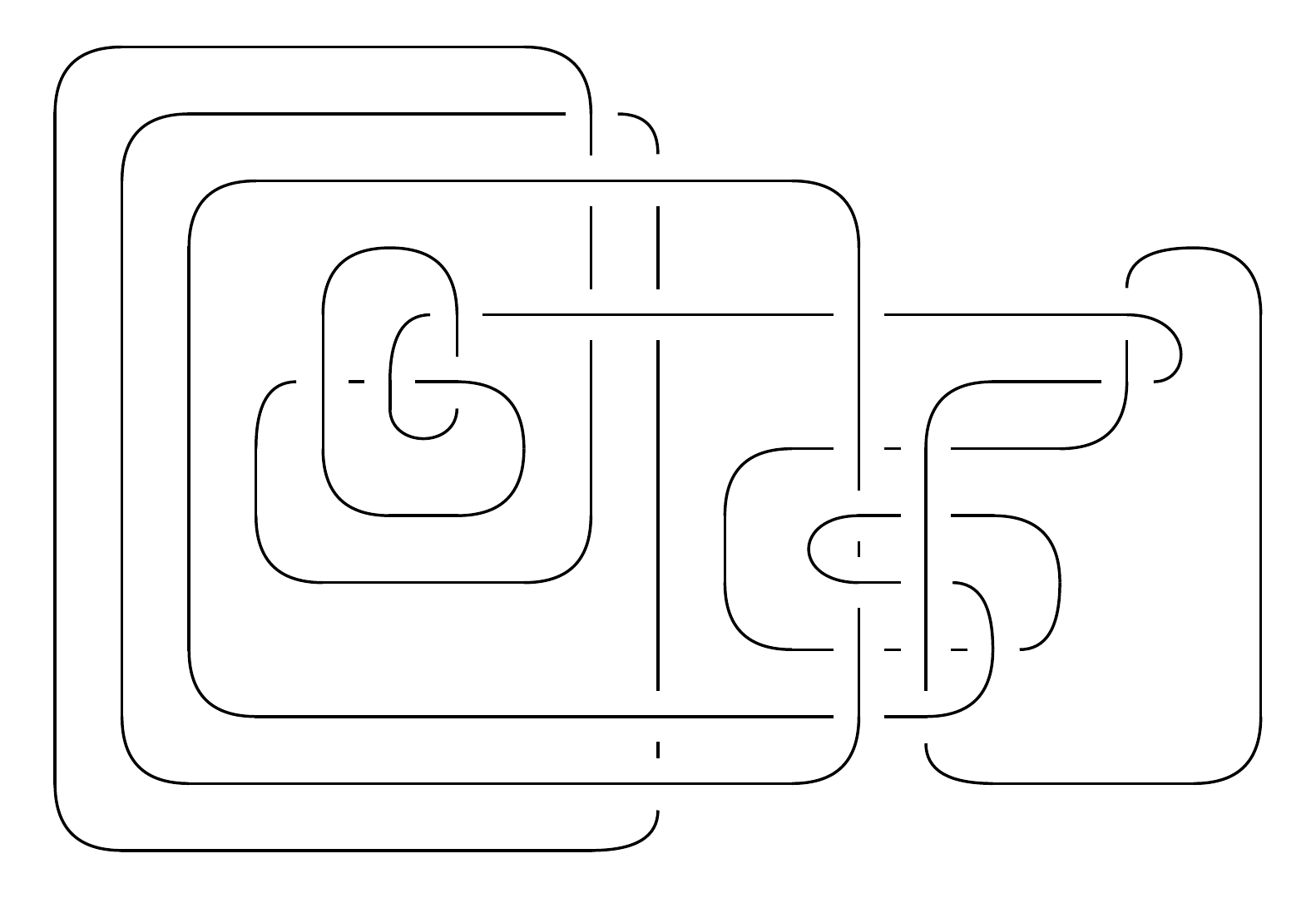}
    } 
    \makebox[\textwidth][c]{
        \includegraphics[scale=.4]{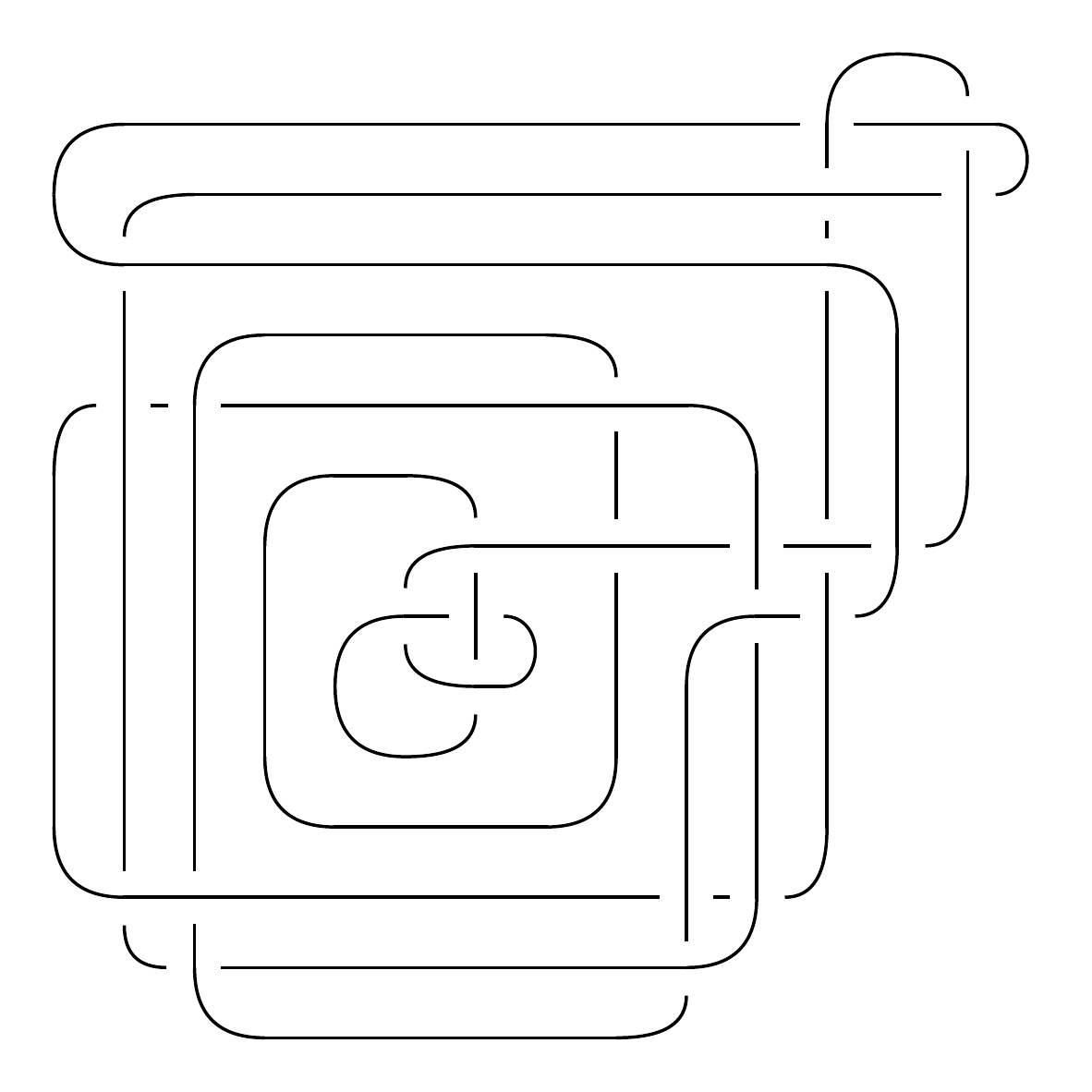}
        \includegraphics[scale=.4]{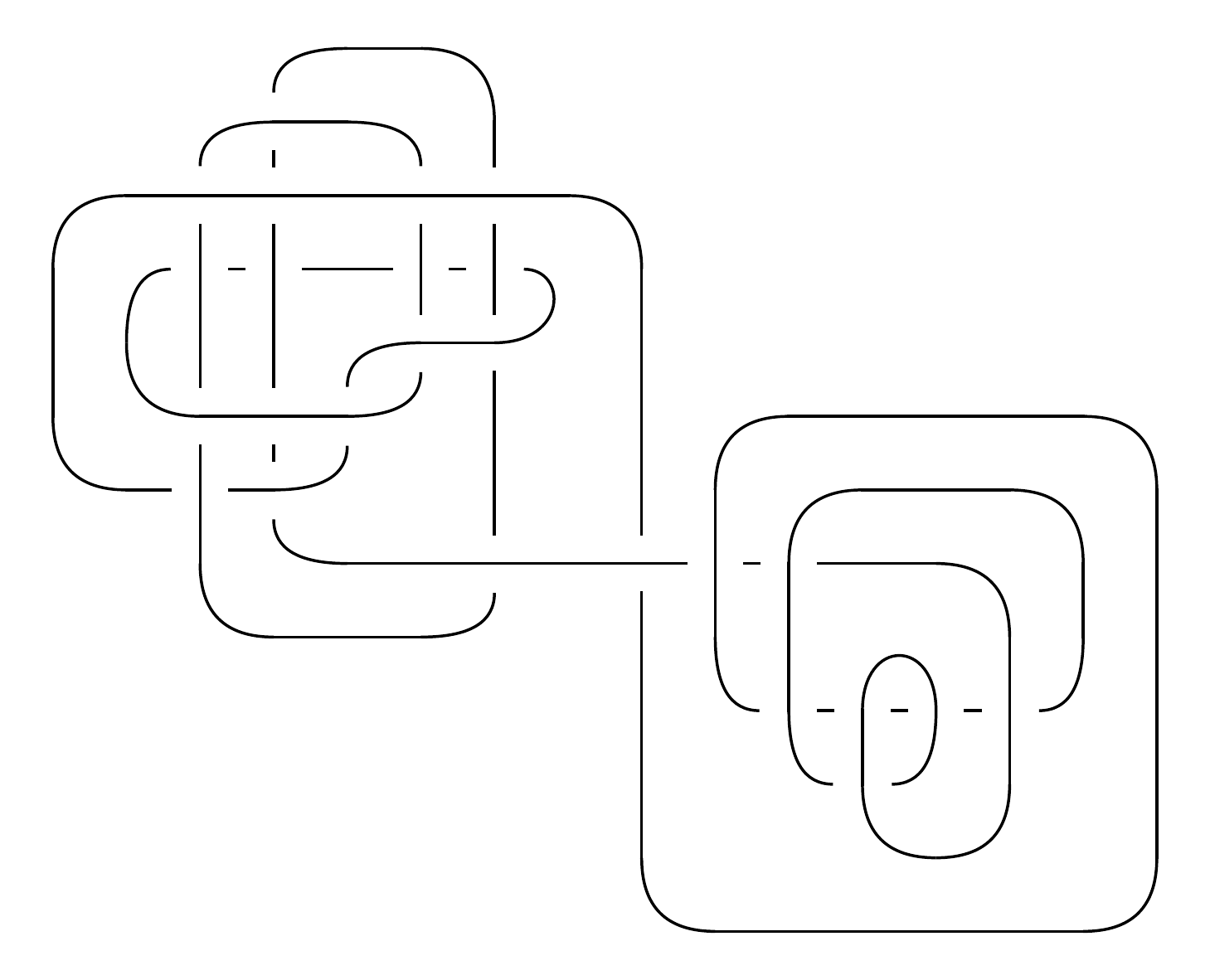}
        \includegraphics[scale=.4]{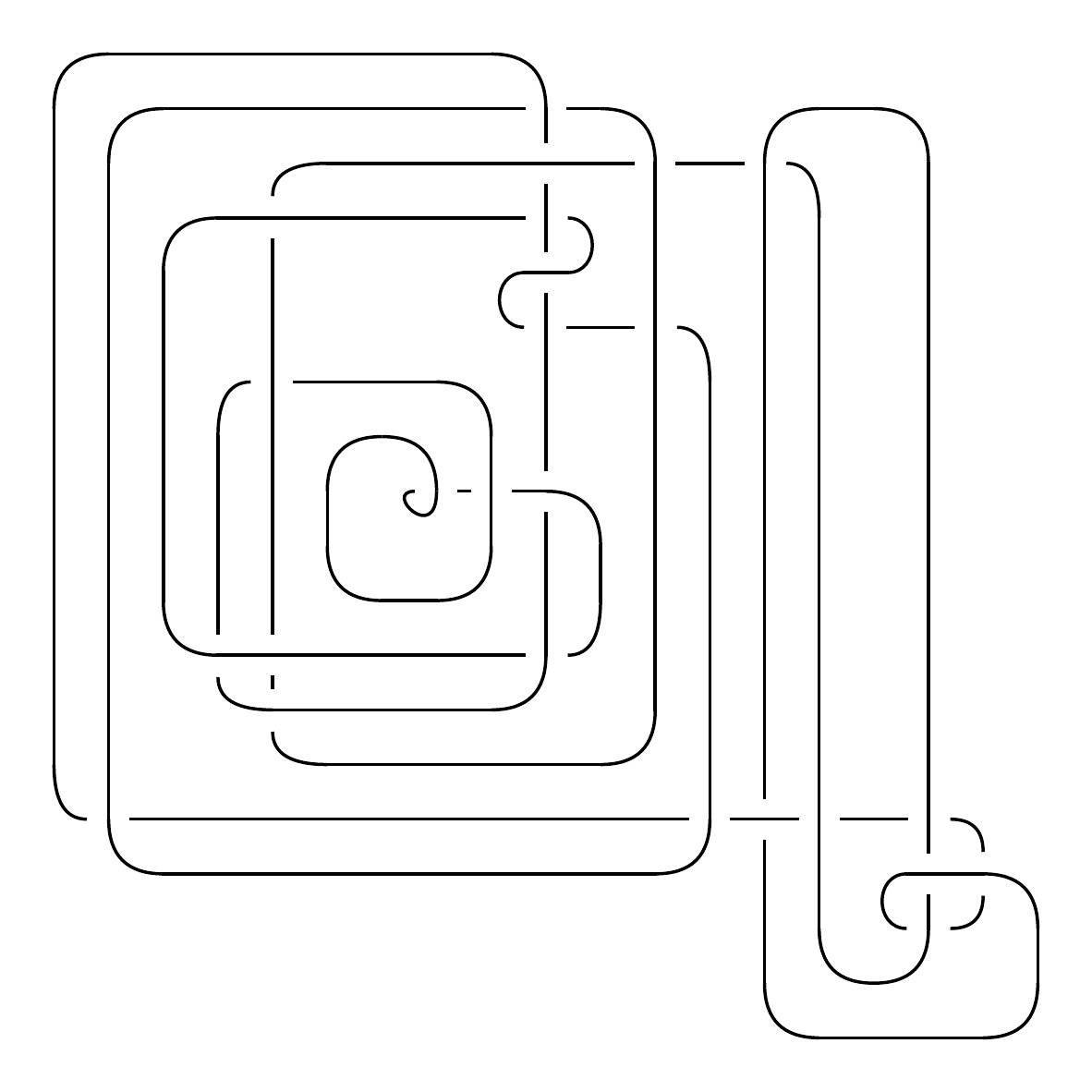}
    } 
    \makebox[\textwidth][c]{
        \includegraphics[scale=.4]{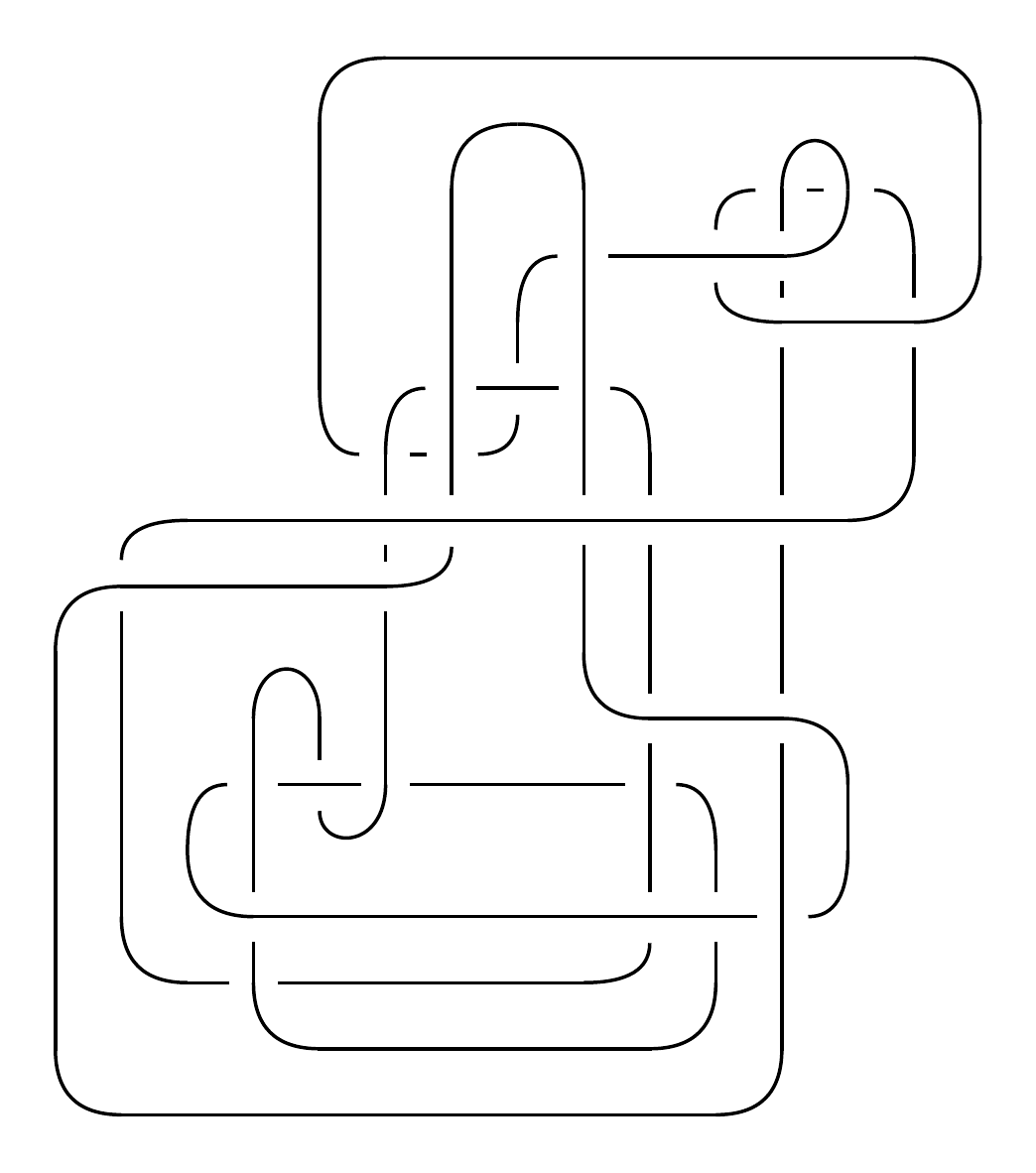}
        \includegraphics[scale=.4]{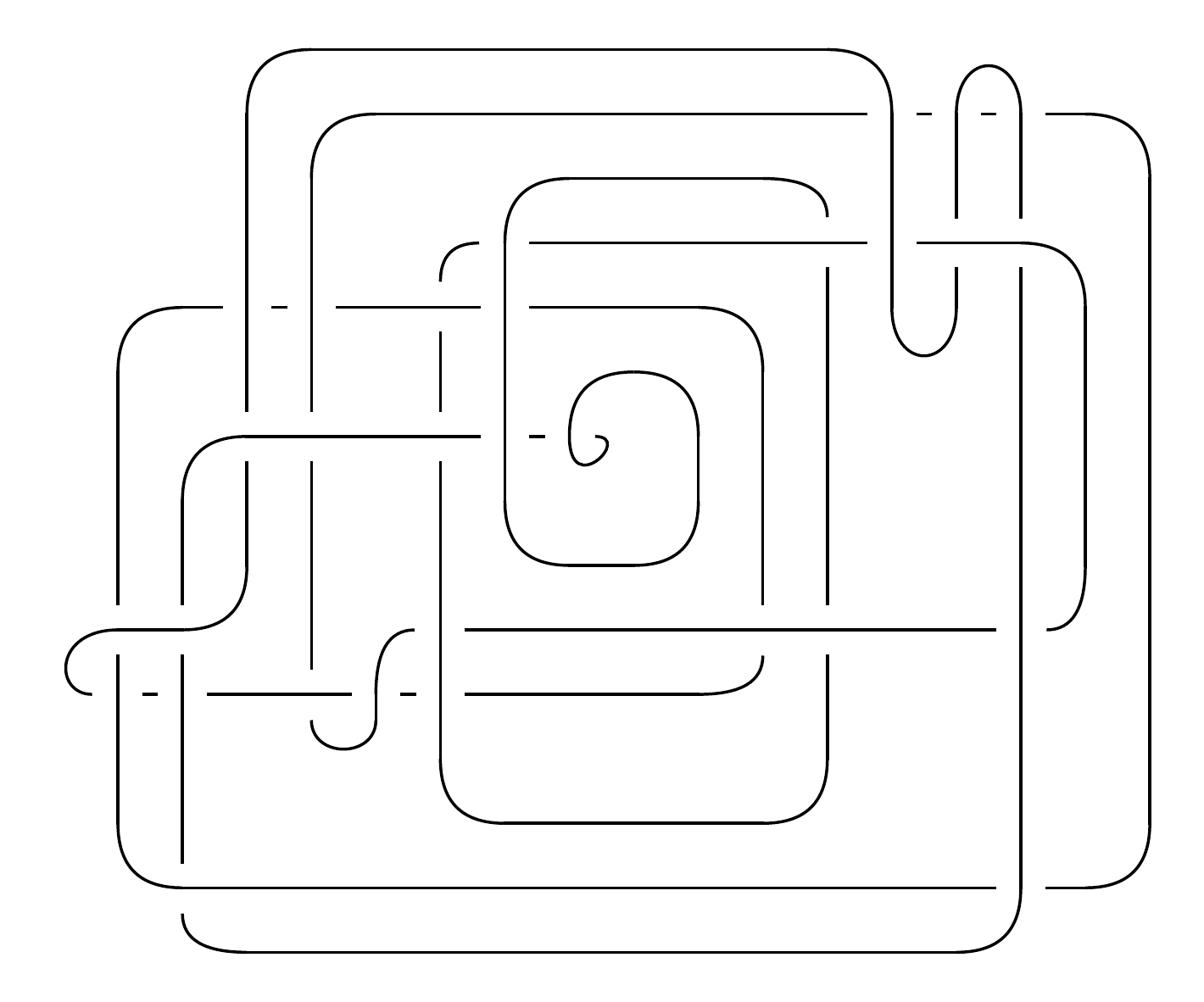}
        \includegraphics[scale=.4]{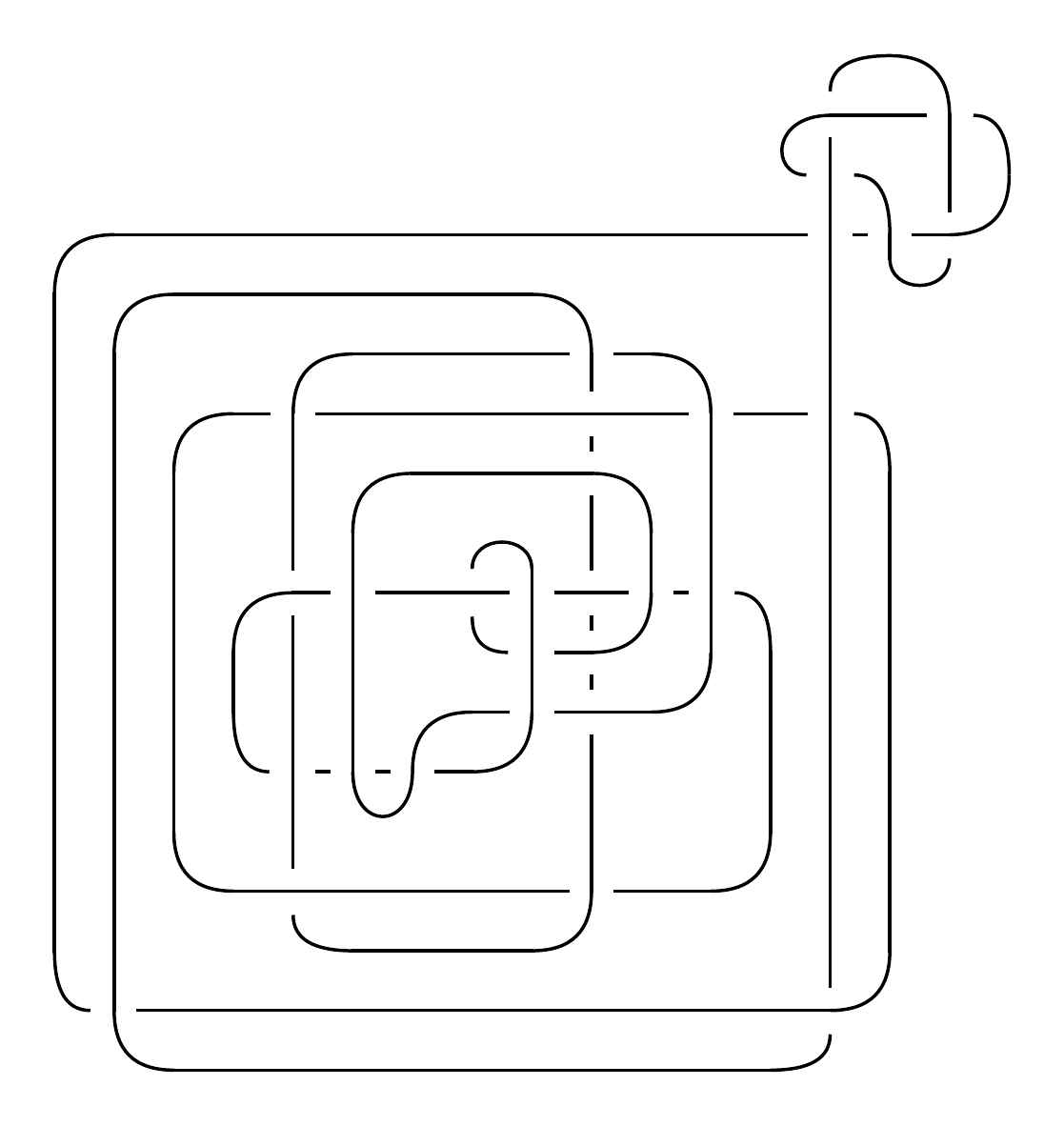}
    }
    \makebox[\textwidth][c]{
        \includegraphics[scale=.4]{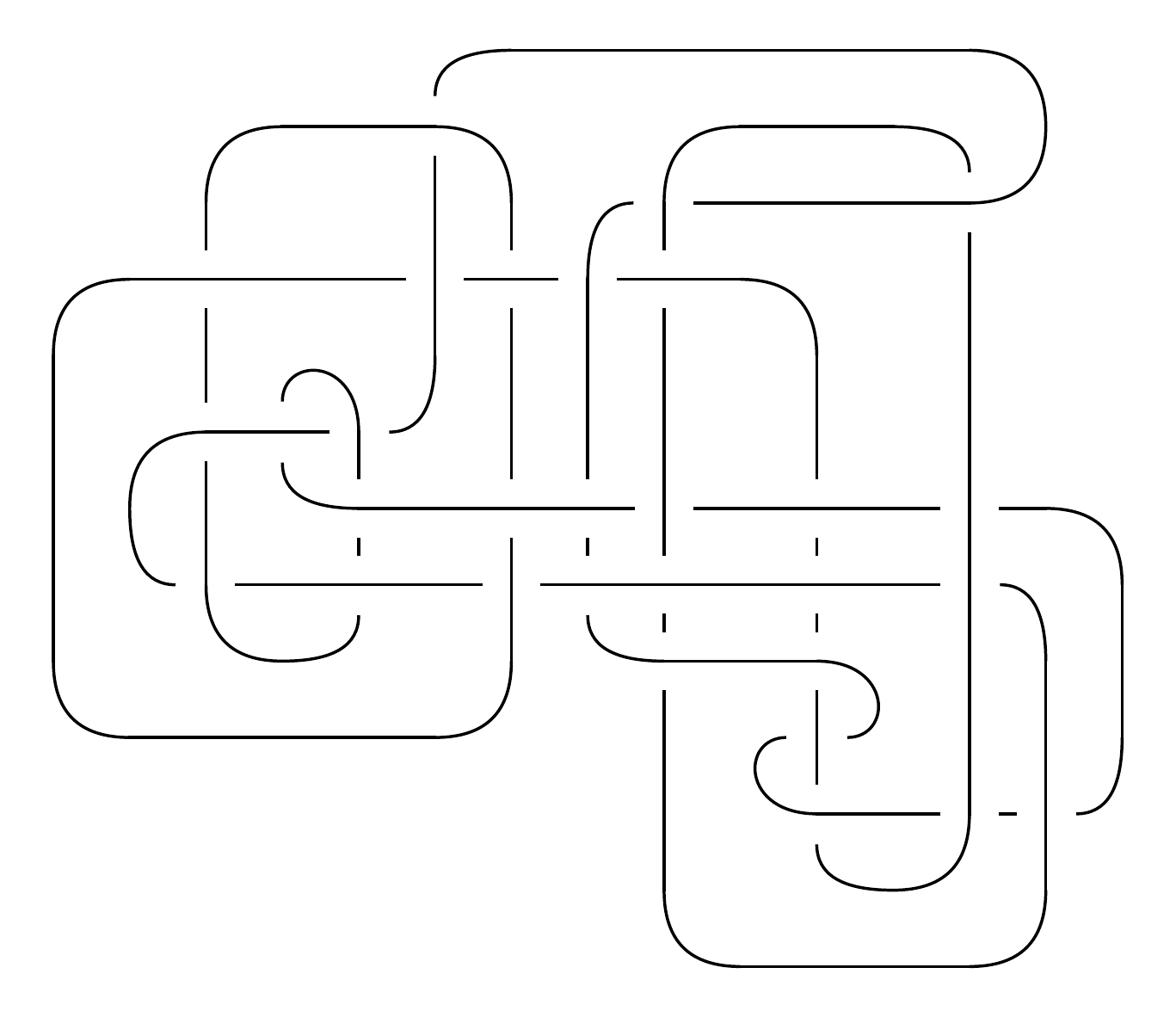}
        \includegraphics[scale=.4]{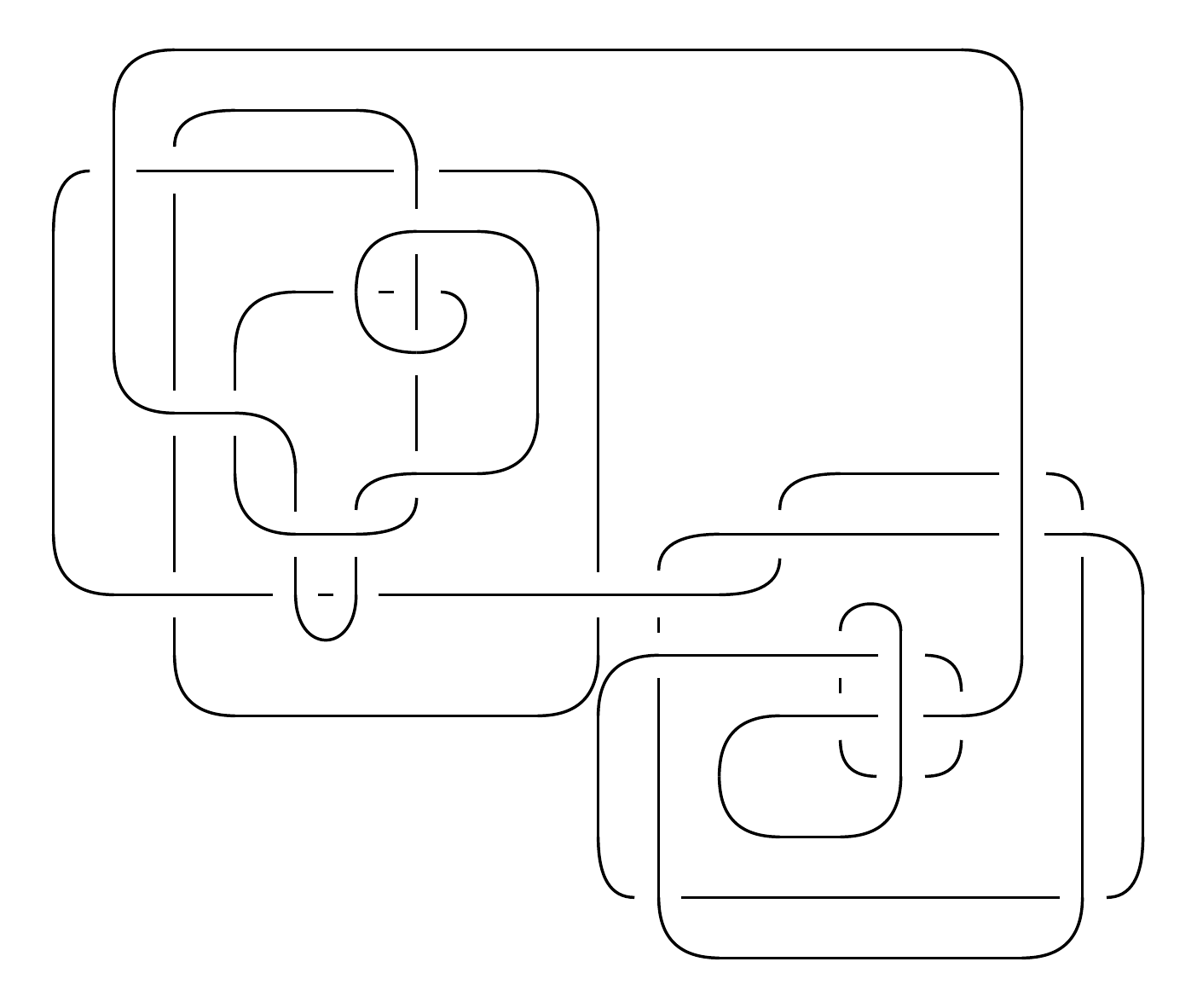}
        \includegraphics[scale=.4]{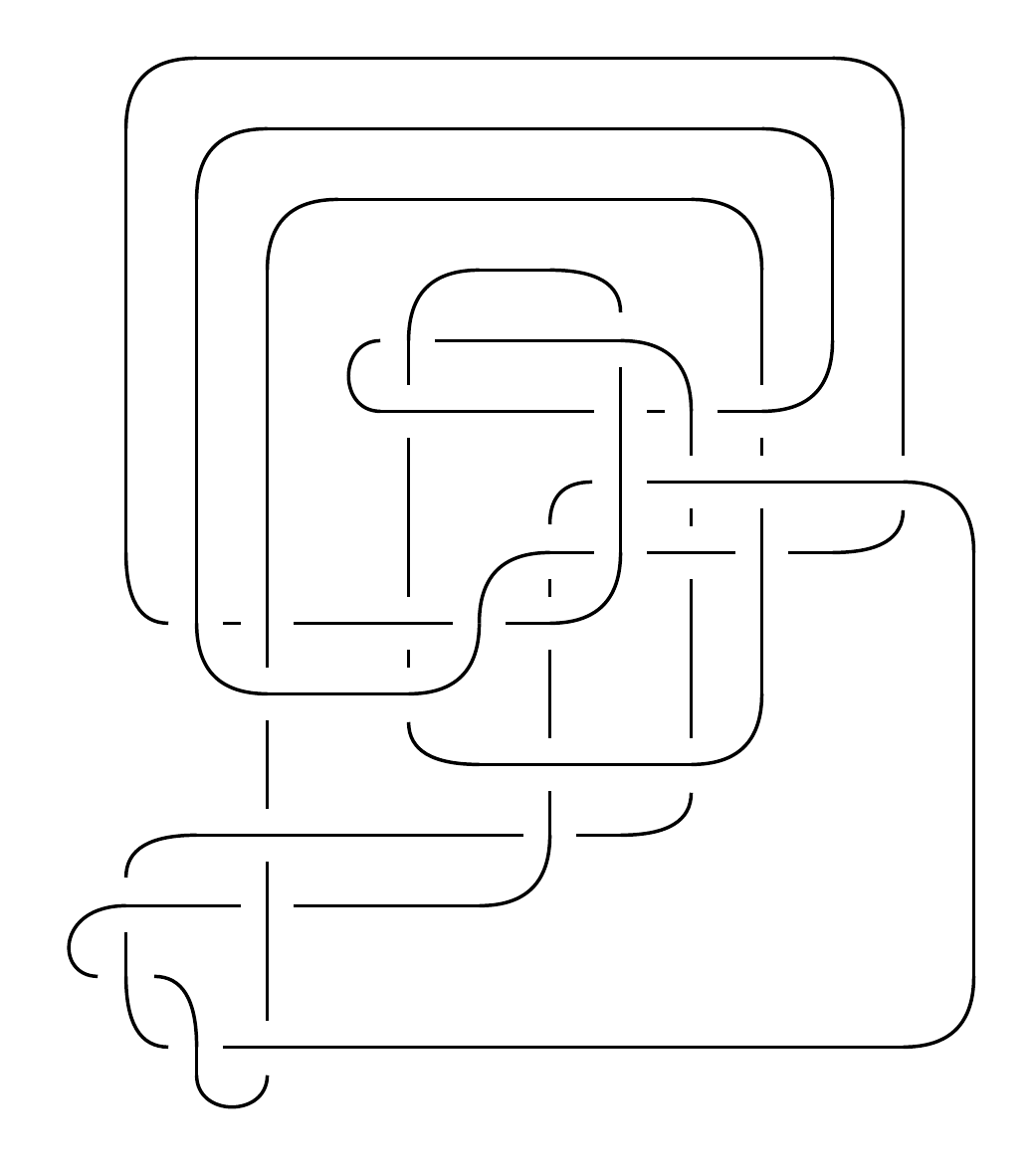}
    }
    \caption{Knot or not? Twenty-five and thirty crossing in rows $1$-$2$ and $3$-$4$, respectively.}
    \label{fig:knotornot3}
\end{figure}

To play, show your child Figures \ref{fig:knotornot1}-\ref{fig:knotornot3} and ask them to determine whether a given knot diagram
may be unknotted. The figures have increasing difficulty due to an 
increasing number of crossings. Solutions may be found in the footnote\footnote{Solutions are presented left-to-right, top-to-bottom, with K and U denoting non-trivial knots and unknots, respectively. Fig. \ref{fig:knotornot1}: KUUKUKUUKUKK. Fig. \ref{fig:knotornot2} UKUKKUKKUKUU. Fig. \ref{fig:knotornot3}
 KUUKUKUUKKUK.}.

\clearpage

\providecommand{\href}[2]{#2}\begingroup\raggedright\endgroup

\end{document}